\documentclass[12pt]{amsart}
\usepackage{amsmath}

\usepackage{url}

 \makeatletter
\renewcommand*\subjclass[2][2010]{%
  \def\@subjclass{#2}%
  \@ifundefined{subjclassname@#1}{%
    \ClassWarning{\@classname}{Unknown edition (#1) of Mathematics
      Subject Classification; using '2010'.}%
  }{%
    \@xp\let\@xp\subjclassname\csname subjclassname@#1\endcsname  
}%
}

 \makeatother

\begin{document}


\vspace*{1mm}

\title[Euclid's theorem on the infinitude  of primes
\ldots]{Euclid's theorem on the infinitude  of primes:
a historical survey of its  200 proofs (300 {\small B.C.}--2022)} 
\author{Romeo Me\v{s}trovi\'c}

\maketitle


\vspace*{3mm}

\hfill{\it ``The laws of nature are but the mathematical thoughts of God."}
\vspace*{3mm}

\hfill Euclid (circa 300 {\small B.C.})

\vspace*{2mm}

\hfill{\it ``If Euclid failed to kindle your youthful enthusiasm,}

\hfill{\it then you were not born to be a scientific thinker."}

\vspace*{3mm}

\hfill Albert Einstein

{\renewcommand{\thefootnote}{}\footnote{2010 {\it Mathematics Subject 
Classification.} Primary 01A05, 11A41, 11A51; Secondary 
11A07, 11A05, 11A25, 11B39, 11B50, 11B65, 11M06, 11N13.

{\it Keywords and phrases}: 
prime, infinitude of primes ($IP$), Elements, Euclid's (second) 
theorem,  Euclid's proof, Fermat numbers, 
Goldbach's proof of $IP$, proof of 
$IP$ based on algebraic number theory arguments, Euler's proof of $IP$,
combinatorial proof of $IP$, Furstenberg's proof of $IP$,
algorithmic proof of $IP$, proof of $IP$ in arithmetic progression,
Dirichlet's theorem, Euclidean proof.
\setcounter{footnote}{0}}

\begin{abstract} 
 
{\bf In the fourth  extended version of this  article, we provide a comprehensive 
historical survey  of 200 different proofs of  famous Euclid's  theorem on the infinitude of  prime 
numbers (300 {\small B.C.}--2022)}. 
 The author is trying
to collect almost all the known  proofs on infinitude 
of primes, including some proofs that can be easily obtained 
as consequences of some known problems or divisibility properties.
Furthermore, here are listed numerous elementary proofs  
of the infinitude of  primes in different arithmetic progressions.

All the references concerning the  proofs of  Euclid's theorem  that use 
similar   methods and ideas are exposed subsequently. Namely, presented proofs 
are divided into the first  five  subsections of Section 2 in dependence of  the methods 
that are used in them. 
  14 proofs which are  proved  from 2012 to 2017 are given in 
 Subsection 2.9, and 18  recent  proofs from 2018 to 2022
 are presented  in Subsection 2.10.

   In Section 3, we mainly  survey  
elementary proofs of the infinitude of primes in different arithmetic 
progressions. Presented proofs are special cases of Dirichlet's theorem.
In Section 4, we give a new simple ``Euclidean's proof" 
of the infinitude of primes. 
   \end{abstract}

\vfill\eject

\centerline{\large{\bf Contents}}

\vspace{5mm}

{\bf 1. Euclid's theorem on the infinitude of primes\hfill \bf 3}

1.1. Primes and the infinitude of primes\dotfill 3

1.2. Euclid's proof of Euclid's theorem\dotfill 4

1.3. Sequences arising from Euclid's proof of $IP$\dotfill 6

1.4. Proofs of Euclid's theorem: a brief history\dotfill 7

\vspace{3mm}

{\bf 2. A survey of different proofs of Euclid's theorem\hfill 9}

2.1. Proofs of $IP$ based on Euclid's idea\dotfill 9

2.2. Proofs of $IP$ based on Goldbach's idea on 
mutually prime 

\qquad integers\dotfill 11

2.3. Proofs of $IP$ based on algebraic number theory arguments\dotfill 15

2.4. Proof of $IP$ based on Euler's idea on the divergence of the sum 

\qquad of prime reciprocals and Euler's formula\dotfill 18

2.5. Proof of $IP$ based on Euler's product for the  Riemann 

\qquad zeta function   and the irrationality of $\pi^2$  
and $e$\dotfill 22

2.6. Combinatorial proofs of $IP$ based on enumerative arguments\dotfill 23

2.7. Furstenberg's topological proof of $IP$ and 
its modifications\dotfill 25

2.8. Another proofs of $IP$\dotfill 26

2.9. 14 proofs of $IP$ (2012--2017) \dotfill 29

2.10. 18 recent proofs of $IP$ (2018--2022) \dotfill 31

\vspace{3mm}

{\bf 3. Proofs of $IP$ in arithmetic progressions: special cases of 

\quad\, Dirichlet's theorem\hfill 34}

3.1. Dirichlet's  theorem\dotfill 34

3.2. A survey of elementary proofs of $IP$ in  special arithmetic

\quad\,\,\,\,  progressions\dotfill 35 

3.3. Elementary proofs of $IP$ 
in arithmetic progressions with small 

\quad\,\,\,\, differences\dotfill 39 

\vspace{3mm}

{\bf 4. Another simple Euclidean's proof of Euclid's theorem\hfill 40}

\vspace{3mm}

{\bf References \hfill 41}

\vspace{3mm}

{\bf Appendixes \hfill 57}

A) External Links on Euclid's theorem and its proofs\dotfill 57

B) Sloane's sequences related to proofs of Euclid's theorem\dotfill 57

C) List of papers and their authors arranged by year of publication 

\quad\, followed by the  main argument(s) of related proof given into 

\quad\, round  brackets\dotfill 58

D)  Author Index\dotfill 66 

E)  Subject Index\dotfill 70

  \section{Euclid's theorem on the infinitude of primes}

\subsection{Primes and the infinitude of primes}

 {\it prime number} (or briefly  in the sequel, a {\it prime}) is an 
integer greater than 1 that is divisible  only by 1 and itself.
Starting from the beginning, prime numbers have always 
been around but the concepts and uniqueness was thought to be first 
considered during Egyptian times. However, 
 mathematicians have been studying primes and their properties for over 
twenty-three centuries. Ancient Greek mathematicians knew that there are 
infinitely many primes. 
Namely, circa 300 {\small B.C.},  Euclid of Alexandria, from the
{\it Pythagorean School} proved ({\it Elements}, Book IX, Proposition 20)
the following celebrated result as rendered into modern language from the 
Greek (\cite{he}, \cite{z}): 

{\it If a number be the least that is 
measured by prime numbers, it will not be measured 
by any other prime number except those originally measuring it.}

Euclid's ``{\it Elements}"  
are one of the most popular and most widely printed mathematicians 
books and they are been translated into many languages. 
{\it Elements} presents   a remarkable collection of 13 books that 
contained much of the mathematical known at the time. Books VII, 
VIII and IX deal with properties of the integers and contain the early 
beginnings of number theory, a body of knowledge that has flourished ever 
since.

Recall that during Euclid's time, integers were understood as lengths of line 
segments and divisibility was spoken of as measuring.
According to G. H. Hardy \cite{ha}, ``{\it Euclid's 
theorem which   states that the number of primes is infinite
is vital for the whole structure of arithmetic. The primes are 
the raw material out of which we have to build arithmetic,
and Euclid's theorem assures us that we have 
plenty of material for the task}." 
Hardy \cite{ha} also remarks  that this proof 
is ``{\it as fresh and significant as when it 
was discovered--two thousand years have not written 
a wrinkle on it}". A. Weil \cite{wei2} also called ``{\it the proof 
for the existence of infinitely many primes represents ubdoubtedly
a major advance, but there is no compelling reason either for attributing
it to Euclid or for dating back to earlier times. What matters for our
purposes is that the very broad diffusion of Euclid in latter centuries,
while driving out all earlier texts, made them widely available to 
mathematicians from then on}". 

Sir Michael Atyah remarked during an interview \cite{rs}: {\it Any 
good theorem should have several proofs, more the better. For two reasons:
usually, different proofs have different strenghts and weaknesses, and they 
generalize in different directions - they are not just repetitions 
of each other.} For example, the {\it Pythagorean theorem} has received more 
than 360 proofs \cite{lo} of all sorts as algebraic, geometric, dynamic and 
so on. The {\it irrationality of} $\sqrt{2}$ is another famous example  of 
a theorem which has been proved in many ways 
(\cite{ti}; on the web page \cite{bog} fourteen different 
proofs appear). C. F. Gauss himself had  10 different proofs for the 
{\it law  of quadratic reciprocity} \cite[Sections 112--114]{ga}.
Surprisingly, here we  present  183 different proofs of {\it Euclid's 
theorem on the  infinitude of primes}, including 44 proofs 
of the infinitude of primes in special arithmetic progressions.

\subsection{Euclid's proof of Euclid's theorem}
Even after almost two and a half 
millennia ago Euclid's theorem on the infinitude of primes
 stands as an excellent model of reasoning. 
Below we follow Ribenboim's statement of Euclid's proof \cite[p. 3]{r2}.
Namely, in Book IX of his celebrated
{\it Elements} (see \cite{he}) we find Proposition 20, which states:
\vspace{1mm}

\noindent{\bf Euclid's theorem.} {\it There are infinitely many primes}. 
\vspace{2mm}

Elegant proof of Euclid's theorem runs as follows. {\it Suppose that 
$p_1=2<p_2=3<\cdots <p_k$ are all the primes. Take
$n=p_1p_2\cdots p_k+1$ and let $p$ be a prime dividing $n$. 
Then $p$ cannot be any of $p_1,p_2,\ldots,p_k$, otherwise $p$ would
divide the difference $n-p_1p_2\cdots p_k=1$.}\hfill$\Box$
\vspace{2mm}

The above proof is actually quite a bit different from what Euclid wrote. 
Since ancient Greeks did not have our modern notion 
of infinity, Euclid could not have written ``there are infinitely 
many primes", rather he wrote: ``{\it prime numbers are more than any assigned 
multitude of prime numbers}."  Below is a proof closer to that 
which Euclid wrote, but still using our modern concepts of 
numbers and proof. An English translation of 
Euclid's actual proof given by D. Joyce in his webpages 
\cite{j} also can be found in 

\noindent{\small{\tt http://primes.utm.edu/notes/proofs/infinite/euclids.html}}.
It is a most elegant proof by {\it contradiction}
({\it reduction ad absurdum}) that goes 
as follows. 

\vspace{2mm}

\noindent{\bf Euclid's theorem.} {\it There are 
more primes than found in any finite list of primes}.
 \begin{proof}
{\it Call the primes in our finite list $p_1,p_2,\ldots,p_k$. Let $P$ be 
any common multiple of these primes plus one 
$($for example $P=p_1p_2\cdots p_k+1$$)$. Now $P$ is either prime 
or it is not. If it is prime, then $P$ is a prime that was not in our
list. If $P$ is not prime, then it is 
divisible by some prime, call it $p$. Notice $p$ cannot be any  of
$p_1,p_2,\ldots ,p_k$, otherwise $p$ would divide $1$, which is impossible.
So this prime $p$ is some prime that was not in our original list.
Either way, the original list was incomplete.}
\end{proof}

The statement of Euclid's theorem together with its proof is given by 
B. Mazur in 2005 \cite[p. 230, Section 3]{maz}  as follows. 

``{\it If you give me 
any finite $($non-empty, of course$!)$ collection of prime numbers, I will
form the number $N$ that is $1$ more than the product of all the primes in the 
collection, so that every prime in your collection has the property that when 
$N$ is divided by it, there is a remainder of $1$. There exists at least one  
prime number dividing this number $N$ and any prime number dividing $N$ is 
new in the sense that it is not in your initial collection.}"

\vspace{2mm}

\noindent{\bf Remarks.} Euclid's proof is often said 
to be ``indirect'' or ``by contradiction'',
but this is unwarranted: given any finite 
set of primes $p_1,\ldots,p_n$, it 
gives a perfectly definite procedure for constructing 
a new prime. Indeed, if we define $E_1=2$, and having defined 
$E_1,\ldots ,E_n$, we define $E_{n+1}$ to be the smallest prime divisor 
of $E_1E_2\cdots E_n+1$, we get a sequence of distinct primes, nowadays 
called the {\it Euclid-Mullin sequence} (of course, we could get a different 
sequence by taking $p_1$ to be a prime different from 2). 
This is Sloane's sequence A000945  whose first few terms are 
 $2,3,7,43,13,53,5,6221671,38709183810571,139,\ldots$.    
The natural question - does every prime occur eventually in the Euclid-Mullin 
sequence remains unanswered. Note that D. Shanks \cite{sha2} 
conjectured on  probabilistic grounds that this sequrncet contains every prime.
This conjecture was supported by computational results up to 
43rd term of yhe sequence $(E_n)$ given in 1993 by S. S. Wagstaff, Jr. \cite{wag}.  
For a discussion on this conjecture, see \cite[Section 2]{boo2},
where it was noticed that N. Kurokawa and T. Satoh \cite{ksa} have 
shown that an analogue of this conjecture for the 
{\it Euclidean domains} $\Bbb F_p[x]$ is false in general.
Notice that the sequence $(E_n)$ and several related sequences were studied
in \cite{gn}. 

Moreover, Mullin \cite{mu2} constructed the second sequence of primes, say
$(P_n)$ similarly as the above sequence $(E_n)$, except that 
we replace the words ``smallest prime divisor'' by ``largest  
prime divisor''. This is the sequence  A000946 in 
\cite{sl}. It was proved in 2013  by A. R. Booker 
\cite[Theorem 1]{boo2} that  the sequence $(P_n)$ omits infinitely many 
primes, confirming a conjecture of C. D. Cox and A. J. Van der Poorten 
\cite{cop}. Notice that in 2014  
P. Pollack and E. Trevi$\stackrel{\sim}{\rm n}$o \cite{pt} 
gave a completely elementary proof of this conjecture. 


Notice also that Euclid's proof actually 
uses the fact that there is a prime dividing given positive 
integer greater than 1. This follows from 
Proposition 31 in Book VII  of his {\it Elements}
(\cite{he}, \cite{ao}, \cite[p.2, Theorem 1]{hw}) which asserts that 
``{\it any composite number is measured by some prime number}",
or in terms of  modern arithmetic, that 
every integer $n>1$ has at least one representation
as a product of primes. Of course,
he also used a unexpressed axiom which states that if $a$ divides $b$
and $a$ divides $c$, $a$ will divide the difference between $b$ and $c$.
\hfill$\Box$

\vspace{2mm}

 {\it  The unique factorization theorem}, otherwise   
known as the ``fundamental theorem of arithmetic," states that any integer 
greater than 1 can, except for the order of the factors, be expressed as a 
product of primes in one and only one way. 
This theorem does not appear in Euclid's {\it Elements} 
(\cite{he}; also see \cite{ao}). However, as noticed in \cite[page 208]{ao}, 
in fact, the unique factorization theorem follows 
from Propositions 30-31 in Book VII (given in  Remarks of Section 4).
More generally, in 1976 W. Knorr \cite{kn} gave a reasonable discussion of 
the position of unique factorization in Euclid's theory of numbers.
Nevertheless, as noticed in \cite{ao}, Euclid played a significant role 
in the history of this theorem (specifically, this concerns 
to some propositions of Books VII and IX). However, the first explicit and 
clear statement and the proof 
of the unique factorization theorem seems to be in C. F. Gauss'
masterpiece {\it Disquisitiones Arithmeticae} 
\cite[Section II, Article 16]{ga}. His Article 16 is given as the following 
theorem:  {\it A composite number can be resolved into prime factors in only 
one way.} After Gauss, many mathematicians provided different 
proofs of this theorem in their work (these proofs are presented and  
classified in  \cite{af}).
In particular, the unique factorization theorem 
was used  in numerous proofs of the infinitude of primes provided below.

Notice also that for any field $F$, Euclid's argument works to show that 
there are infinitely many irreducible polynomials over $F$. This  follows 
inductively taking $p_1(t)=t$, and having produced $p_1(t),\ldots,p_k(t)$, 
consider the irreducible factors of $p_1(t)\cdots p_k(t)+1$.

\subsection{Sequences arising from Euclid's proof of $IP$} As usually,  
for each prime $p$, $p^{\#}$ denotes the product of all  the primes less than 
or equal to $p$ and it is  called the {\it primorial number}
(Sloane's sequence A002110; also see A034386 for the second 
definition of primorial number as a product of primes in the range 
$2$ to $n$). The expressions $p^{\#}+1$ and $p^{\#}-1$ have been 
considered in connection with variants of the Euclid's proof of the 
infinitude of primes.

Further, $n$th {\it Euclid's number} $E_n$
(see e.g., \cite{v})  is defined as a product of first $n$ 
{\it consecutive primes} plus one (Sloane's sequence A006862).   Similarly, 
{\it Kummer's number}  is defined as a product of first $n$ 
consecutive primes minus one (Sloane's sequence A057588).
Euclid's numbers were tested for primality in 1972 by A. Borning \cite{bor},
in 1980 by M. Templer \cite{te}, in 1982 by  J. P. Buhler, R. E. Crandall 
and M. A. Penk \cite{bcp},  and  in 1995 by C. K. Caldwell \cite{cal}. 
Recall also that two interesting conjectures involving the numbers $E_n$
are quite recently proposed by Z.-W. Sun. Namely, for any given $n$, 
if $w_1(n)$ is defined as the least integer $m>1$ such that 
$m$ divides none of those $E_i-E_j$ with $1\le i<j\le n$, then 
Sun \cite[Conjecture 1.5 (i) and (iii)]{sun} conjectured that $w_1(n)$ is a 
prime less than $n^2$ for all $n=2,3,4,\ldots$.
The same conjecture \cite[Conjecture 1.5 (ii) and (iii)]{sun} is proposed 
in relation to the sums $E_i+E_j-2$ instead of $E_i-E_j$
(cf. Sloane's sequences A210144 and A210186). 

The numbers $p^{\#}\pm 1$ (in accordance to the 
first definition given above) and $n!\pm 1$ 
have been frequently checked for {\it primality} 
(see \cite{cg}, \cite{gkp}, \cite{sp} and \cite[pp. 4--5]{r2}).
The numbers $p^{\#}\pm 1$ have  been tested for all $p<120000$
in 2002  by C. Caldwell and Y. Gallot \cite{cg}. They were 
reported that in the tested range there are exactly  19 primes of 
the form $p^{\#} +1$ and 18 primes of the form $p^{\#}-1$
(these are in fact Sloane's sequences A005234 extended with three
new terms and A006794, respectively).
It is pointed out in \cite[p. 4]{r2} 
that the answers to the following 
questions are  unknown: 1) Are there infinitely many primes $p$ for which
$p\# +1$ is  prime?). Are there infinitely many primes $p$ for which
$p\# +1$ is  composite? 

In terms of the second definition of primorial numbers given above,
similarly are defined  Sloane's sequences A014545 
and A057704 (they also called primorial primes).

Other Sloane's sequences related to Euclid's proof and Euclid numbers
are: A018239 (primorial primes),  A057705,  A057713, A065314, A065315,
A065316, A065317, A006794, 
A068488, A068489, A103514, A066266, A066267, A066268, A066269,
A088054, A093804, , A103319, A104350,  A002981, A002982, A038507, A088332,
A005235, A000945 and A000946.\hfill$\Box$ 
\vspace{2mm}

\subsection{Proofs of Euclid's theorem: a brief history}
Euclid's  theorem on the infinitude of primes 
has fascinated generations of mathematicians since its 
first and famous demonstration given by Euclid (300 {\small B.C.}).
Many great mathematicians of the eighteenth  and  nineteenth century
established different proofs of this theorem 
(for instance, Goldbach (1730), Euler (1736, 1737), Lebesgue (1843, 1856, 
1859, 1862), Sylvester (1871, 1888 (4)), 
Kronecker (1875/6), Hensel (1875/6), Lucas (1878, 1891, 1899), 
Kummer (1878/9),   Stieltjes (1890) and Hermite (1915). 
Furthermore, in the last hundred years  
various interesting proofs of the infinitude of primes, including  
the infinitude of primes in different arithmetic progressions,
were obtained by I. Schur  (1912/13), K. Hensel (1913), G. P\'{o}lya (1921),
G. P\'{o}lya and G. Szeg\H{o} (1925), P. Erd\H{o}s (1934 (2), 1938 (2)), 
 G. H. Hardy and E. M. Wright (1938 (2)), L. G.  Schnirelman   
(published posthumously in 1940), R. Bellman (1943, 1947), 
H. Furstenberg (1955), J. Lambek and L. Moser (1957), S. W. Golomb (1963),
A. W. F. Edwards (1964), A. A. Mullin (1964),
 W. Sierpi\'{n}ski (1964, 1970 (4)), S. P. Mohanty (1978 (3)),
A. Weil (1979), L. Washington (1980), S. Srinivasan (1984 (2)), 
M. Deaconescu and J. S\"{a}ndor (1986), J. B. Paris, A. J. Wilkie and 
A. R. Woods (1988), M. Rubinstein (1993), N. Robbins 
(1994 (2)), R. Goldblatt (1998), M. Aigner and G. M. Ziegler (2001 (2)), 
\v{S}. Porubsky (2001), D. Cass and G. Wildenberg (2003),
T.  Ishikawa, N. Ishida and Y. Yukimoto (2004),
R. Crandall and C. Pomerance (2005), 
A. Granville (2007 (2), 2009), J. P. Whang (2010), 
R. Cooke (2011), P. Pollack (2011) and by several other authors.
We also point out that in numerous proofs of Euclid's theorem 
were used methods and arguments  due to Euclid  
(``Euclidean's proofs"), Goldbach 
(proofs based on elementary  divisibility properties of integers) or 
Euler (analytic proofs based on Euler's product). 
Moreover, numerous proofs of Euclid's theorem 
are  based on some of the following methods or 
results: algebraic number theory arguments 
(Euler's totient function, Euler theorem, 
Fermat little theorem, arithmetic functions,
Theory of Finite Abelian Groups etc.), 
Euler's formula for the Riemann zeta function, Euler's factorization,
elementary counting methods (enumerative arguments),  
Furstenberg's topological proof of the  infinitude of primes 
and its combinatorial and 
algebraic modifications etc. All the
proofs of the infinitude of primes exposed  in this articale 
are divided into 8 subsections of Section 2 in dependence of used methods 
in them. In the next section we mainly survey  elementary proofs 
of the  infinitude of primes in different arithmetic progressions. 
These proofs are also based on some of mentioned methods and ideas. Finally, 
in Section 4, we give a new simple proof of the infinitude of primes. 
The first step of our proof is based on Euclid's idea. The remaining of the 
proof is quite simple and elementary  and it does not use the notion of 
divisibility.

In  Dickson's  {\it History of the Theory of Numbers}   
\cite[pp. 413--415]{d} 
and the books by Ribenboim \cite[pp. 3--11]{r}, 
\cite[Chapter 1, pp. 3--13]{r2}, Pollack \cite[pp. 2--19]{pol}, 
Hardy and Wright \cite[pp. 12--17]{hw}, \cite[pp. 14--18]{hw2}, 
Aigner and Ziegler  \cite[pp. 3--6]{az},
and in Narkiewicz's monograph \cite[pp. 1--10]{na}
can be found many different proofs of  Euclid's theorem. Several 
proofs of this theorem  were also explored by  
P. L. Clark  \cite[Ch. 10, pp. 115--121]{cl} and 
T. Yamada \cite[Sections 1-6, 10-12]{y}. 
In Appendix C) of this article we give  a list of all
168    different proofs of Euclid's theorem presented here 
(including elementary proofs related to the infinitude of primes 
in special arithmetic progressions), together with the corresponding 
reference(s), the name(s) of his (their) author(s) and the main method(s)  
and/or idea(s) used in it (them). 
This list is arranged by year of publication. 
We also give a comprehensive (Subject and Author) Index  
to this article. 
 
The Bibliography of this article contains  291 references, 
 consisting mainly of  articles 
(including 47  {\it Notes} and {\it Aricles} published in 
{\it Amer. Math. Monthly}) and mathematical textbooks and monographs.
It also includes a few unpublished works or problems that are available on  
Internet Websites, 
especially on {\tt http:arxiv.org/}, one  Ph.D. thesis, an interview,
one private correspondence, one Course Notes 
 and  {\it Sloane's On-Line Encyclopedia of Integer Sequences}. 
Some of these references does not concern  directly to  
proofs of the infinitude of primes, but results of each of them that are 
cited here give possibilities to simplfy  some of these known proofs.  

We believe that our exposition of different proofs of 
Euclid's theorem may be useful for establishing proofs of
many new and old results in Number Theory via elementary methods.

 \section{A survey of different proofs of Euclid's theorem}

To save the space, in the sequel we will often denote 
by ``$IP$"  ``the infinitude of primes".

\subsection{Proofs of $IP$ based on Euclid's idea}

Ever since  Euclid of Alexandria, sometimes before
300 {\small B.C.},  first  proved that the number of primes is infinite
(see  Proposition 20 in Book IX  of his legendary {\it Elements}
 in \cite{he} (also see  \cite[p. 4, Theorem 4]{hw}) where this result
is called {\it Euclid's second theorem}), 
mathematicians have amused themselves by coming up 
with alternate proofs. 
For more information about the Euclid's proof
of the infinitude of primes see e.g.,   
\cite{co}, \cite[p. 414, Ch. XVIII]{d}, \cite{di}, 
\cite[pp. 73--75]{du}, \cite{haw} and \cite[Section 3]{maz}. 

Euclid's proof of $IP$ is a paragon of simplicity: given a finite list of 
primes, multiply them together and add one. The resulting number, say $N$,
is not divisible by any prime on the list, so any prime factor of $N$ is a
new prime. There are several variants of Euclid's proof 
of $IP$. The simplest of them, which according to  H. Brocard 
\cite{bro}  is due in 1915 to C. Hermite, immediately follows from the 
obvious fact that the smallest prime divisor of $n!+1$ is  greater 
than $n$. Another of these proof, due to E. E. Kummer in 1878/9 \cite{k}
(also see \cite[page 4]{r2} and \cite{y})
is in fact an elegant variant of Euclid's proof. 
In a long paper published in two installments 120 years ago (\cite{pe2},
\cite{pe3}) J. Perott noticed that Euclid's proof works  if we consider 
$p_1p_2\cdots p_k-1$ instead of $p_1p_2\cdots p_k+1$.
 Stieltjes' proof in 1890 given in his work \cite[p. 14]{st}
(also see \cite[p. 414]{d}, \cite{r}, \cite[p. 4]{na}), 
 C. O. Boije af Genn\"{a}s' proof in 1893 
 \cite{bg} (also see \cite[p. 414]{d}, \cite{y}),
Braun's proof in 1899 (\cite{bra}; also see 
\cite[p. 414]{d}, \cite[p. 3]{pol} an \cite[p. 5]{na}), 
L\'{e}vi's proof in 1909/10 (\cite{lev}; see also \cite[p. 414]{d}),
M\'etrod's proof in 1917 (\cite{met}; see also 
\cite[p. 415]{d} and \cite[page 11]{r2}), Thompson's proof 
in 1953 \cite{tho}, Mullin's proof of 1964 \cite{mu}, Trigg's proof 
in 1974 \cite{t} and Aldaz and Bravo's proof (\cite{ab},
\cite[p. 6, Exercise 1.2.6]{pol}) 
in 2003 present refinements of Euclid's proof on $IP$.
For example, supposing that the set of all primes 
 is a finite $\{p_1,p_2,\ldots,p_k\}$ with their product $P$,
then setting $\sum_{i=1}^{k}1/p_i=a/P$ with 
$a=\sum_{i=1}^{k}P/p_i$, we find that $a/P>1/2+1/3+1/5=31/30>1$.
Therefore, Braun \cite{bra} concluded that 
$a$ must have a prime divisor, say $p_j$, but then
$p_j$ must divide $P/p_j$, which is not possible.  

Using algebraic number theory, in 1985 R. W. K. Odoni \cite{od} 
investigated the sequence $(w_n)$ recursively defined by R. K. Guy and 
R. Nowakowski \cite{gn} as $w_1=2$, $w_{n+1}=1+w_1\cdots w_n$ ($n\ge 1$) and 
observed that $w_n\to\infty$ as $n\to\infty$ and the $w_n$ are pairwise 
relatively prime. Clearly, this yields $IP$.

Furthermore, Problem  62 of \cite[pages 5, 42 and 43]{si2},
whose solution uses Euclid's idea,
 asserts that if $a,b$ and $m$ are positive integers 
such that $a$ and $b$ are relatively prime, then the 
arithmetic progression $\{ak+b:\,k=0,1,2,\ldots\}$ 
 contains  infinitely many terms relatively prime to $m$. This together with 
Euclid's argument (i.e., assuming $m$ to be  a product of consecutive 
primes) immediately yields $IP$. 
A proof of $IP$ quite similar to those of Braun is given in 2008 by A. Scimone 
\cite{sci}. Namely,  if $p_1,p_2,\ldots,p_k$ are all the primes with 
a product $N$, then Scimone consider  the divisors of the sum 
$\sum_{i=1}^k{N/p_i}$ to obtain  an immediate contradiction. 
Applying the {\it Chinese remainder theorem},
A. Granville considered more general sum 
in his Course Notes of 2007 \cite[Exercise 1.1b]{gr2} to prove  $IP$. 

If $p_n$ denotes the $n$th prime,  
 then by  \cite[pages 37 and 38, Problem 47; pages 8 and 55, 
Problem 92]{si2}  solved by A. M\c{a}kowski,
$p_{n+1}+p_{n+2}\le p_1p_2\cdots p_n$ for each $n\ge 3$.
This shows that for each $n\ge 3$ 
there are at least two primes between the  $n$th prime  and
the product of the first $n$ primes.
This estimate is in 1998 improved by J. S\'{a}ndor 
\cite{sa} who  showed that 
$p_{n}+p_{p_n-2}+p_1p_2\cdots p_{n-1}\le p_1p_2\cdots p_n$
 for all $n\ge 3$.

In 2008 B. Joyal \cite{jo} proved $IP$ using the 
{\it sieve of Eratosthenes}, devised about 200
{\small B.C.}, which is a 
beautiful and efficient algorithm for finding all the primes less than 
a given number $x$.

Recently, using Euclid's idea and a {\it representation of a rational 
number in a positive integer base}, in \cite{me} the author of this article 
obtained an elementary proof of $IP$. The second similar author's proof of 
$IP$  is given here in Section 4.

We see from Euclid's proof that  $p_{n+1}<p_1p_2\cdots p_n$ for each 
$n\ge 2$, where $p_k$ is the $k$th prime.
In 1907 H. Bonse \cite{bo} gave an elementary proof of a stronger inequality, 
now called {\it Bonse's inequality} \cite[p. 87]{uh}: if $n\ge 4$,
then $p_{n+1}^2<p_1p_2\cdots p_n$. In 2000 M. Dalezman  \cite[Theorem 1]{dal}
gave an elementary proof of stronger inequality 
$p_{n+1}p_{n+2}<p_1p_2\cdots p_n$ with $n\ge 4$.
 J. Sondow \cite[Theorem 1]{so1} exposed a simple proof 
based on the {\it Euler formula} $\zeta(2):=\sum_{n=1}^{\infty}1/n^2=\pi^2/6$
(suggested by P. Ribenboim in 2005), 
that for all sufficiently large $n$,
 $p_{n+1}<(p_1p_2\cdots p_n)^{2\mu}$, where $\mu$ is 
the {\it irrationality measure} for $6/\pi^2$ 
(for this concept and related estimates  see e.g., \cite[pp. 298--309]{r3}).
 Recall also that Bonse's inequality is refined in 1960 by L. P\'{o}sa 
\cite{pos}, in 1962 by S. E.  Mamangakis 
\cite{ma}, in 1971 by S. Reich  \cite{re} and in 
1988 by J. S\'{a}ndor \cite{sa}. 

\vspace{2mm}

\noindent{\bf Remarks.} Euclid's proof of $IP$  may be used to generate a 
sequence $(a_n)$ of primes as follows: put $a_1=2$ and if 
$a_1,a_2,\ldots,a_{n-1}$ are already defined then let $a_n$ be the largest 
prime divisor of $P_n:=a_1a_2\cdots a_{n-1}+1$ (Sloane's sequence A002585). 
This sequence was considered by A. A. Mullin in 1963 \cite{mu2} who asked 
whether it contains all primes and is monotonic. After a few terms of this 
sequence were computed (in 1964 by R. R. Korfhage \cite{ko}, in 1975 
by R. K. Guy and R. Nowakowski \cite{gn} and in 1984 by T. Naur \cite{nau}) 
it turned out that $a_{10}<a_9$. It is still unknown 
whether a sequence $(a_n)$ contains all sufficiently large primes.
Moreover, it can be  constructed the second sequence of primes, 
similarly as the above sequence $(a_n)$, except that 
we replace the expression ``$P_n:=a_1a_2\cdots a_{n-1}+1$'' by 
``$Q_n:=a_1a_2\cdots a_{n-1}-1$''. This is the sequence  A002584 in 
\cite{sl}. \hfill$\Box$

\subsection{Proofs of $IP$ based on Goldbach's idea on 
mutually prime integers} 

Goldbach's  idea consists in the obvious fact that 
any infinite sequence of pairwise relatively prime positive integers 
leads to a proof of Euclid's theorem.
C. Goldbach's proof presented in a letter 
to L. Euler in  July 20, 1730  
(see Fuss \cite[pp. 32--34, I]{fu}, \cite[p. 6]{r2},
\cite[pp. 40--41]{ew}, \cite[p. 4]{pol} or \cite[pp. 85--86]{and}) is based on 
the fact that the {\it Fermat numbers} $F_n:=2^{2^n}+1$, $n=0,1,2,\ldots$ are 
{\it mutually prime} (that is, {\it pairwise relatively prime}). 
{\it Indeed, it is easy to see by induction that 
$F_m-2=F_0F_1\cdots F_{m-1}$. This shows that if $n<m$, then $F_n$ 
divides $F_m-2$. Therefore, any prime dividing both $F_m$ and 
$F_n$ $(n<m)$ must divide the difference $2=F_m-(F_m-2)$.
But this is impossible since $F_n$ is odd, and 
this shows that Fermat numbers are pairwise 
relatively prime. Finally, assuming a prime factor
of each of integers  $F_n$, we obtain an
infinite sequence of distinct prime numbers}.

It seems that this was the first proof of $IP$ which 
essentially differed from that of Euclid.
In 1994 P. Ribenboim \cite{r4} wrote that the previous proof 
appears in an unpublished list of exercises of A. Hurwitz 
preserved in ETH in Z\"{u}rich. A quite similar 
 proof was published in the well known collections of exercises of 
G. P\'{o}lya and G. Szeg\H{o} \cite[p. 322, Problem 94]{ps} in  1925
(see also \cite[p. 14, Theorem 16]{hw}).

Clearly, Goldbach's idea is based on the fact that, in general 
the prime divisors of a sequence of integers greater than 1 form an infinite 
sequence of distinct primes if the integers in the sequence are pairwise 
relatively prime. In other words, Goldbach's proof of $IP$ will work 
with any sequence of positive integers for which any two distinct
terms of the sequence are relatively prime.  

Notice that Fermat numbers  $F_n$ are  Sloane's sequence A000215; 
other  sequences related to Fermat numbers are
  A019434, A094358, A050922, A023394 and A057755 and A080176. 
Today, the Fermat and Mersenne numbers $M_n:=2^n-1$ which are 
 considered in the next subsection, are important topics of discussion in 
many courses devoted to elementary number theory. For more information on 
classical and alternative approaches to the Fermat and Mersenne numbers see 
the article \cite{jr}.

In 1880 J. J. Sylvester (see e.g., \cite{v} 
and {\it Wikipedia}) generalized Fermat 
numbers via  a {\it recursively defined sequence of positive integers} 
in which every term of the sequence is the product of the previous
terms, plus one. This sequence   
is  called {\it Sylvester's sequence} and it is recursively defined as 
$a_{n+1}=a_n^2-a_n+1$ with $a_1=2$ (this is  Sloane's sequence A000058)
and generalized by Sloane's sequences A001543 and  A001544. 
Clearly, choosing a prime factor of each term of Sylvester's sequence 
yields  $IP$.
  
Goldbach's idea is later  used by many authors to prove Euclid's theorem 
by a construction of an infinite sequence of 
positive integers $1<a_1<a_2<a_3<\cdots$ that are 
pairwise relatively prime (i.e., without a common prime factor). 
In particular,  in 1956  V. C. Harris \cite{har}  
(see also \cite[p. 6, Exercise 1.2.5]{pol}, \cite{y})
inductively defined  an increasing sequence  of 
pairwise relatively prime positive integers (cf. Sloane's sequence A001685).
This is the sequence $(A_n)$ recursively defined as 
$A_n=A_0A_1\cdots A_{n-3}A_{n-1}+A_{n-2}$, for $n\ge 3$ 
($A_0,A_1$ and $A_2$ are given pairwise coprime positive integers, and 
$A_n$ is the numerator of approximants of some regular 
infinite continued fraction). 

Euclid's argument and Goldbach's idea are applied in solution of 
Problem 52 \cite[pages 5 and 40]{si2}
to show that there exist arbitrarily long 
arithmetic progressions formed of different positive 
integers such that every two terms of these 
progressions are relatively prime; namely, for any fixed integer $m\ge 1$ 
the numbers $(m!)k+1$ for $k=1,2,\ldots,m$ are relatively prime 
(cf. Sloane's sequence A104189). This yields $IP$. 
This proof was later communicated to P. Ribenboim by P. Schorn 
\cite[pp. 7--8]{r2}.

Several other sequences leading to proofs of $IP$ 
were established in 1957 by J. Lambek and L. Moser \cite{lm} and 
in 1966 by M. V. Subbarao \cite{su}.
 Furthermore, in 1964 A. W. F. Edwards (\cite{e}, \cite[page 7]{r2}) 
indicated various sequences, defined recursively, having this property 
(two related sequences are  Sloane's sequences A002715 and A002716). 
Similarly, in 2003 M. Somos and R. Haas \cite{sh}  proved $IP$ using  
an integer sequence defined recursively  
whose terms are pairwise relatively prime (cf. Sloane's sequences A064526,
A000324 and A007996).
All these sequences (excluding one defined by Harris) 
and several other  sequences of pairwise relatively 
prime positive integers are presented quite recently by A. Nowicki
in his monograph \cite[pp. 50--53, Section 3.5]{no}.
For example, if $f(x)=x^2-x+1$, then for any fixed $n\in\Bbb N$,
a sequence $n,f(n),f(f(n)),f(f(f(n))),\ldots$ has this property    
\cite[p. 51, Problem 3.5.4]{no}. This is also satisfied for the
following sequences $(a_n)$ defined recursively as: $a_{n+1}=a_n^3-a_n+1$;  
 $a_1=a,a_2=a_1+b,\ldots, a_{n+1}=a_1a_2\cdots a_n+b,\ldots$ with any fixed 
$a,b\in\Bbb N$; $a_1=b, a_{n+1}=a_n(a_n-a)+a$
with any fixed $a,b\in\Bbb N$ such that $b>a\ge 1$; $a_1=2,a_{n+1}=2^{a_n}-1$,
 and also for the sequence $a_n:=1+3^{3^n}+9^{3^n}$  
given in \cite[pp. 51--52, Problems 3.5.5, 3.5.6, 3.5.7, 3.5.10 and 3.5.15, 
respectively]{no}. Furthermore, by a problem of 1997
Romanian IMO Team Selection Test \cite[p. 149, Problem 7.2.3]{ana}), 
for any fixed integer $a>1$, the sequence $(a^{n+1}+a^n+1)$ ($n=1,2,\ldots$)
contains an infinite subsequence consisting 
of pairwise relatively prime positive integers. 
By a problem of the training of the German IMO team
\cite[pp. 121--122, Problem E3]{en}, using the factorization
$2^{2^{n+1}}+2^{2^{n}}+1=(2^{2^{n}}-2^{2^{n-1}}+1)(2^{2^{n}}+2^{2^{n-1}}+1)$,
it was proved  that 
$2^{2^{n+1}}+2^{2^{n}}+1$ has at least $n$ different prime factors
for each positive integer $n$. 

In 1965 M. Wunderlich \cite{wu} (also see \cite[p. 9, eleventh 
proof of Theorem 1.1]{na}) indicated that every sequence $(a_n)$
of distinct positive integers having the property that $(m,n)=1$ 
implies $(a_m,a_n)=1$ leads to the proof of $IP$ ($(m,n)$ denotes 
the {\it greatest common divisor} of $m$ and $n$). 
In particular, M. Wunderlich \cite{wu}
noticed that  {\it Fibonacci's sequence} 
$(f_n)$ (defined by conditions $f_1=f_2=1$, $f_{n+2}=f_{n+1}+f_n$ 
with $n=1,2,\ldots$; Sloane's sequence A000045) has this property
 (proved in 1846 by H. Siebeck \cite{sie}; also see \cite[p. 30]{vo}).
Notice that the sequence $(2^n-1)$  also satisfies  
this property because of 
the well known fact that $(2^n-1,2^m-1)=2^{(m,n)}-1$ for all $n,m\in\Bbb N$
(cf. \cite[5]{pol}). 
Using Wunderlich's argument indicated above,  in 1966
R. L. Hemminiger \cite{hem} established $IP$ by proving 
that the terms of the sequence $(a_n)$ defined recursively as 
$a_1=2$, $a_{n+1}=1+\prod_{i=1}^na_i$, are mutually prime.
However, it is easy by induction 
to show that  $a_{n+1}=a_n^2-a_n+1$ for each $n\in\Bbb N$
(cf. Granville's proof in \cite[p. 5, Exercise 1.2.3]{pol}), 
i.e., $(a_n)$ coincides with Sylvester's sequence.   
 
Further, $IP$ obviously follows from Problem 51 of 
\cite[pages 4 and 39]{si2} solved by A. Rotkiewicz which asserts that 
Fibonacci's sequence contains an infinite increasing subsequence such that 
every two terms of this sequence are relatively prime.
This means that the set of all prime divisors of   
Fibonacci sequence is infinite. It was shown in 1921 by G. P\'{o}lya 
\cite{poly} that the same happens for a large class of linear 
recurrences (also cf. related results of H. Hasse \cite{has2} in 1966, 
J. C. Lagarias \cite{la} in 1985, P. J. Stephens  \cite{st} in 1976, 
M. Ward (\cite{war1}) and \cite{war2}) in 1954 and 1961,  and 
H. R. Morton \cite{mor} in 1995).

Proof of $IP$ due to  S. P. Mohanty (\cite[Theorem 1 and Corollary 1]{mo1}; 
also see \cite{mo2}, \cite[pp. 5--6, Exercise 1.2.4]{pol}) 
in  1978, uses  sequences that generalized Sylvester's sequence. 
 By  a  problem of Polish Mathematical Olympiad in 2001/02 
(\cite[Problem 6]{pola}, see also \cite[p. 51, Problem 3.5.3]{no}), 
for any fixed positive integer $k$,  all the terms of a 
 sequence $(a_n)$ defined by $a_1=k+1$, $a_{n+1}=a_n^2-ka_n+k$,
are  pairwise relatively prime. Notice that this sequence is a
generalization of Sylvester's sequence and 
a particular case of a sequence from mentioned  Mohanty's proof.
Motivated by the same idea, in 1947 
R. Bellman \cite{b} (see also \cite[page 7]{r2}) gave a simple 
``{\it polynomial method}" to produce infinite sequences with the mentioned 
property. 
In 1978 S. P.  Mohanty \cite[Theorem 3]{mo1} proved that 
for any prime $p>5$, every prime divisor of Fibonacci number
$f_p$ is greater than $p$. This immediately yields $IP$. 
$IP$ also follows from  Problem 42 of 
\cite[pages 4, 35 and  36]{si2} which asserts that 
there exists an  increasing infinite sequence 
of pairwise relatively prime {\it triangular numbers} 
$t_n:=n(n+1)/2$, with $n=1,2,\ldots$ (Sloane's sequence A000217).
The same statement related 
to the {\it tetrahedral numbers} $T_n:=n(n+1)(n+2)/6$, with $n=1,2,\ldots$,
was given by  Problem 43 of \cite[pages 4 and 36]{si2}
(Sloane's sequence A000292). 

Goldbach's idea is later also applied by some authors. Firstly, notice that 
$IP$ is indirectly  proved by S. W. Golomb in 1963
(\cite[the sequence (1)]{go2}, also see \cite[Section 2.5]{as})
which was constructed a {\it recursive sequence} whose terms are 
pairwise relatively prime and it present 
a generalization of Fermat numbers.  (cf.  Sloane's sequence A000289).
Analyzing the prime factors of $a^n-1$ for given integer $a>1$
and different integer values $n\ge 1$, in 2004 T. Ishikawa, N. Ishida and
Y. Yukimoto \cite[Corollary 3]{iiy} proved that there are
infinitely many primes.
Further, in 2007, for given $n\ge 2$ 
M. Gilchrist \cite{g} constructed the so called 
${}^*$--{\it set of positive integers} $a_1,a_2,\ldots ,a_n$ 
satisfying  $a_j\mid a_i-a_j$ for all distinct $i$ and $j$ with 
$1\le i,j\le n$, and showed that 
 the numbers $b_k:=2^{a_k}+1$, $k=1,2,\ldots ,n$ are mutually prime.
Consequently, the set of primes is infinite. In a similar way,
using the fact that for any integer $n>1$,
$n$ and $n+1$ are mutually prime, and repeating this 
to $n(n+1)$ and $n(n+1)+1$ etc., in 2006 F. Saidak \cite{s}
(for a generaliaztion of this proof, see  \cite[pp. 26--27]{mol}) 
proved the infinitude of primes.  Recently, J. M. Ash and T. K. Petersen 
\cite[Examples 4a)-4e)]{asp} 
proved $IP$ by presenting similar  recursively defined sequences
of positive integers. 
For a construction of some infinite coprime sequences 
see the paper \cite{l} of N. Lord in 2008.


\subsection{Proofs of $IP$ based on algebraic number theory arguments}

In 1736 L. Euler was derived  second proof of  Euclid's theorem 
(published posthumously in 1862 
\cite{eu3} (also see \cite[Sect. 135]{eu5} and  \cite[p. 413]{d}) by 
using the {\it totient function $\varphi(n)$}, defined as the number of 
positive integers not exceeding $n$ and relatively prime to $n$ 
(Sloane's sequence A000010); 
for a proof also see \cite[pp. 134--135]{bu}, \cite[page 3]{pol}. 
As noticed by Dickson \cite[p. 413]{d} (see also \cite[page 80]{sc}),
 this proof is also attributed in 1878/9 by Kummer \cite{k} 
who gave essentially Euler's argument. The proof is based on the 
multiplicativity of the $\varphi$-function. Namely, if $p_1,p_2,\ldots,p_n$ 
is a list of distinct $n\ge 2$ primes with product $P$,  then 
   $$
\varphi(P)=(p_1-1)(p_2-1)\cdots (p_n-1)\ge 2^{n-1}\ge 2.
  $$
This inequality says there exists an integer in the range $[2,P]$ 
that is relatively prime  to $P$, but such an integer has a prime 
factor necessarily different from any of the $p_k$ with $k=1,2,\ldots,n$.
This yields $IP$.  

Euler's idea is in 2009 applied  by J. P. Pinasco \cite{p}.
Assuming that $p_1,p_2,\ldots,p_n$ are all the primes and 
using the  Inclusion-Exclusion Principle, Pinasco 
derived the formula for number of integers in the interval $[1,x]$ that are 
divisible by at least of one of primes $p_i$, which yields
  $$
[x]-1=\sum_{i}\left[\frac{x}{p_i}\right]-
\sum_{i<j}\left[\frac{x}{p_ip_j}\right]+
\sum_{i<j<k}\left[\frac{x}{p_ip_jp_k}\right]-\cdots +(-1)^{n+1}
\left[\frac{x}{p_1p_2\cdots p_n}\right]
  $$
($[\cdot]$ denotes the greatest integer function), 
whence letting $x\to\infty$ easily follows that $1>1$; a contradiction.
Using the  identity
$\sum_{n=1}^{\infty}\mu(n)\left[x/n\right]=1$  
established in 1854 by E. Meissel \cite{mei} (cf. also 
\cite[the formula (3.5.14)]{sh3}),
in 2012 the author of this article 
\cite{me2} presented a very short ``Pinasco's revisted'' proof of $IP$. 
Furthermore, the author  \cite[Remark]{me2} noticed that 
a quite similar proof of $IP$ also follows  using 
 Legendre's formula stated  in the modern form 
\cite[p. 33, Theorem 1.17]{na} as 
$$
\pi(n)-\pi(\sqrt{n})=\sum_{d\mid \Delta}\mu(d)\left[x/d\right]-1
 $$
($\pi(n)$ denotes the number of primes not exceeding $n$).  

Using {\it Theory of Commutative Groups}, in 1888 J. Perott (\cite{pe2}, 
\cite[pp. 303--305]{pe3}; also cf. \cite{c})
showed that, if $p_1,p_2,\ldots, p_n$
are primes, then there exist at least $n-1$ primes between 
$p_n$ and $p_1p_2\cdots p_n$. 

Using {\it Euler theorem}
which asserts that $a^{\varphi(n)}\equiv 1(\bmod{\,n})$ with 
relatively prime integers $a$ and $n\ge 1$, in 1921 G. P\'{o}lya 
\cite[pp. 19--21]{poly} (also see \cite[pp. 131, 324,  Problem 107]{ps})
proved that the set of primes dividing the integer values of
the exponential function $ab^x+c$ $(x=0,1,2,\ldots)$ with integer 
coefficients $a\not=0$, $c\not=0$ and $b\ge 2$ is infinite.

Another proof of $IP$, based on the divisibility property 
$n\mid \varphi(a^n-1)$ ($a,n>1$ are integers),
 is given in 1986 by M. Deaconescu and J. S\'{a}ndor
\cite{ds} (see also \cite{sa}).
Notice that the $\varphi$-function is applied 
by G. E. Andrews \cite[p. 102, Theorem 8-4]{an}  to give an elementary proof
that $\lim_{x\to\infty}\pi(x)/x=0$, where 
$\pi(x)$ is the {\it prime-counting function} 
defined as the number of primes not exceeding $x$ ($x$ is any real number).
In other words, the ``probability" that a randomly chosen 
positive integer is prime is 0.
Using the {\it Inclusion-Exclusion Principle}, this result is
by  an elementary way  also proved by A. M. Yaglom and 
 I. M. Yaglom  \cite[pp. 34, 209--211, Problem 94]{yy}

It was noticed in \cite[p. 4, Exercise 1.2.1]{pol}
that adapting Euclid's proof of $IP$,
it can be proved that for every integer $m\ge 3$, there exist infinitely
many primes $p$ such that $p-1$ is not divisible by $m$. This result
is generalized by A. Granville  
(\cite[p. 4, Exercise 1.2.2]{pol}, \cite[p. 168]{has}; 
also cf. \cite[p. 4, Exercises 1.3 a]{gr2}) to prove that if 
$H$ is a proper {\it subgroup} of the {\it multiplicative group} 
$\Bbb{Z}/m\Bbb{Z}^{*}$ of elements $(\bmod{\, m})$, then there 
exist infinitely many primes $p$ with $p\,(\bmod{\,m})\notin H$.

Similarly, considering {\it order} of $a(\bmod{p})$ in the multiplicative 
group modulo $p$, in 1979 A. Weil  \cite[p. 36, Exercise VIII.3]{wei} 
proved that if $p$ is an odd prime
divisor of $a^{2^n}+1$, with $a\ge 2$  and $n\ge 1$, then $p-1$ is divisible
by $2^{n+1}$. This immediately yields $IP$.

Using {\it Euler's theorem}, it can be proved by induction
that the sequence $2^n-3,$ $n=1,2,\ldots$ contains
an infinite subsequence whose terms 
are pairwise relatively prime (Problem 3 proposed on 
International Mathematical Olympiad (IMO) 1971 
\cite[pages 70 and 392--393]{djp}).
Another less known proof is based on 
{\it Lagrange theorem} on order of subgroup
of a finite {\it group} and {\it Mersenne number} $2^p-1$
with a prime $p$ as follows.  Namely, using Lagrange theorem it can be 
shown that each prime divisor $q$ of $2^p-1$ divides $q-1$,
and so $p<q$, which implies $IP$;
using this fact, we can   inductively obtain  an infinite increasing 
sequence $(p_n)$ 
of primes assuming that $p_{n+1}\mid 2^{p_n}-1$ for each $n=1,2,\ldots$. 
This proof can be found in \cite[p. 3, Second proof]{az}, 
\cite[p. 32, Proposition 1.30 and p. 72, Theorem 1.50]{aaf} 
and at webpage \cite{dp}. 
Mersenne numbers 
(\cite[pp. 75--87, Ch. VII]{r2}, \cite[pp. 109--110]{gkp})
$2^n-1,n=1,2,\ldots$  and the numbers  $2^p-1$ with $p$ prime
  form  Sloane's sequences A000225 and A001348, respectively; 
also see related sequences  A000668,  A000043, A046051 and 
A028335).

Similarly, in 1978 Mohanty \cite[Theorem 2]{mo1} proved 
that for any prime $p>3$, every prime divisor of 
$(2^p+1)/3$ is greater than $p$, and this together with the previous argument
yields $IP$. 

Using the {\it Theory of periodic continued 
fractions} (cf. related Sloane's sequence A003285) and the 
{\it Theory of negative Pell's equations} $x^2-dy^2=-1$, 
in 1976 C. W. Barnes \cite{bar} proved $IP$. Namely, supposing 
that $p_1=2,p_2,\ldots,p_k$ are all the primes with a product $2Q$, Barnes
proved  that  $Q^2+1$ cannot be a power of two; but T. Yamada 
\cite[p. 8]{y} noticed that this fact is obvious since  
$Q^2+1\equiv 2{(\bmod\, 4)}$.

A proof of D. P. Wegener \cite{we} of 1981 based on a study 
of the sums of the legs of {\it primitive Pythagorean triples} 
also contains  Euclid's idea (these triples are triples $(x,y,z)$
of positive integers such that $x^2+y^2=z^2$ and $x$ and $y$
are relatively prime; cf. \cite[Ch. 2, pp. 31--34]{cl}). 

We also point out an interesting result 
established as a solution of advanced problem in  
\cite[pp. 110--111, Problem 37 (a)]{aaf}; namely, this 
result (with two solutions) asserts that if $a$ and $b$ are 
relatively prime positive integers, 
then in the arithmetic progression $a+nb$, $n=1,2,\ldots$ 
there are infinitely many {\it pairwise relatively prime} terms,
which yields $IP$.

Washington's proof of Euclid's theorem from 1980 
(\cite{wa}, \cite[pp. 11--12]{r2}) is 
via {\it commutative algebra}, applying elementary facts of the {\it Theory of
principal ideal domains}, 
{\it unique factorization domains}, {\it Dede}-\
{\it kind domains} and {\it algebraic numbers},  may be found in \cite{sam}. 
Namely, using the fact that 
$(1+\sqrt{-5})(1-\sqrt{-5})=2\times 3$ in the ring of algebraic 
integers $a+b\sqrt{-5}$ ($a,b\in\Bbb Z$) 
(i.e., in the field  of numbers $a+b\sqrt{-5}$ ($a,b\in\Bbb Q$)), 
it follows that 
this ring is not a unique factorization domain. Hence, it is not a
principal domain, whence Washington deduced $IP$. 
The algebraic arguments applied in  this proof are exposed
and well studied  in 2001 by B. Chastek \cite{chas}. 

Quite recently in 2011,  applying two simple lemmas in 
the {\it Theory of Finite Abelian Groups} related to 
the product of some cyclic groups $\Bbb Z_m$, 
R. Cooke \cite{c}  modified Perott's proof noticed above,  
to establish that there are at least $n-1$ primes between the $n$th prime and
the product of the first $n$ primes.

A  ``dynamical systems proof"  due to  S. Srinivasan
(\cite{sr},  also see \cite{y})  in 1984 
uses uses  a polynomial method   and {\it Fermat little theorem}.  Srinivasan
 constructed the sequence $(a_n)$ of positive integers satisfying 
$a_i\mid a_{i+1}$ and $a_i\mid a_{i+1}/a_i$ for each $i=1,2,\ldots$.
Then we immediately see that the sequence $(a_{n+1}/a_n)$ contains 
no two integers which has a nontrivial common divisor. This yields $IP$.

In 2011 P. Pollack \cite{pol2} considered a {\it M\"{o}bius pair} of
arithmetic functions $(f,g)$; that is, functions satisfying 
$f(n)=\sum_{d\mid n}g(d)$ for all $n=1,2,\ldots$, and hence, one can express 
$g$ in terms of $f$ by the {\it M\"{o}bius inversion formula} 
\cite{nzm}. 
Then Pollack deduced  $IP$ by proving the {\it uncertainty principle 
for the M\"{o}bius transform} which asserts that
the functions $f$ and $g$ that become M\"{o}bius pair  
cannot both be of finite support unless they 
both vanish identically. The strategy of Pollack proof goes 
back to J. J. Sylvester \cite{sy2} in 1871, who using certain identities 
between rational functions, gave an argument in the same 
spirit for $IP$  of the form $p\equiv 3(\bmod{\,4})$ and  
$p\equiv 5(\bmod{\,6})$ (cf. Remarks (ii) in \cite{pol2}).

In 2011 R. M. Abrarov and S. M. Abrarov \cite[p. 9]{aa} 
deduced $IP$  applying Euclid's idea to the  identity 
$\mu(n)=-\sum_{i,j=1}^{\sqrt{n}}\mu(i)\mu(j)\delta\left(\frac{n}{ij}\right)$
($n\ge 2$) obtained in their earlier paper \cite[the identity (11)]{aa2}
(also see  \cite[p. 2, the identity (3)]{aa}),
 involving the {\it M\"{o}bius function} $\mu(n)$
(defined so that $\mu(1)=1$, $\mu(n)=(-1)^k$ if $n$ is a product of 
$k$ distinct primes, and $\mu(n)=0$ if $n$ is divisible by the square of 
a prime), and the {\it delta function} $\delta(x)$
(defined as $\delta(x)=1$ if $x\in \Bbb N_0:=\{0,1,2,\ldots,\}$, 
and $\delta(x)=0$ if $x\notin \Bbb N_0$).  
In the same paper, the authors proved $IP$ \cite[p. 9]{aa} as 
an immediate consequence of  \cite[the formula (26)]{aa} 
for the asymptotic density of prime numbers. Their third
 proof \cite[pp. 9--10]{aa} follows from \cite[p. 2, the formula (4)]{aa}
related to the  prime detecting function.  

\subsection{Proof of $IP$ based on Euler's idea on the divergence 
of the sum of prime reciprocals and Euler's formula}

Notice that the proofs of Euclid's theorem presented in 
the previous subsections  are mainly elementary.
On the other hand, there are certain proofs of Euclid's theorem 
that are based on ideas from Analytic Number Theory.
A more sophisticated proof of 
Euclid's theorem  was given many centuries later 
by the Swiss mathematician  Leonhard Euler. 
In 1737 Euler in his work \cite[pp. 172--174]{eu2} (also see \cite{eu})
 showed that by adding the reciprocals 
of successive prime numbers you can attain a sum greater than any 
prescribed number; that is, in terms of modern {\it Analysis},
the sum of the reciprocals of all the primes is {\it divergent}
(cf. \cite[page 8]{r2}, \cite[pp. 8--9]{ew}). 
For more information on Euler's work on infinite series 
see \cite{var}.
{\it Briefly, Euler considered the possibly infinite product 
$\prod 1/(1-p^{-1})$, where the index $p$ runs over all primes.
He expanded the product to obtain the divergent infinite 
{\it harmonic series} $\sum_{n=1}^{\infty}1/n$, 
concluded the infinite product was 
also divergent, and from this concluded that the infinite series 
$\sum 1/p$ also diverges. This can be written symbolically as
$$
\frac{1}{2}+\frac{1}{3}+\frac{1}{5}+\frac{1}{7}+\frac{1}{11}+\frac{1}{13}+
\frac{1}{17}+\frac{1}{19}+\cdots =+\infty.
 $$} 
A result related to this  divergence was refined in 1874 by 
by F. Mertens \cite{mer} (see also \cite[p. 351, Theorem 427]{hw}); namely,  
by {\it Mertens' second theorem},  as $n\to\infty$
the sum $\sum_{p\le n}1/p -\log\log n$ (taken over all primes $p$ not 
exceeding  $n$) converges to the {\it Meissel-Mertens constant} 
$M=0.261497\ldots$ (also known as {\it the Hadamard-de la Vall\'{e}e-Poussin 
constant}).

Using the Euler's idea, in 1888 J.  J. Sylvester \cite{sy3} (also cf. 
\cite[p. 7, Sixth proof of Theorem 1.1]{na})
observed that 
  $$
\prod_{p\le x}\left(1-\frac{1}{p}\right)^{-1}=\prod_{p\le x}
\left(1+\frac{1}{p}+\frac{1}{p^2}+\cdots\right)\ge\sum_{n\le x}\frac{1}{n}
\ge \log x,
 $$
(where the product runs over all primes $p$ not exceeding 
$x$), and since $x$ may be arbitrarily large, the set of primes must be 
infinite. Using the above estimate and the convergence of the series 
$\sum_{n=1}^{\infty}1/n^2$, in the same paper 
J.  J. Sylvester \cite{sy3} 
(also cf. \cite[pp. 11--12, Second proof of Theorem 1.4]{na})
easily proved that the product $\prod_{p\le x}\left(1+1/p\right)$ tends 
to infinity as $x\to\infty$. This impies $IP$. 

A correct realization of Euler's idea was 
presented   by L. Kronecker in his lectures in 1875/76 
(\cite{kro}; also see \cite[pp.269--273]{has} and  
\cite[p. 413, Ch. XVIII]{d}).  
Kronecker noted that ``Euler's'' proof also follows  from the 
{\it Euler's formula}  
$$
\sum_{n=1}^{\infty}\frac{1}{n^s}=\prod\left(1-\frac{1}{p^s}\right)^{-1}
\, (s>1),
  $$ 
where the product on the right is taken  over all primes
$p$ (the first formula in the next  subsection), 
and the fact that the series $\sum_{n=1}^{\infty}1/n^s$ diverges 
for each $s>1$. For some discussion of the history of this formula in relation 
to the infinitude of primes, see \cite{c1}. As noticed 
by Dickson \cite[p. 413]{d}, in 1887/8 L. Gegenbauer \cite{ge}  
proved $IP$ by means of $\sum_{n=1}^{\infty}1/n^s$. 
Dickson \cite[p. 413]{d}) remarked that in 1876 
R. Jaensch \cite{ja} repeated Euler's argument, also ignoring 
convergency. 

Other elementary proofs of the fact that
the sum of reciprocals of all the primes diverges 
were given  in 1943 by R. Bellman \cite{b2}, in 1956 by E. Dux \cite{dux},
 in 1958 by L. Moser \cite{mos}, in 1966 by J. A. Clarkson \cite{cla} and 
in 1995 by D. Treiber \cite{tr}. A survey of some  these 
proofs was given in 1965 by  T. \v{S}al\'{a}t \cite{shal}.
Furthermore, in 1980 C. Vanden Eynden \cite{van}
considered Euler's type product of all the expressions of the form
$\left(1+1/p\right)\sum_{k=0}^{\infty}1/p^{2k}=\sum_{j=0}^{\infty}
1/p^j$, where $p$ ranges over the set of all primes not 
exceeding $x$. This equality together with the divergence of the 
series $\sum_{n=1}^{\infty}1/n$ and the convergence of the series     
$\sum_{n=1}^{\infty}1/n^2$ easily yields the divergence of the
 sum of the reciprocals of all the primes. 

It is interesting to notice that in actual reality, 
Euler never presented his work as a proof of Euclid's theorem,
though that conclusion is clearly implicit in what he did.
Euler's remarkable proof of $IP$, based on formal identity 
$\prod_p 1/\left(1-1/p\right) =\sum_{n=1}^{\infty}1/n$,
amounts to unique factorization, and it is also 
discussed at length by R. Honsberger in his book \cite{h} and modified in 
2003 by C. W. Neville \cite[Theorem 1(a)]{ne1}.  
In 1938 P. Erd\H{o}s (\cite{er3}; also see
\cite[pp. 5--6, Sixth proof]{az}, \cite[p. 17, Theorem 19]{hw} and 
\cite[pp. 12--13]{pol}) gave
an elementary ``counting"  proof of the divergence of the sum 
of reciprocals of primes, and consequently, the set of all
primes is infinite.  P. Pollack \cite[p. 11]{pol} pointed out 
that it is remarkable that this method of proving $IP$ 
(in contrast with Euclid's proof, for instance) is independent 
of the additive structure of the integers. 

\vspace{2mm}

\noindent{\bf Remarks.}
Notice that the the asymptotic behavior of the product of $1-1/p$ 
was given by  {\it Mertens' third theorem}
established in 1874 by F. Mertens \cite{mer} (also see \cite[pp. 351--353, 
Theorem 428]{hw}), which states that  
$\prod_{p\le n}\left(1-1/p\right)\sim e^{-\gamma}/\log n$,
where the product runs over all primes $p$ not exceeding  $n$, and 
$\gamma= 0.577216\ldots$ is 
{\it Euler-Mascheroni constant}. 
An  elementary geometrical proof of Mertens' third theorem with another 
constant $c$ instead of $e^{-\gamma}$, was given in 1954 
by A. M. Yaglom and  I. M. Yaglom  \cite[pp. 41; 194--196, Problem 174]{yy2}. 
 Using Mertens' third theorem 
(with the  constant $e^{-\gamma}$), in \cite[p. 42]{yy2} the authors also derived   
the  curious formula $\prod_{p\le n}\left(1+1/p\right)\sim 
(e^{\gamma}\log n)/\pi^2$, as $n\to\infty$. \hfill$\Box$
\vspace{2mm}

Furthermore, using the classical {\it Chebyshev's argument}
based on the well known {\it de Polignac's formula} 
(attributed by Dickson[p. 263, Ch. IX]{d} to A.-M. 
Legendre \cite[p. 8]{leg} in 1808)
for the exponent $\nu_p(n!)$ of prime $p$ dividing the factorial
$n!$ given as $\nu_p(n!)=\sum_{k=1}^{\infty}\left[n/p^k\right]$,  
a short proof that the sum $\sum_{p}\log p/p$ diverges
due to P. Erd\H{o}s is presented in \cite[$8^{\mathrm th}$  proof]{y}
and this yields $IP$. 
Similarly, using de Polignac's formula, in 1969  
\cite[pp. 613--614, Remark 6]{coh} (cf. also 
\cite[p. 54, Exercise 1.21]{cp}) E. Cohen 
gave a short simple proof that the series $\sum_{p}\log p/p$ diverges
(the sum ranges over all the primes), which yields $IP$. 
This result also follows from {\it Mertens' first theorem}
obtained in 1874 by F. Mertens \cite{mer}, which asserts that the quantity 
$|\sum_{p\le n}\log p/p-\log n|$
is bounded, in fact $<4$ (for an elementary proof, see \cite[pp. 
171, 183--186, Problem 171]{yy2}). Notice that this result 
immediately follows from  Mertens' second theorem.

Further, combining the Euler's idea with the geometrical 
interpretation of {\it definite integral} $\int_{1}^{x}(1/t)\,dt=\log x$
with $n\le x<n+1$ in their Problems book \cite[p. 4, Fourth proof]{az}
A. M. Yaglom and  I. M. Yaglom
proved the inequality $\log x\le \pi(x)+1$, where 
$\pi(x)$ is the  prime-counting function.
This inequality immediately yields $IP$.

Another modification of Euler's proof, involving the 
{\it  logarithmic complex function}, can be found  in book \cite[p. 35]{cp}
of R. Crandall and C. Pomerance. 

\vspace{2mm}

\noindent{\bf Remarks.} Notice that from Euclid's proof  
(see e.g., \cite[p. 12, Theorem 10]{hw})
easily follows that $\pi (x)\ge \log_2\log_2 x$ for each $x>1$,
and the same bound follows more readily from the 
Fermat numbers proof. Of course, this is a horrible bound.
From the Erd\H{o}s's proof \cite{er3}  given above
it can be easily deduced the bound  
$\pi(x)\ge \log x/(2\log 2)=\log_2 x/2$ for each $x\ge 1$
\cite[p. 17, Theorem 20]{hw}. This estimate can be improved using 
{\it Bonse's inequality} presented above.  
Namely, applying induction, it follows from this inequality that 
$p_n\le 2^n$; so, given $x\ge 2$, taking  $x=2^n+y$ with $0\le y<2^n$, we 
find that $\pi(x)\ge \pi(2^n)\ge n\ge \log_2 x-1$.\hfill$\Box$    

\vspace{2mm}
\noindent{\bf Remarks.}
Recall that an extremely  difficult problem in Number Theory is the 
{\it distribution of the primes} among the natural numbers. This problem 
involves the study of the asymptotic behavior of the counting function
$\pi(x)$ which is one of the more intriguing functions in Number Theory. 
For elementary methods in the study of the distribution 
of prime numbers, see \cite{dia}.
Studying tables of primes,
 C. F. Gauss in the late 1700s and A.-M. Legendre in the early 1800s 
conjectured the celebrated {\it Prime Number Theorem}:
$\pi(x)=|\{p\le x:\, p\,\,{\rm prime}\}|\sim x/\log x$
($|S|$ denotes the cardinality of a set $S$). This theorem was proved 
much later (\cite[p. 10, Theorem 1.1.4]{cp}; for its simple analytic 
proof see \cite{new} and \cite{za}, and for its history see \cite{bd}
and \cite{golds}).
Briefly, $\pi (x)\sim x/\log x$ as $x\to\infty$, or in other words, 
the density of primes $p\le x$ is $1/\log x$; that is, the ratio
$\pi(x): \left(x/\log x\right)$ converges to 1 as $x$ grows without bound.
Using L'H\^{o}pital's rule, Gauss showed that the logarithmic integral 
$\int_{2}^x dt\log t$, denoted by $\mathrm{Li}(x)$, is asymptotically
equivalent to $x/\log x$. Recalll that Gauss felt that 
$\mathrm{Li}(x)$ gave better approximations to $\pi(x)$ than 
$x/\log x$ for large values of $x$.
 Though unable to prove the Prime Number Theorem, 
several significant contributions to the proof of Prime Number Theorem 
were given by P. L. Chebyshev in his two important 1851--1852 papers 
(\cite{che1} and \cite{che2}). Chebyshev proved that there exist positive 
constants $c_1$ and $c_2$ and a real number $x_0$ such that 
  $c_1x/\log x\le\pi(x)\le c_1x/\log x$ for $x>x_0$.  In other words, 
$\pi(x)$ increases as $x\log x$. Using methods of complex analysis and 
the ingenious ideas of Riemann (forty years prior), 
this  theorem was first proved in 1896, independently 
 by J. Hadamard and C. de la Vall\'{e}e-Poussin  
(see e.g., \cite[Section 4.1]{pol}).\hfill$\Box$

\subsection{Proof of $IP$ based on Euler's product for the  Riemann 
zeta function   and the irrationality of  $\pi^2$  and $e$}

Proofs of $IP$ presented in this subsection involve 
the {\it Riemann zeta function} (for $\Re (s)>1$, to ensure convergence)
 defined as $\zeta(s):=\sum_{n=1}^{\infty}1/n^s$. Riemann introduced the 
study of $\zeta (s)$ as a function of a complex variable 
in an 1859 memoir on the distribution of primes \cite{ri}.  However, the 
connection between the zeta function and the primes goes back earlier.
Over a hundred years prior, Euler had looked at the same series 
for real $s$ and had shown that  \cite[Theorema 8]{eu2} 
  $$
\sum_{n=1}^{\infty}\frac{1}{n^s}=\prod_{p}\frac{1}{1-\frac{1}{p^s}}\quad (s>1).
  $$
This is the {\it Euler's factorization} which is often called 
an analytic statement of unique factorization (this is a consequence 
of a well known standard uniqueness theorem for Dirichlet series 
\cite[Theorem 11.3]{apo2}). 

  Dickson \cite[p. 414]{d} (also see \cite[p. 10]{pol}) 
noticed that in 1899 J. Braun \cite{bra}
attributed to J. Hacks a proof of $IP$ by means of the Euler's formula 
$\sum_{n=1}^{\infty}1/n^2=\pi^2/6$ (for  elementary proofs of this formula 
see \cite{cho}, \cite{gi} and \cite{li}) and  the  Euler's factorization 
$\prod 1/(1-p^{-2})=\sum_{n=1}^{\infty}1/n^2$ 
(Sloane's sequence A013661) and the {\it irrationality of $\pi^2$} 
proved in 1794 by Legendre \cite{leg2} 
(also see \cite[p. 47, Theorem 49]{hw}, \cite[p. 285]{r3}). 
Namely, if there were only finitely many primes,
then $\zeta(2)$ would be rational; a contradiction. 
Notice also that this proof was reported in 1967 in the 
reminiscences of Luzin's Moscow school of mathematics 100 years ago 
by L. A. Lyusternik \cite[p. 176]{ly} (also cf. \cite[p. 466]{c1})
which ascribed this proof to A. Y. Khinchin. Such proofs attract interest 
because they make unexpected connections. According to Lyusternik, 
 ``exotic" proofs of $IP$ were a routine challenge among Luzin's
students, and many such proofs were found. But apparently no one 
thought of publishing them.
The previous equality  is in fact, the  well known 
{\it Euler's formula} (or {\it Euler's product})
\cite[p. 245]{hw}   for the  Riemann zeta function 
$\zeta(2):=\sum_{n=1}^{\infty}1/n^2$ \cite[p. 246, Theorem 280]{hw}.
The same proof of $IP$ was also presented in  2007 by J. Sondow \cite{so1}.
Notice  that, applying the same argument 
for the product formula 
$\prod 1/(1-p^{-3})=\sum_{n=1}^{\infty}1/n^3:=\zeta(3)$ together with a 
result of R. Ap\'{e}ry in 1979 \cite{ap} that $\zeta(3)$ is irratioanl,
we obtain $IP$. 
  
Further, using the Euler's formulas  for $\zeta(2)$ and 
$\zeta(4)=\sum_{n=1}^{\infty}1/n^4=\pi^4/90$ \cite[p. 245]{hw} 
(Sloane's sequence A0013662), it 
can be easily obtained that $5/2=\prod_{p}\left((p^2+1)/(p^2-1)\right)$,
where the product is taken  over all the primes  \cite[p. 11]{pol}.
 In 2009 P. Pollack \cite[p. 11]{pol}  observed that if the set of 
all primes is finite, then the numerator of the ratio on
 the right of this formula is not divisible by 3, but its denominator 
is divisible by 3. This contradiction yields $IP$. 
We  recall Wagstaff's (open) question \cite[B48]{gu}  as to 
whether there exists an elementary proof of the previous formula.

Notice that $IP$ can be proved using the formula 
$\prod_{n=1}^{\infty}(1-x^n)^{\mu(n)/n}=e^{-x}$ ($|x|<1$ and 
$\mu(n)$ is the M\"{o}bius function) 
 proposed as a Monthly's Problem in 1943 by R. Bellman \cite{b3} 
and solved in 1944 by R. C. Buck \cite{buc}. 
If we suppose that $p_1,p_2,\ldots,p_k$ are all the primes with a product $P$,
then obviously $\mu(n)=0$ for each $n\ge P$, and so, the previous
formula for $x=-1/2$ becomes 
$\left(\prod_{n=1}^{P}\left(1-(-1/2)^n\right)^{\mu(n)/n}\right)^2=e$. 
The previous equality and the fact  that the number $e$ is irrational
(a result due to J. Fourier in 1815; see e.g., 
\cite[pp. 27--28]{az}) give a contradiction which yields $IP$.  

  \vspace{2mm}

\noindent{\bf Remarks.}
Notice that the above formulae  for $\zeta(2)$ and $\zeta(4)$
are two special cases of the following classic formula 
 discovered by  Euler in 1734/35 \cite{eu}, 
which express  $\zeta(2n)$  as a rational 
multiple of $\pi^{2n}$ involving 
{\it Bernoulli number} $B_{2n}$:   
$\zeta(2n)=(-4)^{n-1}B_{2n}\pi^{2n}/(2\cdot (2k)!)$ ($n=1,2,3,\ldots$).
An elementary proof of this formula 
 for $n=1$ is given by I. Papadimitriou \cite{pa} in 1973  and 
for arbitrary $n$ by T. M. Apostol \cite{apo} in the same year
  (for another elementary evaluations of $\zeta(2n)$ see \cite{be} and 
\cite{os}). For instance, since $B_2=1/6$ and $B_4=1/30$, we find that 
$\zeta(2)=\pi^2/6$ and $\zeta(4)=\pi^4/90$, respectively.\hfill$\Box$

\subsection{Combinatorial proofs of $IP$ based on enumerative arguments}

Several combinatorial proofs of $IP$ involve simple counting 
arguments. More precisely, 
these proofs are mainly based
on counting methods  which are used in them to count
the cardinality  of integers less than a given integer $N$ 
and which satisfy certain  divisibility properties.  
The first such proof, given by   J. Perott  
in 1881 (\cite{pe}, \cite[p. 10]{r2} and \cite[p. 8]{na}) 
is based on the facts that the series $\sum_{n=1}^{\infty}1/n^2$ is 
{\it convergent} with the sum smaller than 2 
and  that there exist exactly  $2^n$ divisors of the product of $n$ distinct 
primes. In his proof Perott also established the estimate 
$\pi(n)>\log_2(n/3)$, where $\pi(n)$ is the number of primes less than 
or equal to $n$.
  Perott's proof was modified in \cite[pp. 11--12]{pol} by eliminating 
use of the formula $\zeta(2)=\pi^2/6$. 
Using Perott's method, in 2006
L. J. P. Kilford \cite{ki} presented a quite similar 
 proof of $IP$ based on the fact that 
for any given $k\ge 2$, the sum $\sum_{n=1}^{\infty}1/n^k$
converges to a real number which is strictly between 1 and 2.

A classical  proof of $IP$  which is combinatorial  in spirit
and entirely elementary,   
was given by Thue in  1897  in his work \cite{th} 
(also see \cite[414]{d} and \cite[page 9]{r2}).
This proof  uses a ``{\it counting method}" and 
the {\it fundamental theorem of unique factorization} of positive
integers as a product of prime numbers as follows. {\it Choose integers
$n,k\ge 1$ such that $(n+1)^k<2^n$ 
and set $m=2^{e_1}\cdot 3^{e_2}\cdots p_r^{e_r}$, where
we assume that  $2<3<\cdots <p_r$ is a set of all the primes 
and $1\le m\le 2^n$. Suppose that $r\le k$. 
Since $m\le 2^n$, we have $0\le e_i\le n$ for each $i=1,2,\ldots, r$. 
Then counting all the possibilities, it follows that 
$2^n\le (n+1)n^{r-1}<(n+1)^r\le (n+1)^k<2^n$. 
This contradiction yields $r\ge k+1$. Now taking $n=2k^2$,
then since $1+2k^2<2^{2k}$ for each $k\ge 1$, it follows that
$(1+2k^2)^k\le 2^{2k^2}=4^{k^2}$, and so there at least 
$k+1$ primes $p$ such that $p<4^{k^2}$. Thus, letting $k\to\infty$
yields $IP$}. 

Applying a  formula for the number of positive integers less than $N$ 
given in \cite[Ch. XI]{cah}, in 1890 J. Hacks \cite{ha2} 
(see also \cite[p. 414]{d}) proved $IP$.

In order to prove $IP$, 
similar {\it enumerating arguments} to those of Thue 
were used in a simple Auric's proof, 
which appeared in 1915 \cite[p. 252]{au} (also see \cite[p. 414]{d}, 
\cite[page 11]{r2}), as well by P. R. Chernoff in 1965 \cite{ch}, 
M. Rubinstein \cite{ru}  in  1993 and 
M. D. Hirschorn \cite{hi} in 2002.
A proof of $IP$ similar to that of Auric is
given in 2010 by M. Coons \cite{coo}.

Using a combinatorial argument, the  unique factorization theorem 
and the {\it pigeonhole principle}, 
$IP$ is recently proved by D. G. Mixon \cite{mi}.

A less known  elementary result of P. Erd\H{o}s \cite[p. 283]{er2}
(also see \cite{er1}) in  1934,
 based on {\it de Polignac's formula} (actually due to A.-M. Legendre), 
asserts that there is a prime between $\sqrt{n}$ and $n$ for each positive integer $n>2$.
In the same paper Erd\H{o}s proved that if $n\ge 2k$, then 
${n\choose k}$ contains a prime divisor greater than $k$.
In particular, this fact for $n=2k$ obviously yields $IP$. 
Notice also that $IP$ follows by two results of W. Sierpi\'{n}ski
from his monograph in 1964 \cite{si}. Namely, if we suppose that there are 
a total of $k$ primes, then by \cite[page 132--133, Lemmas 1 and 4]{si}, 
we have $4^n/{2\sqrt{n}}<{2n\choose n}\le (2n)^k$ for each positive 
integer $n>1$. This contradicts the fact that $4^n/({2\sqrt{n}})\ge (2n)^k$
for sufficiently large $n$.

In 2010 J. P. Whang \cite{w} gave a short proof of $IP$
 by using  de Polignac's formula.

\subsection{Furstenberg's topological proof of $IP$ and 
its modifications}

A  proof of Euclid's theorem due to
H. Furstenberg  in 1955 (\cite{f}; also see
 \cite[pp. 12--13]{r2}, \cite[p. 12]{pol} or 
\cite[p. 5]{az}) is a short ingenious proof based on {\it topological ideas}.
In order to achieved a contradiction,
Furstenberg introduced a {\it topology} on 
the set of all integers, namely the smallest topology 
in which any set of all terms of a nonconstant {\it arithmetic 
progression} is {\it open}. Here we quote this proof in its 
entirety: ``{\it In this note we would like to offer an elementary 
``topological" proof of the infinitude of prime numbers. We introduce
a topology into the space of integers $S$, by using the arithmetic 
progressions $($from $-\infty$ to $+\infty$$)$ as a basis. It is not difficult
to verify that this actually yields a topological space. In fact under 
this topology $S$ may be shown to be normal and hence metrizable.
Each arithmetic progression is closed as well as open, since its complement
is the union of other arithmetic progressions $($having the same difference$)$.
As a result the union of any finite number of arithmetic progressions is 
closed. Consider now the set $A=\bigcup A_p$, where $A_p$ consists 
of all multiples of $p$, and $p$ runs though the set of primes $\ge 2$.
The only numbers not belonging to $A$ are $-1$ and $1$, and since the
set $\{-1,1\}$ is clearly not an open set, $A$ cannot be closed. Hence $A$ 
is not a finite union of closed sets which proves that there are an infinite
of primes}.''

In 1959 S. W. Golomb \cite{go}  developed further the idea
of  Furstenberg and gave another prooof of Euclid's theorem
 using a topology $\mathcal D$ on the set $\Bbb N$ of natural 
numbers with the base $\mathcal B=\left\{\{an+b\}:\, (a,b)=1\right\}$
($(a,b)$ denotes the  greatest common divisor of $a$ and $b$),
defined in 1953 by M. Brown \cite{br}. In the same paper Golomb proved 
that the topology $\mathcal D$ is {\it Hausdorff}, {\it connected}
and {\it not regular}, 
$\Bbb N$ is 
${\mathcal D}$-{\it connected}, and the {\it Dirichlet's theorem}
(on primes in arithmetic progressions) is equivalent to the 
$\mathcal D$-density of the set of primes in $\Bbb N$.
Moreover,  in 1969 A. M. Kirch \cite{kir} proved that the {\it topological
space} $(\Bbb N,\mathcal D)$ is not {\it locally connected}.

In 2003 D. Cass and G. Wildenberg \cite{cw} 
(also cf. \cite{kl}) have shown that Furstenberg's proof
can be reformulated in the language of {\it periodic 
functions on integers}, without reference to topology.
This is in fact, a beautiful combinatorial version of 
Furstenberg's proof. Studying arithmetic properties 
of the multiplicative structure of {\it commutative rings} and related
topologies, in  2001  \v{S}. Porubsky \cite{por} established 
new variants of Furstenberg's topological proof.
Notice also that Furstenberg's proof of $IP$ 
is well analyzed  in 2009 by A. Arana \cite{ar}, in 
2008 by M. Baaz, S. Hetzl, A. Leitsch, C. Richter and H. Spohr 
\cite{bhlrs}, and also  discussed in greater detail in 2011 by 
M. Detlefsen and A. Arana \cite{da}. Furthermore, C. W. Neville 
\cite[Theorem 1(a)]{ne1} 
pointed out that this proof has been extended in various directions,
for example, to the setting of {\it Abstract Ideal Theory} see \cite{ks}
and \cite{por}. 

More than 50 years later, in 2009 using Furstenberg's ideas but 
rephrased without topological language, 
I. D. Mercer \cite{m} provided a new short proof that the number of primes 
is infinite.  Finally, notice that Furstenberg's proof is an important 
beginning example in the {\it Theory of profinite groups} (see  book  reviews
by A. Lubotzky \cite{lub} in 2001).

\subsection{Another proofs of $IP$}

Euclid's proof of $IP$  was revisted in 1912/13 
by I. Schur \cite{sch} (see also \cite[pp. 131, 324,  Problem 108]{ps}) 
who showed that the set of primes dividing the integer values of a 
nonconstant integer polynomial is infinite. 
Suppose that $Q$ is  a polynomial with integer coefficients such that 
 $\{p_1,p_2,\ldots ,p_k\}$ is a set of all primes with this property is 
finite. Then assuming that $Q(a)=b\not=0$, we will consider the
integer value $c=\left(Q(a+bp_1p_2\cdots p_k)\right)/b$. 
Then obviously $c\equiv 1(\bmod{\, p_1p_2\cdots p_k})$ and  therefore,
 $c$ has at least one prime divisor, say $p$,  
distinct from every element of the set  $\{p_1,p_2,\ldots p_k\}$.
It follows that the value $Q(a+bp_1p_2\cdots p_k)=bc$ 
is also divisible by $p$; a contradiction.  
In particular, for $Q(x)=x+1$  the previous proof is a copy of Euclid's 
proof of $IP$. If $Q(x)=\Phi_m(x)$ is the $m$th {\it cyclotomic polynomial}, 
then the above proof yields that there are  infinitely many primes 
which are congruent to $1(\bmod{\,m})$ (cf. Section 3).

\vspace{2mm}

\noindent{\bf Remarks.}
In 1990 P. Morton \cite{mort} considered a related problem for 
an integer sequence $(a_n)$ for which there is an integer constant $c$
such that for all $i\in \Bbb Z=\{\ldots -2,-1,0,1,2,\ldots\}$ 
$a_n=i$ holds for almost $c$ values of $n$. If for such   
an  integer sequence $(a_n)$, the so called {\it almost-injective},
define the set ${\mathcal S}(a_n)=\{p\,\,{\rm prime}:\,\, p\mid a_n\,\,
\,\,{\rm for\,\, at \,\, least\,\, one}\,\, n\in\Bbb N\}$,
then Morton \cite{mort} proved that ${\mathcal S}(a_n)$ is infinite 
if $(a_n)$ has at most polynomial growth, i.e., $|a_n|\le an^d$ for 
some positive constants $a$ and $d$. This result is extended quite recently in
2012  by C. Elsholtz \cite{el} for  almost-injective integer sequences 
of  {\it subexponential growth}, i.e., for  
 almost-injective integer sequences $(a_n)$ for which $a_n=o(\log n)$.  
As noticed in \cite[p. 333]{el}, another way to look at this theorem is to 
study  ``primitive divisors'' of integer sequences. Given an integer sequence 
$(a_n)$, a divisor $d$ is called {\it primitive} if $a_i$ is divisible by 
$d$, but  $a_j$ is not divisible by $d$ for any $j<i$. For a good survey 
of this topic, see Chapter 6 of the book \cite{evsw}.
  
However, it is not known whether there are polynomials of degree greater than 
1 with integer coefficients representing infinitely many primes for 
integer argument. Using Chebyshev's 
estimate $\pi(x)\ge x/\log x$ and a simple counting argument, 
in 1964 W. Sierpi\'{n}ski \cite{si1} 
(also see \cite[p. 35, Theorem 1.6.1]{pol}) 
proved that  for every   $N$ there exists an integer $k$ for which 
there are more than $N$ primes represented 
by $x^2+k$ with $x=0,1,2,\ldots$. 
In 1990  B. Garrison \cite{gar} (cf. \cite[p. 36, 
Exercise 1.6.2]{pol})  generalized Sierpi\'{n}ski's result 
to polynomials $x^d+k$ of degree $d\ge 2$ and proved 
that  for any such $d$ and any  $N$ there exists a positive integer $k$  
such that $x^d+k$ ($x=0,1,2,\ldots$) assumes more than $N$ prime values.
P. Pollack \cite[p. 36, Exercise 1.6.2 b)]{pol} noticed 
that the previous assertion remains true if ``positive" is replaced
by ``negative". This obviously implies $IP$. Modifying Garrison's proof, 
in 1992 R. Forman \cite{fo} extended Garrison's result to a large class of 
sequences. Forman \cite[Proposition]{fo} proved that if $f(x)$ is a 
nonconstant polynomial with positive leading coefficient
(the coefficients need not be integers), then 
 for any $N$ there are infinitely many  nonnegative 
integers $k$ such that the sequence  
$[f(n)]+k$ ($n=0,1,2,\ldots$) contains at least  $N$ primes
($[\cdot]$ denotes the greatest integer function). 
Furthermore, in  1993  U. Abel and H. Siebert 
\cite{abs}  also extended Garrison's result. 
They proved that if $f(x)\in \Bbb Z[x]$ is a polynomial of degree $d\ge 2$
with positive leading coefficient, then 
 for every $N$ there exists an integer $k$ for which 
$f(x)+k$ ($x=0,1,2,\ldots$) assumes more than $N$ prime values. 
Their argument of proof depends on counting the number of solutions of 
certain inequalities and shows that no arithmetical properties of 
polynomials are needid other than rate of growth. 
In particular, in  \cite[p. 167, proof of Theorem]{abs} 
it was applied the well known {\it Sylvester's version of the Chebyshev 
inequalities} $0.9\le \pi(x)\log x/x\le 1.1$ (for sufficiently large $x$)
(\cite{sy4}, see also \cite[p. 555, (1.7)]{dia}). 

However, the problem of characterizing the prime divisors 
of a polynomial of degree $>2$ is still unsolved, 
except in certain special cases. We see that if $p$ is any prime
that  does not divide $a$, then $p$ divides each polynomial 
$Q_1(x)=ax+b$ with arbitrary  $b\in \Bbb Z$. Similarly, 
the set of all prime divisors of $Q_2(x)=x^2-a$ can be determined 
by using law  of quadratic reciprocity. 
 Some known and new related results 
for various classes of integer polynomials were presented by 
I. Gerst and J. Brillhart \cite{gb}. \hfill$\Box$

\vspace{2mm}

By Problem 3 proposed on International Mathematical Olympiad
(IMO) 2008 \cite[pages 336 and 776]{djp}, there exist 
infinitely many positive integers $n$ such that $n^2+1$ has a prime 
divisor greater than $2n+\sqrt{2n}$. This immediately yields $IP$.

In  \cite{cha} (see also \cite{cha2} and 
\cite[page 118, Section 10.1.5]{cl}), in 1979 the computer scientist 
G. J. Chaitin gave a  
proof of $IP$ using algorithmic information theory.  
If $p_1,p_2,\ldots ,p_k$ are all the primes, then for a fixed 
$N=p_1^{a_1}p_2^{a_2}\cdots p_k^{a_k}$ Chaitin defines {\it algorithmic
entropy} $H(N):=\sum_{i=1}^ka_i\log p_i$ of $N$, and uses 
various properties, such as {\it subadivity of algorithmic entropy} 
expressed as $H(N)\le \sum_{i=1}^k H(n_i)+O(1)$.  
In order to prove this property, Chaitin 
estimates how many integers $n$ with $1\le n\le N$, could possibly be expressed in the form
$p_1^{b_1}p_2^{b_2}\cdots p_k^{b_k}$. In order for this expression to be at 
most $N$, every exponent has to be much smaller than $N$: precisely, we need 
$0\le b_i\le \log_{p_i}N$; the latter quantity 
is at most $\log_2 N$, so there are at most $\log_2N+1$ choices for each
exponent, or $(\log_2N+1)^k$ choices overall. However, this latter quantity
is much smaller than $N$ for sufficiently large $N$; a contradiction which
implies $IP$. 

We also notice that Chaitin's proof is quite similar to those 
of $IP$ due to L. G.  Schnirelman's   book \cite[pp. 44--45]{shn} 
published posthumously in 1940.
Moreover, a more sophisticated version of Chaitin's proof which uses 
an obvious representation $n=m^2k$ of a positive integer $n$ 
where $k$ is squarefree,  can be found in the book \cite[pp. 16--17]{hw} 
of Hardy and Wright (which was first written in 1938).
A similar idea is used in 2008 by E. Baronov \cite[p. 12, Problem 5]{ibk}
to show that if a sequence of positive integers $(a_n)$ satisfies 
$a_n<a_{n+1}\le a_n+c$, with a fixed $c\in\Bbb N$ and for each $n\in \Bbb N$,
then the set of prime divisors of this sequence is infinite.
This immediately yields $IP$. Similarly, the same author    
\cite[pp. 12--13, Problem 6]{ibk} proved that if $m$ and $n$ are 
positive integers such that $m>n^{n-1}$, then there exist 
distinct primes $p_i$, $i=1,2,\ldots,n$ such that $p_i\mid m+i$ 
for each   $i=1,2,\ldots,n$. This also implies $IP$.

\vspace{2mm}
\noindent{\bf Remarks.}
The  argument in Chaitin's proof also shows that the percentage of 
nonnegative integers up to $N$ which we can express as 
a product of any $k$ primes tends to 0 as $N$ approaches infinity.
Notice that this proof  gives a lower  bound on 
$\pi(x)$ which is between $\log\log x$ and $\log x$ (but much closer to 
$\log x$). Using the same method, the lower bound $\pi (x)\ge (1+o(1))
\log x/(\log\log x)$ was established in 
\cite[p. 15, Proof of Lemma 1.2.5]{pol} 
(cf. also \cite[pp. 15--17, Lemma 0.3 and Exercise 0.5]{hil}).   
In revisted Chaitin's proof H. N. Shapiro  
\cite[pp. 34--35, Theorem 2.8.1]{sh3} obtained the estimate 
$\pi (x)> \log x/(3\log\log x)$ for each $x>e^2$.
\hfill$\Box$
\vspace{2mm}

   In his dissertation, in 1981 A. R. Woods \cite{wo} proved $IP$    
by adding $PHP\Delta_0$ to a {\it a weak system of arithmetic} 
$I\Delta_{0}$, where  $PHP\Delta_0$ stands for the pigeonhole principle 
formulated for functions defined by $\Delta_0$-formulas. 
($I\Delta_{0}$ is the theory over the vocabulary 
$0,1,+,\cdot ,<$ that is axiomatized by basic properties
of this vocabulary and induction axioms for all bounded formulas).
In 1988 J. B. Paris, A. J. Wilkie and A. R. Woods (\cite{pww}; 
also see \cite[pp. 162--164]{daq} and \cite{pw}) replaced Woods' earlier proof 
with one using an even weaker version of the pigeonhole principle. 
They showed that a considerable part of elementary number 
theory, including $IP$, is provable in
 a weak system of arithmetic $I\Delta_{0}$ with the {\it weak pigeonhole 
principle} for $\Delta_0$-{\it definable functions} added as an axiom 
scheme. It is a longstanding open question  \cite{wi}
 whether or not one can dispense with the weak pigeonhole principle,
by proving the existence of infinitely many primes within $I\Delta_0$.
Studying the problem of proving in 
{\it weak theories of Bounded Arithmetic} that there are infinitely many 
of primes, in 2008 P. Nguyen \cite{ng} showed that $IP$ can be proved by some 
``minimal" reasoning (i.e., in the theory 
$<$Emphasis Type=``Bold''>I$\Delta</$Emphasis$><$
\noindent Subscript><Emphasis Type=``Bold''$>0</$Emphasis$><$
Subscript$>)$ using concepts such as (the logarithm) of a binomial coefficient.  

Euclid's  revisted proof of $IP$ via methods of  
{\it nonstandard Analysis} was given by R. Goldblatt \cite{gold} in 1998 
(also see \cite[p. 16, Section 1.2.6]{pol}).  


\subsection{14 proofs of $IP$ (2012--2017)} 

{\bf 1)} In 2012  the author of this article by \cite[Theorem 1]{me4} 
improved Cooke's result \cite[Theorem]{c} (see page 17 of this article), 
 refining the Euler's proof of $IP$ by the following  result: 
``{\it Let $\alpha$ be a real number such that
$1< \alpha <2$ and let $x_0=x_0(\alpha)$ be 
a {\rm(}unique{\rm)} positive solution of the equation 
   $$
 x^{\alpha-1}-\frac{\pi }{e^2\sqrt{3}}x+1=0.
  $$
Then for each  positive integer $n>x_0$ 
there exist at least  $\lfloor n^\alpha \rfloor$ primes between 
the $(n+1)$th prime and the product of the first $n+1$ primes,
where $\lfloor a\rfloor$ denotes the greatest integer 
less than or equal to  $a$. 

Moreover,  for each positive integer $n$ there are at least $n$ primes 
between the $(n+1)$th prime and the product of the first $n+1$ primes.}''

{\bf 2)} In 2015  L. Alpoge \cite{al} establihed  $IP$ as the amusing 
 consequence of the following 
(called by Khinchin \cite{kh} beautiful) theorem of van der Waerden (\cite{wae}; also see 
\cite[Theorem 1]{al}): ``{\it Suppose the positive integers are 
colored with finitely many colors. Then there are arbitrarily many 
arithmetic progressions containing integers all of the same color.}''

More formally, if $f:\Bbb Z^+\to S$ is any function to  a finite set $S$, 
then for each $k>0$, there are $n$ and $d$ for which 
    $$
f(n)=f(n+d)=\cdots = f(n+kd).
     $$

{\bf 3)} Motivated by the previous Alpoge's proof of $IP$, in 
2017 A. Granville \cite[Theorem 1]{gr5} proved $IP$  
combining van der Waerden's theorem with a famous result of Fermat which 
asserts that {\it there are no four-term arithmetic progresssions of distinct 
integer squares} (see, e.g., \cite{sil}).   
    \vspace{2mm}

{\bf 4)} Proceeding in a similar way as in Saidak's proof of $IP$ 
(see Subsection 22, p. 14 of this article),  
in 2015  B. Maji \cite{maj} constructed an infinite 
sequence of pairwise relatively prime positive integers.
This fact immediately yields $IP$.
  \vspace{2mm}  

{\bf 5)} Assuming that the set of all primes is finite,
  in 2015 S. Northshield \cite{nor} proved $IP$ by
considering the product 
  $$
\prod_{p}\sin\left(\frac{\pi}{p}\right),
 $$
where $p$ runs over all primes (``a one-line proof''). 
   \vspace{2mm}

{\bf 6)} In 2016  A. R. Booker \cite{boo} considered a generalization of 
Euclid's proof of $IP$ and showed that it leads to variants 
of the Euclid-Mullin sequence that provably contain 
every prime number. Namely, given a finite set 
$\{p_1,\ldots, p_k\}$ of primes, let $p_{k+1}$ be a 
prime factor of $1+p_1\cdots p_k$. Then, as Euclid showed, 
$p_{k+1}$ is necessarily distinct from $p_1,\ldots,p_k$. Iterating 
this procedure, we thus obtain an infinite sequence of distinct 
primes. For instance, beginning with $k=0$ (with the convention that 
the empty product is 1) and choosing $p_{k+1}$ as small as possible
at each step, one obtains the {\it Euclid-Mullin sequence}
given as the Sloane's sequence A000945   in \cite{sl} 
(cf. Remarks on pages 5 and 11 of this paper). Following 
\cite{boo}, any sequence resulting from this construction is called
a {\it generalized Euclid sequence with seed} $\{p_1,\ldots, p_k\}$ (for 
such a particular sequence, see the sequence A167604 in \cite{sl};
for related sequences, see \cite{boo2} and \cite{booi}).
 More precisely,  Booker in \cite{boo} considered 
a generalization of Euclid's construction described as follows.
If $\{p_1,\ldots, p_k\}$ is a set  of primes, then for any 
$I\subseteq \{1,\ldots,k\}$, the number 
$N_I:=\prod_{i\in I}p_i+ \prod_{i\in \{1,\ldots,k\}\setminus I} p_i$ is 
coprime to $p_1\cdots p_k$ and has at least one prime factor. 
Iteratively choosing a set $I$ and a prime $p_{k+1}\mid N_I$,
we obtain an infinite sequence $p_1,p_2,\ldots$ of distinct 
primes, as in Euclid's proof. It was proved in \cite[Theorem 1]{boo}
that for any finite set $P$ of primes, there is a generalized 
Euclid sequence with seed $P$ containing every prime.            
Notice that in 2016 A. R. Booker and S. A. Irvine \cite{booi} introduced 
the so-called the {\it Euclid-Mullin graph} which encodes all
instances of Euclids's proof of $IP$.

{\bf 7)} In 2016  P. L. Clark \cite{cl2} recast Euclid's proof of $IP$ as 
a {\it Euclidean Criterion} for a domain to have infinitely many 
atoms. It is showed that there is a connection with Furstenberg's 
topological proof of $IP$ (see Subsection 2.7 of this article, p. 25) and
that the presented criterion applies even in certain domains in 
which not all nonzero nonunits factor into products of irreducibles.

{\bf 8), 9)} In 2017  A. Sadhukhan \cite{sad} introduced a partition of the positive 
integers and used it to give two proofs of the infinitude of primes. 
The first proof is a slight variant of the various known  combinatorial 
 proofs. The second is similsr to Euler's proof but it makes no use 
of Euler's product formula.

{\bf 10), 11)} In 2017  S.-I. Seki \cite{sek} gave two proofs of $IP$ via 
valuation theory and gave a new proof of the divergence of the  sum of prime 
reciprocals by Roth'stheorem and Euler-Legendre's theorem for arithmetic 
 progressions.

{\bf 12), 13)} In 2017  S. Northshield \cite{nor2} presented two new proofs 
of $IP$. The first proof uses the 
basic idea of  Furstenberg's celebrated topological proof of $IP$ 
(see Subsection 2.7 of this article, p. 24) but without using topology.
Namely, while Furstenberg's proof is in terms of topological 
space, this proof is in terms of the continuous functions on the space.
 The second proof in \cite{nor2}  uses probability theory.
Namely, this proof is built on the difficulty of defininig a random 
integer.

{\bf 14)} Finally, in 2017 the author of this article in the short note 
\cite{me5} supposed that $\{p_1,p_2,p_3,\ldots, p_k \}$ is a set of all primes with $p_1=2$.
Then by considering the  set of all positive integers that are 
relatively prime to the product $p_2p_3\cdots p_k$, we easily 
obtain a contradiction which implies $IP$.


\subsection{18 recent proofs of $IP$ (2018--2022)}

{\bf 1)} In 2018 S. Silwal \cite[Theorem 1]{sil} proved that  the following 
 inequality holds for sufficiently large $n$:
     $$
\sum_{p\le n}\frac{1}{\log p}>\frac{1}{3}\log n,
     $$ 
where the summation ranges over all primes $p$ such that $p\le n$.
Clearly, the above inequality implies $IP$.

{\bf 2)} In 2018 K. Saito \cite{sai} gives a short  proof of $IP$ by using 
the upper box dimension, which is one of fractal dimensions.

{\bf 3)} In 2020 V.J.-Vera and C.S.-\'{A}vila \cite[Theorem 2]{veav}
gave a new proof of the divergence of the sum 
of the reciprocals of primes using the number of distinct 
prime divisors of a positive integer $n$, and the 
placement of lattice points on a hyperbola given by $n=pr$
with a prime $p$. This immediately yields $IP$.

{\bf 4)} By applying a geometric approach, in 2020 H. G\"{o}ral
\cite{gor} provided  a proof of 
$IP$ via $p$-adic metrics. Notice that 
this is a novel approach to a well-known and quite old result.

{\bf 5), 6, 7)} In 2020 H. G\"{o}ral and H.B. \"{O}czan \cite{goro}
 provided three proofs of $IP$ by considering the properties 
of the Jacobson radical of the ring of integers $\Bbb Z$.
In all of these proofs, the authors supposed that the 
 set $\Bbb P$ of  all primes is finite, i.e.,
$\Bbb P=\{p_1,p_2,\ldots,p_n\}$. 
Let $P=p_1p_2\cdots p_n$ be the product of all primes.
Then by considering the sum $aP^2+P$ with a fixed integer $a\ge 0$,
it is proved in  \cite[Theorem 2.2]{goro} that the fundamental theorem of
 arithmetic implies that $aP^2+P=P$, i.e., $a=0$. This contradiction 
implies $IP$.

Recall that the Jacobson radical of a commutative ring $R$, denoted 
by $J(R)$, is defined as the intersection of all maximal ideals of $R$.
Since all maximal ideals of the ring 
of integers are of the form $p\Bbb Z$ for a prime $p$, the Jacobson 
radical $J(\Bbb Z)$ is the intersection of maximal ideals
 $p_1\Bbb Z,\ldots,p_n\Bbb Z$, and hence,  
       $$
J(\Bbb Z)=\prod_{i=1}^{n}p_i \Bbb Z=(p_1p_2\cdots p_n)\Bbb Z.
 $$
The second proof  and the third  proof of $IP$ given in \cite[Theorem 2.2]{goro}
are based  on consideration of Jacobson radical $J(\Bbb Z)$. The authors 
noticed that there are similarities between Furstenberg's 
topological proof of $IP$ \cite{f} and  their second proof 
 and third proof of $IP$.

{\bf 8)} In 2020 F. Lemmermeyer \cite{lem} provided a short simple 
proof of $IP$. This proof is a simplification of the proof of $IP$
using continued fractions given in 1976 by Barnes \cite{bar}. Asssume 
that there are only infinitely many primes, namely, 2, $p_1=3,\ldots,p_n$,
Let $q=p_1\cdots p_n$ be the product of all odd primes. Then 
$q^2+1$ is not divisible by any odd prime, hence must be a power of two.
Since $q^2+1\equiv 2\pmod{4}$, must be $q^2+1=2$ and therefore, 
$q=1$, which is a contradiction. Since no odd prime 
$p\equiv  3\pmod{4}$ can divide $q^2+1=2$, the proof 
actually shows that there are infinitely many primes
 $p\equiv 1\pmod{4}$.

{\bf 9), 10)} By considering the notion of a realizable integer 
sequence, in 2020   P. Moss and T. Ward \cite[Lemma 1]{mow} 
proved that the set of primes dividing a denominator of 
$\frac{1}{n}\sum_{d\mid n}\mu\left(\frac{n}{d}\right)f_d$ for some 
positive integer $n$ is infinite, where $\mu(n)$ is the classical
 M\"{o}bius 
arithmetic function and $(f_n)$ is the well known  Fibonacci's sequence 
defined by the conditions $f_1=f_2=1$, $f_{n+2}=f_{n+1}+f_n$ 
with $n=1,2,\ldots$. This implies $IP$.

 Furthermore, P. Moss and T. Ward \cite[Corollary 3]{mow} 
proved that if $j$ is an arbitrary odd positive integer,
then the set of primes dividing  denominators of 
$\frac{1}{n}\sum_{d\mid n}\mu\left(\frac{n}{d}\right)f_{d^j}$. 
 This implies $IP$.

{\bf 11)} Let $p_1,p_2,\ldots,p_{\pi(n)}$ be all the primes less than 
or equal to $n$. Using the inclusion-exclusion principle,  in 2020 S. Laad   
\cite{laa} proved the inequality
  $$
n\prod_{i=1}^{\pi(n)}\left(1-\frac{1}{p_i}\right)<1+2^{\pi(n)-1}.
  $$
Assuming that there are only $k$ primes, 
then clearly, the left hand side of the above inequality is unbounded,
while the right hand side is a constant. This contradiction impies $IP$.  


{\bf 12), 13, 14)} In 2021 C. Elsholtz \cite[Theorem 1]{els} showed that 
Fermat's last theorem 
and a combinatorial theorem of Schur 
 on monochromatic solutions of $a+b=c$
(Lemma 1 in \cite{els}) implies $IP$. In particular, since 
there exist elementary proofs of Fermat's last theorem 
for $n=3$ $n=4$ and $n=5$ (concerning the Diophantine equations
in positive integers $x^3+y^3=z^3$, $x^4+y^4=z^4$ and $x^5+y^5=z^5$, 
respectively; (see \cite{r3}), Theorem 1 in \cite{els} implies the elementary 
proof of $IP$. 

It follows from Theorem 2 of \cite{els} that Roth's theorem 
(Lemma 2 in \cite{els})  implies $IP$.

It was also proved in \cite[Theorem 3]{els} that Hindman's theorem 
(Lemma 4 in \cite{els}) implies $IP$.

{\bf 15)} In 2021 L. Haddad \cite{hadd} simplified the above 
mentioned proof of $IP$ due to C. Elsholtz \cite[Theorem 1]{els}. 
Namely, this proof is greatly simplified, making no use at all Fermat's last 
theorem, and using only a weak form of the theorem of Schur   
 on monochromatic solutions of $a+b=c$.

{\bf 16)} In 2022 J. Mehta \cite[Theorem 1]{meh}  generalized 
M\'{e}trod's proof of 
$IP$ given in 1917 \cite{met} (also see \cite[p. 11]{r2}). Assume that 
$p_1,\ldots,p_n$ are distinct primes whose product is $P$, and choose 
any factorization of $P$ into $k\ge 2$ termse, say $P=d_1\cdots d_k$, and 
put \cite[Proof of Theorem 1]{meh}  
   $$
M=\frac{P}{d_1}+\cdots + \frac{P}{d_k}.
  $$
Then it is easy to show that there exists a  prime 
$p\not\in\{p_1,\ldots,p_n\}$ dividing $M$. This implies $IP$.  

         Notice  that Stieltjes' proof of $IP$ given in his work in 1890  
(\cite[p. 14]{st}; also see \cite[p. 4]{na}) is a particular case of Theorem 1 of 
\cite{meh} with $k=2$. Furthermore, M\'{e}trod's proof of $IP$ is recovered 
by taking $k=n$ and $d_i$ for $i=1,2,\ldots,n$, i.e., 
by considering the divisors of $M=p_1p_2\cdots p_n\left(1/p_1+ 1/p_2+\cdots 
+1/p_n\right)$.

{\bf 17)} Using prime factorization theorem of a positive integer,  in  2022   
R. Me\v{s}trovi\'{c} \cite{me9} gave a short proof by contradiction of $IP$
(this is in  fact the proof in Section 4 of this article, pp. 40--41).

{\bf 18)} Using M\"{o}bius inversion formula \cite{nzm},  in  2023   
R. Me\v{s}trovi\'{c} \cite{me10} gave a very short proof of the formula 
due in 2009  to J. Pinasco \cite{p} which is applied in his proof of $IP$.
Consequently, using a simpler argument than those of Pinasco's proof, 
it follows $IP$.

\section{Proofs of $IP$ in arithmetic progressions: special cases of 
Dirichlet's theorem} 

\subsection{Dirichlet's  theorem}
In 1775 L. Euler \cite{eu4} 
(also cf. \cite[p. 415]{d}, \cite[p. 108, Section 3.6]{ta}) 
stated  that an arithmetic progression with the first term equals 1 
and the difference $a$ to be a positive integer,  contains 
infinitely many primes.  More generally,  in 1798 
in the second edition of his book A.-M. Legendre \cite{leg} 
(cf. \cite[p. 415]{d} and \cite[p. 108, Section 3.6]{ta}) conjectured that 
for relatively prime positive integers $a$ and $m$ there are
infinitely many primes which leave a remainder of $m$ when 
divided by $a$. In other words, if $a$ and $m$ are relatively
prime positive integers,  then the arithmetic progression $a, a+m,a+2m,
a+3m,\ldots$ contains infinitely many primes. The condition 
that $a$ and $m$ are relatively prime is essential, for otherwise there would 
be no primes at all in the progression. However, Legendre gave a 
proof that was faulty.
In 1837 Peter Gustav Lejeune Dirichlet, Gauss's successor of G\"{o}ttingen 
and father of analytic number theory, gave a correct proof.
Namely,  Dirichlet \cite{dir}  
proved the following theorem  which is 
a far-reaching extension of Euclid's theorem on the infinitude of primes 
and is one of the most beautiful results in all of Number Theory. 
It can be stated as follows.
\vspace{2mm}

\noindent{\bf Dirichlet's  theorem.}  {\it Suppose $a$ and $m$ are relatively 
prime positive integers. Then there are infinitely many primes of the form 
$mk+a$ with $k\in\Bbb N\cup\{0\}$.}

\vspace{2mm}

Dirichlet's proof is derived by means of $L$-{\it functions} and analysis.
The main strategy is, as in Euler's proof of $IP$ 
(which in fact shows that the sum of reciprocals of primes diverges),
 to consider the function 
 $$
P_m(s):=\sum_{p\equiv a(\bmod{\, m})}\frac{1}{p^s},
  $$
(where the sum is only over those primes $p$ 
that are congruent to $a(\bmod{\, m})$)
which is defined say for real numbers $s>1$, and to show that
$\lim_{s\to 1+}P_m(s)=+\infty$. Of course this suffices, because 
a divergent series must have infinitely many terms.
The function $P_m(s)$ will in turn be related to a finite 
linear combination of logarithms of {\it Dirichlet $L$-series},
and the differing behavior of the Dirichlet series for 
principal and non-principal characters is a key aspect of the proof.
Dirichlet used an ingenious argument to show that the sum 
$\sum_{p\equiv a(\bmod{\, m})}1/p$ diverges, where the 
sum ranges  over all primes $p$ that are congruent to $a(\bmod{\, m})$.

\vspace{2mm}

\noindent{\bf Remarks.}  
As it is pointed out by P. Pollack \cite{pol3}, 
there exist proofs of Dirichlet's theorem which minimize 
analytic prerequisities (e.g., those of A. Selberg \cite{se}
in 1949, A. Granville \cite{gr} in 1989
and  H. N. Shapiro (\cite{sh1} and \cite{sh2}) in 1950). For example, 
Selberg \cite{se} gave a proof that is, he wrote
``more elementary in the respect that we do not use the complex 
characters mod $k$, and also in that we consider only finite sums." 
An ``elementary proof" of Dirichlet's theorem 
 in the sense that it does not use complex analysis is given 
by  M. B. Nathanson \cite[Ch. 10]{nat}.
Nevertheless, all these ``elementary" proofs exhibit at
least as complicated a structure as Dirichlet's original argument.  
This  is well discussed and considered  in 2010 by A. Granville  in his 
expository article \cite[Sections 2 and 3]{gr3}.\hfill$\Box$
\vspace{2mm}


\subsection{A survey of elementary proofs of  $IP$ in  
special arithmetic progressions}
For many arithmetic progressions  with small differences one can obtain
simple elementary (i.e. not using analytic means) proofs of Dirichlet's 
theorem. Several of them are listed by Dickson 
\cite[pp. 418--420, Chapter XVIII]{d} 
and Narkiewicz \cite[pp. 87--96, Section 2.5]{na}.
In \cite{mt} M. R. Murty and N. Thain asked 
``how far Euclid's proof can be pushed to yield Dirichlet's theorem''.
 The existence of such a ``{\it Euclidean proof}'' 
(precised in \cite{mt}) for certain arithmetic progressions is well 
known. For example, considering the product $k(2\cdot 3\cdots p_n)$,
Euclid's elementary proof can be used to prove that for any fixed  
positive integer $k>2$ there are infinitely many primes which are not 
congruent to $1(\bmod{\,k})$.  
This result was proved in 1911 by H. C. Pocklington \cite{po}
(also see \cite[p. 116, Theorem 114]{cl} and \cite[p. 419]{d}).

 Further,  we expose other proofs of $IP$ in special arithmetic 
progressions of the form $1(\bmod{\,k})$ and $-1(\bmod{\,k})$. 
An excellent source for this subject is Narkiewicz's monograph 
\cite[pp. 87--93, Section 2.5]{na}. 
An elementary proof of $IP$ in every progression $1(\bmod{\,2p})$, 
where $p$ is any prime, was established in 1843  by  V. A. Lebesgue 
(\cite[p. 51]{le2}, \cite[p. 418]{d}) who showed the fact 
that $x^{p-1}-x^{p-2}y+\cdots +y^{p-1}$ has besides the possible factor $p$ 
only prime factors of the form $2kp+1$ ($k=1,2,\ldots$). 
Using a quite similar method, in 1853 F. Landry (\cite{lan}, \cite[p. 418]{d}) 
considered prime divisors of $(n^p+1)/(n+1)$ 
to prove $IP$ for the same progressions. This proof can be found in 
\cite[p. 121, Ch. 24, Exercise 24.1]{aru}. 
By a quite similar method, the same result can 
be  obtained using the fact that for any prime $q$
 every prime divisor $p$ of $(n^q-1)/(n-1)$ coprime with $q$ satisfies
$p\equiv 1(\bmod{\,q})$ (see e.g.  \cite[p. 34, Section 2.3]{ik}
or \cite[pp. 151--152, Problem 7.3.3]{ana}).
The analogous  method is also applied by Lebesgue in 1862 (\cite{le4}, 
\cite[p. 418]{d}) 
for the progression $-1(\bmod{\,2p})$ with a prime $p$.  
Using the rational and irrational parts of $(a+\sqrt{b})^k$,
in  1868/9 A. Genocchi (\cite{gen}, \cite[p. 418]{d}) 
proved $IP$ in both progressions $1(\bmod{\,2p})$ and $-1(\bmod{\,2p})$, 
where $p$ is an arbitrary prime. Furthermore, in lectures of 1875/6
L. Kronecker  (\cite{kro}, \cite[pp. 440--442]{has}) 
gave another proof of $IP$ 
in the progression $1(\bmod{\,2p})$ with a prime $p$.
Another simple proof of the same result  based on 
Euler's totient function and Fermat little theorem 
is  recently given in \cite{me3}.

Using the fact that $(2^{mp}-1)/(2^m-1)$ 
($p$ a prime and $m$ a positive integer) has at least one prime divisor of 
the form $p^nk+1$ (\cite[p. 107, proof of Theorem 47]{sha}; 
also cf. \cite[pp. 178--179, Theorem 11]{ers} or 
\cite[p. 209, Exercise 1.5.28]{mur2}), 
in 1978 D. Shanks \cite{sha} proved that for every prime power $p^n$ there 
are infinitely primes  $\equiv 1(\bmod{\,p^n})$.
Another  elementary proof of $IP$ in the progression $1(\bmod{\,p^n})$ for 
any prime $p$ and $n=1,2,\ldots$ was given in 1931 by F. Hartmann \cite{hart}.

Using divisibility properties of cyclotomic polynomials,  
in 1888  J. J. Sylvester \cite{sy} proved $IP$ in the progressions 
$-1(\bmod{\,p^n})$, where $p^n$ is any prime power.   
In 1896 R. D. von Sterneck \cite[p. 46]{ste} (cf. \cite[p. 90]{na}) 
considered a product $F(n):=\prod_{d\mid n}f(n/d)^{\mu(d)}$, where
$\mu$ is the  M\"{o}bius function,  
$f(n)$ is an integer-valued function satisfying $f(1)=1$ and 
two divisibility  properties. Then every prime dividing $F(n)$ divides 
$f(n)$ but does not divide $f(i)$ for each $i=1,2,\ldots, n-1$. 
   Von Sterneck remarked that a recursive sequence $f(n)$ defined as 
$f(n)=f(n-1)+cf(n-2)$ with $f(1)=1$ and a positive integer $c$, satisfies 
these conditions, and used this it can be obtained an elementary proof of    
infinitely many primes  $\equiv -1(\bmod{\,p^m})$ for any fixed  
prime power $p^m$. 
The same result for powers of odd primes 
and the infinitude of primes $\equiv -1(\bmod{\,3\cdot 2^n})$ 
were  proved in 1913 by R. D. Carmichael \cite{car}.

As remarked by Dickson \cite[p. 418]{d}, using  cyclotomic polynomials
$\Phi_m(x)$, in 1886 A. S. Bang (\cite{ba}, \cite[p. 418]{d})  and in 1888 
Sylvester (\cite{sy}, also cf. \cite[p. 418]{d})
obtained proofs of $IP$ in arithmetic progressions $1(\bmod{\,k})$, 
where  $k$ is any integer $\ge 2$. Both these proofs  are based 
on the fact that if $p$ is a prime not dividing $m$, then $p$ divides 
$\Phi_m(a)$ if and only if the order of $a(\bmod{\,}p)$ is $m$.
(Here  $\Phi_m(x)$ is the $m$th cyclotomic polynomial).
Such a simple classical proof 
of $IP$ in arithmetic progressions $1(\bmod{\,k})$ which is in spirit 
``Euclidean'' can be found in (\cite{gt} and \cite[pp. 116--117]{cl}; 
also cf. \cite[pp. 97--99]{ke} and \cite[pp. 12--13]{wa2}). 
Considering the least common multiple of polynomials 
$\{x^d-1:\, d\mid n\}$, in 1895 E. Wendt \cite{wen} (cf. \cite[p. 89]{na}) 
gave a simple proof of the same result. 
Moreover, Narkiewicz \cite[p. 88]{na} noticed that, according 
to a theorem of Kummer \cite{k2} (also see \cite[Theorem 4.16]{na2}),
a rational prime $p$ splits in the $k$th {\it cyclotomic field} 
$\Bbb Q(\zeta_k)$ (where $\zeta_k$ denotes a primitive $k$th 
root of unity) if and only if it is congruent to  $1(\bmod{\,k})$.
Using this and the fact that in any given finite extension of $\Bbb Q$ there 
are infinitely many splitting primes, we obtain $IP$ in every arithmetic 
progression $1(\bmod{\,k})$ with $k\ge 2$.
Studying the existence of {\it primitive prime divisors} 
of integers $a^n-b^n$, where $n\in\Bbb N$ and $a$ and $b$ 
are relatively prime integers, in 1903/04 G. D. Birkhoff and H. S. Vandiver
\cite{bv} gave an elementary proof of this result.
A variation of this proof has been given in 1961 by A. Rotkiewicz \cite{rot},
whose proof was simplified in 1962/3 by T. Estermann \cite{es}
and in 1976 by I. Niven and B. Powell \cite{np}.
In their proof Niven and  Powell use only elementary divisibility 
properties and the fact that the number of roots of a non-zero
polynomial cannot exceed its degree. Applying 
{\it Birkhoff-Vandiver theorem} (see e.g., \cite[p. 88]{na}), 
the same result was proved in 1981 by R. A. Smith \cite{sm} 
(see also \cite[Chapter 1]{na3} and \cite[pp. 88--89]{na}).
Another two elementary proofs were given in 1984 by S. Srinivasan \cite{sr}
and in 1998 by N. Sedrakian and J. Steinig \cite{ss}.   
An elementary proof of this assertion was provided in 2004 by 
J. Yoo \cite{yo} without using cyclotomic polynomials.
Another two old proofs of this result are due  to K. Th. Vahlen 
\cite{va} in 1897 by using {\it Gauss' periods of roots of unity}
and \'{E}. Lucas \cite[p. 291, Ch. XVII]{lu3} in 1899 applying his (Lucas)
sequence $u_n$. 

A short but not quite elementary proof of $IP$ in 
the progressions $-1(\bmod{\,k})$ for each $k\ge 2$
was given by M. Bauer  \cite{bau} in 1905/6.
In 1951 T. Nagell \cite[pp. 170--173]{nag} gives an elementary  proof of
$IP$ in arithmetic progression  $-1(\bmod{\,k})$ with $k\ge 2$.

Applying a similar argument to those of Niven and Powell 
for $IP$ in the progressions $\equiv 1(\bmod{\,k})$,
in 1950 by M. Hasse \cite{has} proved $IP$ in the 
progressions $-1(\bmod{\,k})$ for each $k\ge 2$.

Euclidean's proofs of $IP$ in  various arithmetic progressions can be
found in Problems book of Murty and Esmonde \cite[Section 7.5]{mes} in 2005.
For example,  the known facts  that every prime divisor of the Fermat number
$F_n:=2^{2^n}+1$ is of the form $2^{n+1}k+1$ 
(see e.g., \cite[p. 8, Exercise 1.2.8]{mes}) and that $F_n$ and $F_m$ are 
relatively prime if $m\not=n$ (see  Subsection 2.2) yield that 
there are infinitely primes $\equiv 1(\bmod{\,2^n})$ for any given 
$n$ (\cite[p. 11, Exercise 1.4.13]{mes}, also cf. \cite[p. 151, 
Problem 7.3.2]{ana}).

As noticed by K. Conrad \cite{con}, a Euclidean proof of 
Dirichlet's theorem for $m({\bmod\, a})$ involves, at 
the very least, the construction of a nonconstant polynomial 
$h(T)\in\Bbb Z[T]$ for which any prime factor 
$p$ of any integer $h(n)$ satisfies, with finitely many exceptions,
either $p\equiv 1(\bmod{\, a})$  or $p\equiv m(\bmod{\, a})$,
and infinitely many primes of the latter  type occur.
For example \cite{con}, Euclidean proofs of 
Dirichlet's theorem exist for arithmetic progressions 
$1(\bmod{\, a})$ with any $a\ge 2$, $3(\bmod{\, 8})$,
$4(\bmod{\, 5})$ and $6(\bmod{\, 7})$.

A characterization of  arithmetic 
progressions for which Euclidean proof exist 
is given by I. Schur \cite{sch} and M. R. Murty \cite{mur}.
In 1912/13 I. Schur \cite{sch} proved that if $m^2\equiv 1(\bmod{\, a})$,
then a Euclidean proof of Dirichlet's theorem exists for the
arithmetic progression $m(\bmod{\, a})$. 
In particular, Schur extended Serret's approach 
based on law  of quadratic reciprocity to establish
 proofs  of $IP$ for the progressions $2^{m-1}+1(\bmod{\, 2^m})$,
$2^{m-1}-1(\bmod{\, 2^m})$ ($m\ge 1$), and $l(\bmod{\, k})$ for 
$k=8m$ (with $m$ being an odd positive squarefree integer) and 
$l=2m+1$, $l=4m+1$ or $l=6m+1$ (cf. \cite[p. 91]{na}. 
A similar method was used in 1937 by A. S. Bang \cite{ba2}
(cf. \cite[p. 91]{na} who proved $IP$ in the progressions 
$2p^m+1(\bmod{\, 4p^m})$ with prime $p\equiv 3(\bmod{\, 4})$,
$2p^{2n+1}+1 (\bmod{\, 6p^{2n+1}})$ with prime $p\equiv 2(\bmod{\, 3})$,
and $4p^{2n}+1(\bmod{\, 6p^{2n}})$ with prime $p\equiv 2(\bmod{\, 3})$.

\vspace{2mm}

\noindent{\bf Remarks.}  
In 1988 Murty  (\cite{mur}; also see \cite{mt}) proved the 
converse of Schur's result, i.e., he showed 
that a Euclidean proof exists for the arithmetic progression 
$m(\bmod{\, a})$  only if $m^2\equiv 1(\bmod{\, a})$.
This means that it is impossible to prove Dirichlet's theorem for
certain arithmetic progression by Euclid's method.
The proof due to Murty is not difficult, but involves some {\it Galois Theory}.
For example, since $2^2\equiv 4\not\equiv 1(\bmod{\, 5})$,
there is no proof of Dirichlet's theorem for $2(\bmod{\, 5})$ 
which can mimic Euclid's proof of $IP$.  Notice also that 
 Dirichlet's theorem can be proved by Euclidean's methods for all
the possibilities modulo $a=24$ (cf. \cite{bl}).  
Recently, P. Pollack \cite{pol3} discussed  Murty's definition
 of a ``Euclidean proof" and Murty's converse of Schur's result.  
Finally, we point out an interesting expository article of A. Granville
\cite{gr4} in 2007 in which are compared numbers of primes in different
arithmetic progressions with the same small difference.\hfill$\Box$

\subsection{Elementary proofs of $IP$ 
in arithmetic progressions with small differences}
In this subsection,
we expose several Euclidean proofs of $IP$ in different arithmetic 
progressions with small differences. 
Dickson's {\it History} records several further attempts at giving 
Euclidean proofs for particular progressions 
(see the listing on \cite[pp. 418--420]{d}).
Considering the product $2^2\cdot 3\cdot 5\cdot p_n-1$,  
 Euclid's idea is used by V. A. Lebesgue \cite{le} in 1856 
(also cf. \cite[p. 13, Theorem 11]{hw})
for the progression $3(\bmod{\,4})$.
A. Granville \cite[p. 3, Section 1.3]{gr2} remarked that 
a similar proof works for primes $\equiv 2(\bmod{\,3})$. 
The same idea  that involves the product  $2\cdot 3\cdot 5\cdot p_n-1$
was also used by V. A. Lebesgue in 1859  
(\cite{le3}; also see  \cite[p. 13, Theorem 13]{hw} and \cite[p. 419]{d})
for the proof of $IP$ in the progressions  $5(\bmod{\,6})$ and 
$1(\bmod{\,2^n})$ with a fixed $n=1,2,\ldots$.   
The situation is more complicated 
for the progression $1(\bmod{\,4})$ and related proof 
is based on the consideration of the product 
$N:=(5\cdot 13\cdot 17\cdots p_n)^2+1$
and the fact that if integers $a$ and $b$ have no common factor,
then any odd prime divisor of $a^2+b^2$ is congruent to 
$1(\bmod{\,4})$ \cite[Theorem 13]{hw}. In fact,
using this property of {\it quadratic residues}
and Euclid's idea, Hardy and Wright proved in his book
 (\cite[Theorem 14]{hw} which was first written in 1938)
that the progression $5(\bmod{\,8})$ contains infinitely many primes. 
Dickson \cite[p. 419]{d}
noticed that this result and proofs of $IP$ in  progressions 
$1(\bmod{\,8})$, $3(\bmod{\,8})$ and $7(\bmod{\,8})$
were firstly proved in 1856 also by A. V. Lebesgue \cite{le4}. 
Using some properties of Fermat and Fibonacci numbers,
 two constructive proofs of $IP$ in progression $1(\bmod{\,4})$ 
were presented in 1994 by N. Robbins \cite{ro}.  

Dickson \cite[p. 419]{d}  pointed out the  proofs 
of $IP$ also in the following arithmetic progressions: 
$9(\bmod{\,10})$ due to J. A. Serret \cite{ser} in 1852,
$2(\bmod{\,5})$ and $7(\bmod{\,8})$ due to  \'{E}. Lucas \cite[p. 309]{lu} in 
1878, $1(\bmod{\,4})$, $5(\bmod{\,6})$ and 
$5(\bmod{\,8})$  due to  \'{E}. Lucas \cite[pp. 353--354]{lu2} in 1891, 
 $1(\bmod{\,4})$, $1(\bmod{\,6})$ and $5(\bmod{\,8})$  due to
E. Cahen  \cite[pp. 318--319]{cah} in 1900
and also \cite[pp. 438--439]{has} in 1875/6, $1(\bmod{\,4})$, 
$1(\bmod{\,6})$,  $3(\bmod{\,8})$, $7(\bmod{\,8})$, $9(\bmod{\,10})$ and
  $11(\bmod{\,12})$,  due to K. Hensel \cite{has} in 1913.
Furthermore, using law  of quadratic reciprocity
\cite[Sections 112--114]{ga}, in 1852 J. A. Serret 
\cite{ser} (also cf. \cite[pp. 90--91, Theorem 2.19]{na}) proved $IP$ 
in the progressions $3(\bmod{\,8})$, $5(\bmod{\,8})$ and $7(\bmod{\,8})$.

Considering divisors of integer $(11\cdot 31\cdot 41\cdot 61\cdots p_n)^5-1$,
it was proved in 1962 \cite[pages 60, 371--373, Problem 254(c)]{scy}
$IP$ in the  progression $1(\bmod{\,10})$. The analogous idea 
was used in 2007 by A. Granville \cite[p. 4]{gr2}  to show $IP$ 
in the progression $1(\bmod{\,3})$. 
 
There are also elementary arguments in spirit of Euclid's idea 
showing that there are infinitely many primes in other arithmetic progressions
with small differences, such as $4(\bmod{\,5})$, $1(\bmod{\,8})$
and $3(\bmod{\,8})$. In 1965 P. Bateman and M. E. Low \cite{bl} give  a 
proof similar to Euclid's that  for every coprime residue class 
$a(\bmod{\, 24})$ there are infinitely many primes in  
progression $a(\bmod{\, 24})$.
Their proof makes use of the interesting fact that every integer 
$a$ relatively prime to 24 has the property $a^2\equiv 1(\bmod{\, 24})$.
Using a couple of observations about the polynomial 
$f(x)=x^4-x^3+2x^2+x+1$ and the law  of quadratic reciprocity,
a Euclid-type proof for the progression $4(\bmod{\,15})$ 
is presented in 2005 by M. R. Murty and J. 
Esmonde \cite[pp. 92--64, Example 7.5.4]{mes}.

Considering the linear second order recurrence $u_n=u_{n-1}+3u_{n-2}$ with 
$u_0=1$, $u_1=1$, in 2005 R. Neville \cite{n2} 
gave a simple proof of $IP$ in progression $1\bmod{\,3}$.
The author \cite[Remarks]{n2} also noticed that 
if $q\ge 5$ is a given prime, then considering 
the {\it Lucas sequence} $u_n=u_{n-1}+3u_{n-3}$ with $u_0=0$, $u_1=1$,
similarly one can  prove that there are infinitely many 
primes $p$ such that $\left(\frac{-q}{p}\right)=1$
($\left(\frac{\cdot}{p}\right)$ denotes the {\it Legendre symbol}).
In particular, for $q=5$ this yields $IP$ in 
all progressions $a(\bmod{\,20})$ with $a\in\{1,3,7,9\}$.
In book \cite[p. 15]{r3} P. Ribenboim noticed that 
in 1958 D. Jarden \cite{jar} proved $IP$ in the progression 
$1(\bmod{\,20})$.
 
\section{Another simple Euclidean's proof of Euclid's theorem}

\begin{proof}[Proof of Euclid's theorem]
Suppose that $p_1=2<p_2=3<\cdots <p_k$ are all the primes. Take
$n=p_1p_2\cdots p_k+1$ and let $p$ be a prime dividing $n$. 

The first step is a ``shifted" first step of Euclid's proof.
Suppose that $p_1=2<p_2=3<\cdots <p_k$ are all the primes. Take
$n=p_1p_2\cdots p_k$. 
Then $n-1=p_1^{e_1}p_{2}^{e_2}\cdots p_{k}^{e_k}(\ge 5)$
for some $k$-tuple of nonnegative integers $(e_1,e_2,\ldots,e_k)$,
and so taking $s=\max\{e_1,e_2,\ldots,e_k\}$, we find that
   $$
n-1=p_1^{e_1}p_{2}^{e_2}\cdots p_{k}^{e_k}=
\frac{p_1^{s}p_{2}^{s}\cdots p_{k}^{s}}
{p_1^{s-e_1}p_{2}^{s-e_2}\cdots p_{k}^{s-e_k}}=\frac{n^s}{a},
  $$
where $a=p_1^{s-e_1}p_{2}^{s-e_2}\cdots p_{k}^{s-e_k}$ and $s$ 
are positive integers. The above equality yields
$$
a=\frac{n^s}{n-1}= \frac{(n^s-1)+1}{n-1}=
\sum_{i=0}^{s-1}n^i+\frac{1}{n-1},
  $$
whence it follows that $1/(n-1)=a-\sum_{i=0}^{s-1}n^i$ 
is a positive integer. This contradicts the fact that $n-1\ge 4$,
and the proof is completed. 
  \end{proof}

\noindent{\bf Remarks.} Unlike most other proofs of the Euclid's
theorem, Euclid's proof and our proof  
does not require Proposition 30 in Book VII of {\it Elements} (see 
\cite{z}, \cite{hw},  where this result
is called {\it Euclid's first theorem}; sometimes called 
``{\it Euclid's Lemma}") 
that states into modern language from the Greek \cite{he}: 
{\it that if two numbers, multiplied by one another make some number,
and any prime number measures the product, then it also 
measures one of the original numbers},
or in terms of modern {\it Arithmetic}:
if $p$ is a prime such that $p\mid ab$ then either $p\mid m$ 
or $p|b$. It was also pointed in \cite[page 10, Notes on Chapter 1]{hw}
that this result does not seem to have 
been stated explicitly before Gauss  of 1801 who gave the 
first correct proof of this assertion \cite[Sections 13--14]{ga}. 
The only divisibility property  used in our  proof and Euclid's  proof
is the fact that  every integer $n>1$ has at least
one representation as a product of primes.
This is in fact, Proposition 31 in Book VII of {\it Elements}
(see above Remarks).
 
In order to achieved a contradiction, in the second step of his proof
Euclid take a prime  that divides a product $P$ of
all the primes plus one, and further consider 
two cases in dependence  on whether $P$ is prime or not. 
But in the second step of our  proof 
we obtain directly a contradiction dividing $n^s$ by $n-1$.\hfill$\Box$

\vspace{5mm}

%
\centerline{A{\scriptsize{PPENDIXES}}}

\vspace{5mm}

\centerline{\bf{A) External Links on Euclid's theorem and 
its proofs}}

\vspace{3mm}

Wikipedia {\tt http://en.wikipedia.org/wiki/Euclid's$\_$theorem}

\noindent{\tt http://mathworld.wolfram.com/EuclidsTheorems.html}, 
from MathWorld. 

{\tt http://primes.utm.edu/notes/proofs/infinite/euclids.html}

{\tt http://mathforum.org/}

{\tt http://aleph0.clarku.edu/$\sim$djoyce/java/elements/elements.html}

{\tt http://planetmath.org/encyclopedia/}

{\tt http://mathoverflow.net}

{\tt http://tech.groups.yahoo.com/group/primenumber/}\\

\centerline{\bf{B) Sloane's sequences related to 
proofs of Euclid's theorem}}
\vspace{3mm}

{\small A000040, A002110, A034386, A210144, A210186,
A006862, A005234, A006794, 
A014545, A057704, A057713, A065314, A065315, A065316, A065317, 
A018239, A057588, A057705  A006794, A002584, A002585,
A051342, A068488, A068489, A103514, A066266, A066267, A066268, A066269,
A088054, A093804,  A103319, A104350, A002981, A002982, A038507, 
A007917, A007918, 
A088332, A005235, A000945, A000946, A005265, A005266, A0084598, A0084599, 
A005266; A000215, A019434, A094358,
A050922, A023394, A057755, A080176, A002715; 
A000668, A001348, A000225, A000043, A046051, A028335;
A002716; A104189; A001685; A000045; A000217; A000292; 
A064526, A000324, A007996; A000289;  
A000058, A001543, A001544, A126263; A005267; 
A0013661, A0013662; A003285; A000010; A000984; A167604.}    

 In ``{\it The On-Line Encyclopedia of Integer Sequences}.'' 
(published electronically at {\tt www.research.att.com/$\sim$njas/sequences/})
\cite{sl}.

\vfill\eject

\begin{center}
\bf{C) List of papers and their authors arranged by year of publication
followed by the  main argument(s) of related proof given into round brackets}
\end{center}

\footnote{${}^1$ $*$ denotes that a a related proof of $IP$ concerns 
 a particular arithmetic progression}

For brevity, into round brackets after  a reference in the following list
${}^1$  we denote the method(s) and/or  idea(s) that are used in  related proof by:  

$AP$--an arithmetic progression/arithmetic progressions;

$C$--a combinatorial method;

$CM$--a counting method, based on some combinatorial
enumerating arguments; 

$CS$-an idea  based on a convergence 
of sums $\sum_{n=1}^{\infty}1/n^s$ with $s>1$ etc;

$DS$--Euler idea, that is an idea  based on the divergence 
of reciprocals of primes and related series;

$E$--Euclid's idea of the proof of the infinitude of primes,
that is, a consideration of product $P:=p_1p_2\cdots p_k+1$
or some analogous product;

$FT$--a factorization (not necessarily to be unique) 
of a positive integer as a product of prime powers;

$MPI$--the idea based on a construction of 
sequences consisting of mutually prime positive integers;

$T$--a topological method;

$UFT$--the unique factorization theorem of a positive integer 
as a product of prime powers.  

\vspace{2mm}

{\small

\cite[$\sim$ 300 B.C.]{he}, \cite[p. 4, Theorem 4]{hw},  Euclid of Alexandria 
($E$)

\cite[1730, pp. 32--34, I]{fu}, \cite[p. 6]{r2}, \cite[pp. 40--41]{ew}, 
\cite[p. 4]{pol} C. Goldbach, ($MPI$, especially 
Fermat numbers $F_n:=2^{2^n}+1$)

\cite[1736]{eu3} (posthumuous paper), \cite[pp. 134--135]{bu}, 
 \cite[p. 413]{d}, \cite[p. 3]{pol} L. Euler 
(multiplicativity of Euler's totient function $\varphi$) 

\cite[1737, pp. 172--174]{eu2}, \cite{eu}, \cite[p. 8]{r2}, 
\cite[pp. 8--9]{ew}, 
L. Euler ($UFT$, $DS$; especially, the series of the reciprocals of 
the primes is divergent) 

${}^*$\cite[1843]{le2}, \cite[p. 418]{d}  V. A. Lebesgue (prime factors 
of $x^{p-1}-x^{p-2}y+\cdots +y^{p-1}$, Fermat little theorem, and
 $IP$ in $AP$ $1(\bmod{\,2p})$ with a $p$ a prime)

${}^*$\cite[1852]{ser} J. A. Serret ($E$,  and $IP$ in $AP$
 $9(\bmod{\,10})$) 

${}^*$\cite[1852]{ser}  \cite[pp. 90--91, Theorem 2.19]{na}
J. A. Serret, (law  of quadratic reciprocity, and $IP$ in $AP$ $3(\bmod{\,8})$, 
$5(\bmod{\,8})$ and $7(\bmod{\,8})$)

${}^*$\cite[1853]{lan},  \cite[p. 418]{d}, \cite[p. 34, Section 2.3]{ik}
F. Landry (prime divisors of $(n^p+1)/(n+1)$, Fermat little theorem, and
 $IP$ in $AP$ $1(\bmod{\,2p})$ with a prime $p$)

${}^*$\cite[1856]{le}, \cite[p. 13]{hw} V. A. Lebesgue ($E$, and
$IP$ in $AP$ $3(\bmod{\,4})$)

${}^*$\cite[1859]{le3},  \cite[p. 418]{d}, \cite[p. 13]{hw},
V. A. Lebesgue ($E$, $IP$ in $AP$ $5(\bmod{\,6})$ 
and $IP$ in $AP$ $1(\bmod{\,2^nk})$ with some fixed 
$k,n\in\Bbb N $)

${}^*$\cite[1862]{le4}, \cite[p. 418]{d} V. A. Lebesgue,
(prime factors of an integer  polynomial in two variables, 
Fermat little theorem, and 
 $IP$ in  $AP$ $-1(\bmod{\,2p})$ with a $p$ a prime)

${}^*$\cite[1868/9]{gen},   \cite[p. 418]{d}) A. Genocchi
(rational and irrational parts of $(a+\sqrt{b})^k$, and 
$IP$ in $AP$ $\pm 1(\bmod{\,2p})$, where $p$ is an arbitrary prime)

${}^*$\cite[1871]{sy2} J. J. Sylvester 
(certain identities between rational functions, 
 $IP$  in  $AP$ $3(\bmod{\,4})$ and $5(\bmod{\,6})$)

\cite[1874]{mer}, \cite[pp. 171, 183--186]{yy2} F. Mertens
($DS$, the boundedness of the quantity $|\sum_{p\le n}\log p/p-\log n|$
as $n\to\infty$)

\cite[1875/6, pp. 269--273]{has}, \cite[p. 413]{d} 
L. Kronecker ($DS$ and $CS$)

${}^*$\cite[1875/6, pp. 438--439]{has}  K. Hensel  
($E$, and $IP$ in $AP$ $1(\bmod{\,4})$, 
$1(\bmod{\,6})$ and $5(\bmod{\,8})$)

${}^*$\cite[1878]{lu} \'{E}. Lucas ($E$, Lucas sequences, and $IP$ 
in $AP$ $2(\bmod{\,5})$ and $7(\bmod{\,8})$)

\cite[1878/9]{k},  \cite[p. 4]{r2}, \cite{y} E. E. Kummer 
($E$ and Euclid's proof revisted with $p_1p_2\cdots p_n-1$ instead 
of $p_1p_2\cdots p_n+1$)

 \cite[1881]{pe}, \cite{y}  J. Perott ($CS$, $UFT$, $CM$, 
the fact that $\sum_{n=1}^{\infty}1/n^2<2$, the estimate of upper bound of 
number of integers $\le N$ by some square)

${}^*$\cite[1886]{ba}, \cite[p. 418]{d} A. S. Bang 
($E$, cyclotomic polynomials, and $IP$ in $AP$ 
$1(\bmod{\,k})$ with $k\ge 2$)

\cite[1887/8]{ge}, \cite[p. 413]{d} L. Gegenbauer ($CS$ and the 
the convergent series $\sum_{n=1}^{\infty}1/n^s$)

 \cite[1888]{pe2}, \cite[p. 414]{d} J. Perott (Theory of Commutative Groups)

\cite[1888]{sy3}, \cite[p. 7]{na} J.  J. Sylvester,
(evaluation of Euler's product $\prod_{p\le x}\left(1-1/p\right)^{-1}$ and  
the estimate $\sum_{n\le x}1/n\ge \log x$)

 \cite[1888]{sy3}, \cite[pp. 11--12]{na} J.  J. Sylvester 
($DS$, the series $\sum_{n=1}^{\infty}1/n$ is divergent and the series 
$\sum_{n=1}^{\infty}1/n$ is convergent)

${}^*$\cite[1888]{sy}, \cite[p. 418]{d} J. J. Sylvester 
($E$, cyclotomic polynomials, and $IP$ in $AP$ 
$1(\bmod{\,k})$ with $k\ge 2$)

${}^*$\cite[1888]{sy} J. J. Sylvester ($E$, and $IP$ in $AP$
$-1(\bmod{\,p^n})$ with $p$ any fixed prime)  

\cite[1889]{pe2}, \cite{pe3} J. Perott ($E$, 
Euclid's proof revisted, with  
$p_1p_2\cdots p_k-1$ instead of $p_1p_2\cdots p_k+1$)

\cite[1890]{ha2}, \cite[p. 414]{d} J. Hacks
(formula for the number of positive integers less than $N$ 
from  \cite[Ch. XI]{cah}) 

\cite[1890, p. 14]{st}, \cite[p. 414]{d}, \cite{r}, \cite{y} 
T. J. Stieltjes  ($E$ and the fact that the sum 
$p_1p_2\cdots p_k+p_{k+1}p_{k+2}\cdots p_{k+r}$  is not divisible 
by any  $p_i$ $(i=1,2,\ldots,k+r)$)

${}^*$\cite[1891]{lu2} \'{E}. Lucas ($E$, Lucas sequences,  $IP$ 
in $AP$ $1(\bmod{\,4})$, $5(\bmod{\,6})$ and in $AP$
$5(\bmod{\,8})$) 

\cite[1893]{bg}, \cite[p. 414]{d}, \cite{y}
C. O. Boije af Genn\"{a}s ($E$, $FT$ and 
the representation 
$Q=P/a-a>1$, where $a$ and $P/a$ are relatively prime factors of
$P:=p_1^{e_1}p_2^{e_2}\cdots p_n^{e_n}$) 

${}^*$\cite[1895]{wen}, \cite[p. 89]{na}  E. Wendt (the factorization 
$x^n-1=f(x)g(x)$, where $g(x)$ is  the least common multiple of polynomials 
$\{x^d-1:\, d\mid n\}$, common divisors 
of integers $f(x)$ and $g(x)$ with $x\in\Bbb Z$,
and $IP$ in $AP$  $1(\bmod{\,k})$) 

${}^*$\cite[1896]{ste} R. D. von Sterneck ($E$, 
 and $IP$ in $AP$ 
$-1(\bmod{\,k})$ with $k=2,3,\ldots$)

 \cite[1897]{th}, \cite[p. 9]{r2} A. Thue ($CM$ and $UFT$)

${}^*$\cite[1897]{va} K. Th. Vahlen (Gauss' periods of roots of unity, and
$IP$ in $AP$ $1(\bmod{\,k})$ with $k\ge 2$)

\cite[1899]{bra}, \cite[p. 414]{d} \cite[p. 3]{pol}, \cite{y} J. Braun 
($E$ and a prime divisor of $\sum_{i=1}^{k}(p_1p_2\cdots p_k)/p_i$)

\cite[1899]{bra}, \cite[p. 414]{d}, 
\cite{so1} J. Hacks (Euler's formula 
$\prod 1/(1-p^{-2})=\sum_{n=1}^{\infty}1/n^2=\pi^2/6$   
and the irrationality of $\pi^2$) 

${}^*$\cite[1899, p. 291]{lu3} \'{E}. Lucas (Lucas sequence, and 
$IP$ in $AP$ $1(\bmod{\,k})$ with $k\ge 2$)

\cite[1915]{bro} C. Hermite ($E$, a prime divisor of $n!+1$)

${}^*$\cite[1900, pp. 318--319]{cah} E. Cahen  
($E$, and $IP$ in $AP$ $1(\bmod{\,4})$, 
$1(\bmod{\,6})$ and $5(\bmod{\,8})$)  

${}^*$\cite[1903/04]{bv} G. D. Birkhoff and H. S. Vandiver 
(the existence of  primitive prime divisors 
of integers $a^n-b^n$, where $n\in\Bbb N$ and $a$ and $b$ 
are relatively prime integers)

${}^*$\cite[1905/6]{bau} M. Bauer  
($E$, and $IP$ in $AP$  $-1(\bmod{\,k})$ with $k\ge 2$)

\cite[1907]{bo}, \cite[p. 87]{uh} H. Bonse ($E$)

\cite[1909/10]{lev}, \cite[p. 414]{d} A. L\'{e}vy ($E$)

\cite[1911]{po}, \cite[p. 116, Theorem 114]{cl}, 
\cite[p. 419]{d}) H. C. Pocklington ($E$, and 
$IP$ which are not congruent to $1(\bmod{\,k})$)

${}^*$\cite[1912/13]{sch}, \cite[p. 91]{na} I. Schur  
($E$, law  of quadratic reciprocity, and $IP$ in $AP$ 
$2^{m-1}+1(\bmod{\, 2^m})$, $2^{m-1}-1(\bmod{\, 2^m})$ ($m\in\Bbb N $))

${}^*$\cite[1912/13]{sch}, \cite[p. 91]{na} I. Schur  
($E$, law  of quadratic reciprocity, and $IP$ in $AP$
$l(\bmod{\, k})$ for 
$k=8m$ (with $m$ being an odd positive squarefree integer) and 
$l=2m+1$, $l=4m+1$ or $l=6m+1$)

\cite[1912/13]{sch}, \cite[pp. 131, 324, Problem 108]{ps} 
I. Schur ($E$, and $IP$ of primes dividing the integer values of a
nonconstant  integer polynomial) 

${}^*$\cite[1913]{car} R. D. Carmichael ($IP$ in $AP$
$-1(\bmod{\,p^n})$ with $p$ any fixed odd prime, 
and $IP$ in $AP$ $-1(\bmod{\,3\cdot 2^n})$)

${}^*$\cite[1913]{has}  K. Hensel  
($E$, and $IP$ in $AP$ $1(\bmod{\,4})$, 
$1(\bmod{\,6})$ and $7(\bmod{\,8})$, $3(\bmod{\,8})$,
$9(\bmod{\,10})$ and $11(\bmod{\,12})$)

\cite[1915]{au}, \cite[p. 414]{d}, \cite[p. 11]{r2}  A. Auric 
($CM$, $FT$ and the estimate of number of positive integers 
$m=p_1^{e_1}p_2^{e_2}\cdots p_r^{e_r}$ less than $N$) 

\cite[1917]{met}, \cite[p. 415]{d} \cite[p. 11]{r2}  G. M\'etrod 
($E$ and a prime divisor of $\sum_{i=1}^nN/p_i$, where $N=p_1p_2\cdots p_n$)

\cite[1921]{poly},  \cite[pp. 131, 324,  Problem 107]{ps}
G. P\'{o}lya and G. Szeg\H{o} (Euler theorem, primes dividing the integer values of
the function $ab^x+c$ $(x=0,1,2,\ldots)$ with integral 
coefficients $a\not=0$, $c\not=0$ and $b\ge 2$)

\cite[1925, pp. 130, 322, Problem 108]{ps}, \cite[p. 14]{hw},
G. P\'{o}lya and G. Szeg\H{o} ($MPI$, Fermat numbers $F_n:=2^{2^n}+1$)

 \cite[1931]{hart} F. Hartmann 
($IP$ in $AP$ $\equiv 1(\bmod{\,p^n})$)) 

\cite[1934, p. 283]{er2}, \cite{er1} P. Erd\H{o}s ($C$ and de 
 Polignac's formula)

\cite[1934]{er2}, \cite{er1} P. Erd\H{o}s (de Polignac's formula
and inequalities for  central binomial coefficients)

${}^*$\cite[1937]{ba2}, \cite[p. 91]{na} A. S. Bang,
($E$, and $IP$ in $AP$ $2p^m+1(\bmod{\, 4p^m})$ with prime 
$p\equiv 3(\bmod{\, 4})$,
$2p^{2n+1}+1\bmod{\, 6p^{2n+1}})$ with prime $p\equiv 2(\bmod{\, 3})$,
and $4p^{2n}+1\bmod{\, 6p^{2n}})$ with prime $p\equiv 2(\bmod{\, 3})$)

 \cite[1938]{er3}, \cite[p. 17, Theorem 19]{hw}, 
\cite[pp. 5--6, Sixth proof]{az} P. Erd\H{o}s ($UFT$, $CM$ 
and the divergence of the sum $\sum_{p\,{\mathrm prime}}1/p$)

 \cite[1938,  $8^{\rm th}$ proof]{y} P. Erd\H{o}s (Chebyshev's argument, 
de Polignac's formula and $DS$) 

\cite[1938, pp. 16--17]{hw} G. H. Hardy and E. M. Wright ($FT$, 
a representation $n=m^2k$ where $k$ is squarefree and $CM$)

${}^*$\cite[1938, p. 13]{hw} G. H. Hardy and E. M. Wright ($E$, 
prime divisor of $a^2+b^2$, and $IP$ in $AP$ $5(\bmod{\,8})$)

\cite[1940, pp. 44--45]{shn} (published posthumously) L. G.  Schnirelman  
(the estimates $\lim_{x\to\infty}(\log_ax)^k/x=0$ for $a>1$ and $k>0$ and an 
enumerative argument)

\cite[1943]{b2} R. Bellman ($DS$ and the sum of prime reciprocals)

 \cite[1947]{b}, \cite[p. 7]{r2} R. Bellman ($MPI$ and a polynomial method) 

${}^*$\cite[1951]{nag} T. Nagell, 
($IP$ in $AP$ $-1(\bmod{\,k})$ for all $k\ge 2$)

\cite[1953]{tho} J. G. Thompson ($E$)

 \cite[1955]{f}, \cite[pp. 12--13]{r2}, \cite[p. 12]{pol}, \cite[p. 5]{az}
H. Furstenberg ($T$) 

\cite[1956]{dux} E. Dux ($DS$ and the sum of prime reciprocals)

\cite[1956]{har}, \cite[p. 6]{pol} V. C. Harris ($MPI$ and 
the numerators of finite continued fractions)

\cite[1957]{lm} J. Lambek and L. Moser  ($MPI$)

\cite[1958]{jar} D. Jarden 
(recurring sequences and $IP$ in $AP$ $1(\bmod{\,20})$)) 

\cite[1958]{mos} L. Moser ($DS$ and the sum of prime reciprocals)

 \cite[1959]{go} S. W. Golomb ($T$)

${}^*$\cite[1961]{rot} A. Rotkiewicz (Birkhoff-Vandiver theorem, 
the order of $k$ modulo a prime $p$, and $IP$ in $AP$ $1(\bmod{\,k})$)

${}^*$\cite[1962, pp. 60, 371--373]{scy} D. O. Shklarsky, N. N. Chentzov 
and I. M. Yaglom ($E$, divisors of $a^5-1$, 
Fermat little theorem, and $IP$ in $AP$ $1(\bmod{\,10})$)

\cite[1962/3]{es} T. Estermann  
(prime divisors of $nm$, where $n/m:=\prod_{d\mid k}(k^{k/d}-1)^{\mu(d)}$
with relatively primes integers $n$ and $m$, 
the order of $k$ modulo a prime $p$, and $IP$ in $AP$ 
$1(\bmod{\,10})$)
 
 \cite[1963]{go2} S. W. Golomb ($MPI$)

 \cite[1964]{e}, \cite[p. 7]{r2} A. F. W. Edwards ($MPI$)

 \cite[1964]{mu} A. A. Mullin ($E$)

\cite[1964, pp. 132--133]{si} W. Sierpi\'{n}ski (inequality for 
central binomial coefficient, mathematical induction, 
$UFT$, de Polignac's formula)

${}^*$\cite[1965]{bl} P. Bateman and M. E. Low ($E$, law  of quadratic 
reciprocity, and $IP$ in $AP$ $1(\bmod{\,24})$)

 \cite[1965]{ch} P. R. Chernoff ($CM$, $FT$, the estimate of upper bound
of number of $k$-tuples $(e_1,e_2,\ldots,e_k)$ satisfying 
$p_1^{e_1}p_2^{e_2}\cdots p_k^{e_k}\le N$)

\cite[1965]{wu}, \cite[p. 9]{na}, 
M. Wunderlich ($MPI$, Fibonacci sequence $(f_n)$, 
 the property $(m,n)=1$ implies $(f_m,f_n)=1$ and the factorization 
$f_{19}=113\cdot 37$)

\cite[1966]{hem} R. L. Hemminiger ($MPI$, a sequence $(a_n)$
with the property: $(m,n)=1$ implies $(f_m,f_n)=1$, 
the sequence $(a_n)$ defined recursively as 
$a_1=2$, $a_{n+1}=1+\prod_{i=1}^na_i$)

\cite[1966]{su} M. V. Subbarao  ($MPI$)

 \cite[1969]{coh} E. Cohen (de Polignac's formula and $DS$)

 \cite[1970, Problems 47 and 92]{si2} A. M\c{a}kowski
($E$ and relatively prime numbers)

 \cite[1970, Problem 50]{si2} A. Rotkiewicz ($MPI$ and Fibonacci numbers)

\cite[1970, Problem 52]{si2} W. Sierpi\'{n}ski 
(attributed to  P. Schorn by P. Ribenboim \cite[pp. 7--8]{r2})
($E$, $MPI$ and $AP$ $(m!)k+1$ for a fixed $k=1,2,\ldots,m$)

 \cite[1970, Problem 62]{si2} W. Sierpi\'{n}ski 
($E$, $MPI$ and $AP$) 

 \cite[1970, Problem 36]{si2} W. Sierpi\'{n}ski 
($MPI$ and triangular numbers) 

 \cite[1970, Problem 36]{si2} W. Sierpi\'{n}ski 
($MPI$ and tetrahedral numbers) 

\cite[1971]{djp} Problem 3 on IMO 1971 ($FT$ and Euler's theorem) 

 \cite[1974]{t} C. W. Trigg ($E$)

\cite[1976]{bar} C. W. Barnes ($E$, Theory of periodic continued fractions and 
Theory of negative Pell's equations $x^2-dy^2=-1$)

${}^*$\cite[1976]{np} I. Niven and B. Powell (the induction, the order of 
$k$ modulo a prime $p$, a polynomial equation, and $IP$ in $AP$ 
$1(\bmod{\,k})$)

\cite[1978, Theorem 1]{mo1}, \cite{mo2}, \cite[pp. 5--6]{pol} S. P. Mohanty 
($MPI$ and the induction)

\cite[1978, Theorem 2]{mo1}, S. P. Mohanty  ($MPI$ and Fermat little theorem)

\cite[1978, Theorem 3]{mo1}, S. P. Mohanty   ($MPI$ and prime divisors 
of Fibonacci numbers $f_p$)

${}^*$\cite[1978, p. 107]{sha}, \cite[pp. 178--179]{ers},
\cite[p. 209]{mur2},  D. Shanks  
(a prime divisor of $(2^{mp}-1)/(2^m-1)$ 
of the  form $p^nk+1$, and $IP$ in $AP$ $\equiv 1(\bmod{\,p^n})$)

\cite[1979]{ap} R. Ap\'{e}ry  (Euler's formula 
$\prod 1/(1-p^{-3})=\sum_{n=1}^{\infty}1/n^3:=\zeta(3)$   
and the irrationality of $\zeta(3)$) 

\cite[1979]{cha}, \cite[p. 118, Section 10.1.5]{cl}
 G. Chaitin (algorithmic information theory and an enumerative argument)

\cite[1979, p. 36]{wei},  A. Weil ($E$ and Group Theory)

\cite[1980]{van} C. Vanden Eynden ($DS$, the divergence of the 
series $\sum_{n=1}^{\infty}1/n$ and the convergence of the series     
$\sum_{n=1}^{\infty}1/n^2$)

 \cite[1980]{wa}, \cite[pp. 11--12]{r2}, \cite{chas} L. C. Washington 
(Theory of principal ideal domains, 
and the  factorizations $(1+\sqrt{-5})(1-\sqrt{-5})=2\times 3$ of 6 in the 
ring $\Bbb Z[a+b\sqrt{-5}]$)

${}^*$\cite[1981]{sm} R. A. Smith (Birkhoff-Vandiver idea, 
the solvability of the congruence $x^k\equiv 1(\bmod{\, p})$ 
with an integer of order $k$ modulo a prime $p$, and  $IP$  
in $AP$ $1(\bmod{\,k})$)

\cite[1981]{we} D. P. Wegener ($E$ and  primitive Pythagorean triples) 

\cite[1981]{wo} A. R. Woods (weak system of arithmetic $I\Delta_{0}$,  
$\Delta_0$-definable functions, the pigeonhole principle 
$PHP\Delta_0$ formulated for functions defined by $\Delta_0$-formulas)

 \cite[1984]{sr},  \cite{y} S. Srinivasan 
($MPI$, ``dynamical systems proof" and  the sequence 
$\left((2^{2^{m+1}}+2^{2^{m}}+1)/(2^{2^{m}}+2^{2^{m-1}}+1)\right)$)

 \cite[1984]{sr},  \cite{y} S. Srinivasan 
($MPI$, "dynamical systems proof", 
Fermat little theorem and  the sequence $\left((2^{p^{n+1}}-1)/(2^{p^n}-1)
\right)$)

\cite[1985]{od} R. W. K. Odoni ($E$, $MPI$ and a
 sequence $w_n$ recursively defined 
as $w_1=2$, $w_{n+1}=1+w_1\cdots w_n$ ($n\ge 1$)) 

\cite[1986]{ds}, \cite{sa} M. Deaconescu and J. S\'{a}ndor 
(divisibility property $n\mid \varphi(a^n-1)$, $a,n>1$)

\cite[1988]{pww}, \cite{pw} J. B. Paris, A. J. Wilkie and A. R. Woods  
(weak system of arithmetic $I\Delta_{0}$,  weak pigeonhole 
principle, $\Delta_0$-definable functions)

\cite[1993]{ru} M. Rubinstein ($CM$, $UFT$ and the asymptotic formula 
for the cardinality of a set $\{(e_1,\ldots,e_k)\in\Bbb N^k:\,
x_1\log p_1+x_2\log p_2+\cdots +x_k\log p_k\le\log x\}$)

${}^*$\cite[1994]{ro} N. Robbins ($MPI$, prime divisors of Fermat numbers, and
$IP$ in $AP$ $1(\bmod{\,4})$)

${}^*$\cite[1994]{ro} N. Robbins ($MPI$, prime divisors 
of Fibonacci numbers, and $IP$ in $AP$ $1(\bmod{\,4})$)

\cite[1995]{tr} D. Treiber ($DS$ and the sum of prime reciprocals)

\cite[1997, Problem 7.2.3]{ana} Problem on 
1997 Romanian IMO Team Selection Test, ($MPI$, the induction, Euler theorem
and a subsequence  of the sequence 
$(a^{n+1}+a^n+1)$ for a fixed integer $a>1$)

\cite[1998, Problem E3]{en} Problem of the training of the German
 IMO team, ($MPI$, the induction, the factorization 
$2^{2^{n+1}}+2^{2^{n}}+1=(2^{2^{n}}-2^{2^{n-1}}+1)(2^{2^{n}}+2^{2^{n-1}}+1)$
and $2^{2^{n+1}}+2^{2^{n}}+1$ has at least $n$ different prime factors
for each $n=0,1,2,\ldots$)

\cite[1998]{gold}, \cite[p. 16]{pol}  R. Goldblatt ($E$ and nonstandard 
Analysis)

${}^*$\cite[1998]{ss} N. Sedrakian and J. Steinig (a prime divisor 
of $(k^k-1)/[k^{k/p_1}-1,\ldots,k^{k/p_s}-1]$, 
where $p_1,\ldots,p_s$ are all distinct prime divisors of $k$ and 
$[a_1,...,a_s]$ denotes the greatest common divisor of $a_1,...,a_s$, 
and $IP$ in $AP$ $1(\bmod{\,k})$)

\cite[2000]{dal} M. Dalezman ($E$, $CM$)

\cite[2001, p. 4]{az} M. Aigner and G. M. Ziegler
($CM$, definite integral of the function $1/t$, $DS$, $UFT$)

\cite[2001, p. 3]{az}, \cite[p. 72]{aaf}, \cite{dp} 
M. Aigner and  G. M. Ziegler,  (Lagrange's theorem of Group Theory and Mersenne numbers)

\cite[2001]{por} \v{S}. Porubsky ($T$ and Theory of commutative rings)

\cite[2001/2, Problem 6]{pola}, \cite[p. 51, Problem 3.5.3]{no}
  Problem on Polish Mathematical Olympiad ($MPI$ and recursive sequence)

\cite[2002]{hi} M. D. Hirschorn ($CM$ and $FT$)

\cite[2003]{ab}, \cite[p. 6]{pol} J. M. Aldaz and A. Bravo  
($E$, a sequence $(P-2^n)$ with $P=\prod_{i=1}^rp_i$ and  $MPI$)

\cite[2003]{cw}, \cite{kl}  D. Cass and G. Wildenberg
($C$, periodic functions on integers)

 \cite[2003, p. 2]{ne1} C. W. Neville ($DS$)

\cite[2003]{sh} M. Somos and R. Haas ($MPI$)

 \cite[2004]{iiy}  T. Ishikawa, N. Ishida and Y. Yukimoto ($MPI$)

 \cite[2005, p. 35]{cp} R. Crandall and C. Pomerance ($DS$, the harmonic sum) 

${}^*$\cite[2005, pp. 92--64, Example 7.5.4]{mes}
M. R. Murty and J. Esmonde ($E$, properties of polynomial 
$f(x)=x^4-x^3+2x^2+x+1$, law  of quadratic reciprocity, and $IP$ 
in $AP$ $4(\bmod{\,15})$)  

${}^*$\cite[2005, p. 11]{mes} M. R. Murty and J. Esmonde  
(prime divisor of Fermat number 
$F_n:=2^{2^n}+1$ is of the form $2^{n+1}k+1$,
$F_n$ and $F_m$ are relatively prime if $m\not=n$, and  $IP$ in $AP$  
$\equiv 1(\bmod{\,2^n})$)

${}^*$\cite[2005]{n2} R. Neville (sequence $u_n=u_{n-1}+3u_{n-2}$, the induction,  
and $IP$ in $AP$ $1(\bmod{\,3})$)

${}^*$\cite[2005]{n2}, 

R. Neville (Lucas sequence $u_n=u_{n-1}+qu_{n-2}$
with a prime $q\ge 5$, the induction, Legendre symbol,  and
$IP$ in $AP$ $a(\bmod{\,20})$ for $a\in\{1,3,7,9\}$)

\cite[2006]{s} F. Saidak ($MPI$)

 \cite[2006]{ki} L. J. P. Kilford ($CS$)

 \cite[2007]{g} M. Gilchrist ($MPI$)

\cite[2007, pp. 110--111]{aaf} T. Andrescu, D. Andrica and Z. Feng
(first proof via induction; 
the second proof due by Sherry Gong via induction using Euler's theorem)

${}^*$\cite[2007, p. 4]{gr3} A. Granville ($E$, divisors of $a^2+a+1$, 
Fermat little theorem, and  $IP$ in $AP$ $1(\bmod{\,3})$)

\cite[2007, p. 2]{gr3} A. Granville ($E$, $FT$ and 
Chinese remainder theorem)

\cite[2008]{djp} Problem 3 on IMO 2008 (quadratic residues modulo 
a prime and infinitely many positive integers $n$ such that $n^2+1$ has a 
prime divisor greater than $2n+\sqrt{2n}$) 

E. Baronov \cite[2008, p. 12, Problem 5]{ibk} E. Baronov
($UFT$ and an enumerative argument)

E. Baronov \cite[2008, pp. 12--13, Problem 6]{ibk} E. Baronov
($UFT$ and an enumerative argument)

 \cite[2008]{jo} B. Joyal, (sieve of Eratosthenes and the formula for the 
proportion of the positive integers which are divisible by one 
of the first $n$ primes)

\cite[2008]{ng} P. Nguyen (weak theories of Bounded Arithmetic
"minimal" reasoning using concepts 
such as (the logarithm) of a binomial coefficient).

\cite[2008]{sci} A. Scimone ($E$ and $FT$)

\cite[2009, p. 4]{pol}) A. Granville ($E$ and Group Theory)

 \cite[2009]{m} I. D. Mercer ($C$)

\cite[2009]{p} J. P. Pinasco (Inclusion-Exclusion Principle,
 $CM$ and $DS$)

\cite[2009, p. 11]{pol} P. Pollack 
(the formula $\frac{5}{2}=\prod_{p}\frac{p^2+1}{p^2-1}$ and divisibility by 3)

 \cite[2010]{w}  J. P. Whang (de Polignac's formula)

\cite[2010]{coo} M. Coons ($UFT$, $CM$)

\cite[2011, p. 9]{aa} R. M. Abrarov and S. M. Abrarov 
($E$, M\"{o}bius function and delta function)

\cite[2011, p. 9]{aa} R. M. Abrarov and S. M. Abrarov 
(formula for the asymptotic density of primes and frequencies or 
probabilities)

\cite[2011, p. 9]{aa} R. M. Abrarov and S. M. Abrarov 
(prime detecting function, frequencies or probabilities)

 \cite[2011]{c}  R. Cooke (Theory of Finite Abelian Groups, 
the product of cyclic groups $\Bbb Z_{2n_1}\cdots\times Z_{2n_m}$ cannot be 
generated by fewer than $m$ elements and the isomorphism of the rings 
$\Bbb Z_{ab}$ and $\Bbb Z_{a}\times Z_{b}$) 

\cite[2011]{pol2} P. Pollack 
(uncertainty principle for the  M\"{o}bius transform
of arithmetic functions, entire function, pole of a rational function)

\cite[2011]{asp} J. M. Ash and T. K. Petersen ($MPI$ and $FT$)

\cite[2011]{mi} D. G. Mixon ($C$, $UFT$ and the pigeonhole principle)

\cite[2012]{me} R. Me\v{s}trovi\'{c} ($FT$ and 
representation of a rational number in a positive integer base)

\cite[2012]{me2} R. Me\v{s}trovi\'{c}
 (Meissel's identity $\sum_{n=1}^{\infty}\mu(n)\left[x/n\right]=1$ 
 and Pinasco's revisted proof)

\cite[2012]{me2} R. Me\v{s}trovi\'{c} (M\"{o}bius inversion formula and  
Legendre's formula 
$\pi(n)-\pi(\sqrt{n})=\sum_{d\mid \Delta}\mu(d)\left[x/d\right]-1$}

${}^*$\cite[2012]{me3} R. Me\v{s}trovi\'{c} ($MPI$, Euler's totient function, 
and $IP$ in $AP$ $1(\bmod{\,p})$ with a $p$ a prime)

[this article, Subsec. 2.5, 2012]  
(the formula $\prod_{n=1}^{\infty}(1-x^n)^{\mu(n)/n}=e^{-x}$ with 
$|x|<1$ and the irrationality of $e$)

[this article, Sec. 4, 2012] ($E$, $FT$ and the 
formula $(n^s-1)=(n-1)(\sum_{i=0}^{s-1}n^i)$)}


\cite[2012]{me4} R. Me\v{s}trovi\'{c} ($C$, $CN$ and $UFT$)

\cite[2015]{al} L. Alpoge ($C$ and $UFT$)

\cite[2015]{maj} B. Maji ($MPI$)

\cite[2015]{nor} S. Northshield (the estimates of the product 
$\prod_{p}\sin\left(\frac{\pi}{p}\right)$)

\cite[2016]{boo} A. R. Booker (analytic number theory, $UFT$ and $C$)

\cite[2016]{cl2} P. L. Clark, (the Euclidean criterion 
for irreducibles, $T$ and $UFT$)

\cite[2017]{gr5} A. Granville ($C$ and $UFT$)

\cite[2017]{sad} A. Sadhukhan ($C$)

\cite[2017]{sad} A. Sadhukhan ($C$)

\cite[2017]{sek} S.-I. Seki (valuation theory and approximation theorem)

\cite[2017]{sek}($DS$, Roth'stheorem and Euler-Legendre's theorem for 
arithmetic  progressions)

\cite[2017]{sek} ($DS$, $C$ and Euler-Legendre's theorem)

\cite[2017]{nor2} S. Northshield ($UFT$ and the idea of Furstenberg's proof)

\cite[2017]{nor2} S. Northshield ($UFT$ and the random integer)

\cite[2017]{me5} R. Me\v{s}trovi\'{c} ($E$ and $UFT$)

\vspace{5mm}

The following  Author${}^2$ and Subject Indices contain
names of all authors of references of this article
related to the proofs of $IP$, and mathematical concepts (notions) 
and notations that appear  in this article, respectively.
\footnote{${}^2$ The data about authors are owned from 
Wikipedia (List of mathematicians):
{\tt http://en.wikipedia.org/wiki/List\_of\_mathematicians}}

\vspace{3mm}

\centerline{{\bf D) Author Index}}

\vspace{3mm}

{\small
 
{\it {\bf A}bel, N. H.} (Norway, 1802--1829), 8, 18, 61;
{\it Abel, U.}, 27(2); 
{\it  Abrarov, R. M.} 18(7), {\bf 61(2), 58}; 
{\it Abrarov, S. M.,} 18(7), {\bf 61(3)};
{\it  A\v{g}arg\"{u}n, A. G.}, 6;   
{\it Aho, A. V.}, 14;  
{\it Aigner, M.} (Austria, born 1942), 8, 16, 19(2), 23, 24, {\bf 53, 
54, 56(2)};
{\it  Aldaz, J. M.}, 9, {\bf 56};
{\it  Alppoge, L.}
{\it Andji\'c, M.}  
  {\it Andrescu, T.} (Romania/USA, born 1956), 13, 16, 17, 31, {\bf 56(2), 57}; 
{\it Andrica, D.}, 13, 16, 17, 31, {\bf 56(2), 57}; 
{\it  Andrews, G. E.} (USA, born 1938), 11, 16;
{\it Ap\'{e}ry, R.} (France/Greece, 1916--1994), 22, 55;
{\it  Apostol, T. M.} (USA/Greece, born 1923), 22, 23;
{\it D'Aquiono, P.}
{\it Arana, A.}, 25;
{\it Arun-Kumar, S.} 
{\it Ash,  J. M.}, 14, {\bf 58}; 
{\it A\u{g}arg\"{u}n, A. G.}, 3, 3, 3; 
{\it Auric, A.} (France, 18??-19??), 24, {\bf 53}.

{\bf B}{\it aaz, S.}
{\it Bang, A. S.}, 31, 33, {\bf 52, 53};
{\it Barnes, C. W.}, 16, {\bf 55}; 
{\it Baronov, E.}, 28(2), {\bf 57(2)};
{\it Bateman, P. T.}, (USA, born 1919), 21, 33, 34, {\bf 54}; 
{\it Bauer, M.} (1874--1945), 32, {\bf 53}; 
{\it Bellman, R. E.} (USA, 1920--1984), 14,  19, 22, {\bf 54(2)}; 
{\it  Berndt, B. C.} (USA, born 1939), 23;
{\it Bernoulli, J.} (Switzerland, 1654--1705) 23; 
{\it Birkhoff, G. D.} (USA, 1884--1944), 32, {\bf 53};
{\it Bogomolny, A.}
{\it Boije af Genn\"{a}s, C. O.} (Sweden, 18??-19??), 9, {\bf 52}; 
{\it Bonse, H.} , 10, {\bf 53};
{\it Booker, A. R.}
{\it Borning, A.}, 6;
{\it Braun, J.}, 9(2), 22, {\bf 52(2)};
{\it Bravo, A.}, 9, {\bf 56}; 
{\it Brillhart, J.} (USA, ?), 27;
{\it Brocard, H.} (France, 1845--1922), 9, {\bf 53};
{\it Brown, M.} (USA, born 1931), 25;
{\it   Buck, R. C.} (USA, 1920--1998), 22;
{\it Buhler, J. P.}, 6;
{\it Burton, D. M.}, 15, {\bf 51}.

{\it {\bf C}ahen, E.}, 24, 34, {\bf 52, 53}; 
{\it  Caldwell, C. K.}, 6(3); 
{\it Carmichael, R. D.} (USA, 1879--1967), 31, {\bf 53}; 
{\it Cass, D.}, 25, {\bf 56};
{\it Chaitin, G.} (USA/Argentina, born 1947), 27(2), {\bf 55}; 
{\it Chastek, B.}, 17, {\bf 55}; 
{\it Chebyshev, P. L.} (Russia, 1821--1894), 21(2);
{\it Chentzov, N. N.}, 34, {\bf 54};
{\it Chernoff, P. R.}, 24, {\bf 54}; 
{\it  Choe, B. R.}, 22;
{\it Clark, P. L.}, 8, 17, 27, 30, 31, {\bf 53, 55};
{\it Clarkson, J. A.}, 19  ;
{\it Cohen, E.}, 20, {\bf 54};
{\it Conrad, K.}, 33(2);
{\it Cooke, R.}, 15, 17, 19, 22, {\bf 58};
{\it Coons, M.}, 24, {\bf 57};
{\it Cosgrave, J. B.} (Ireland, born 1946), 9;
{\it {\bf C}ox, C.D.}
{\it Crandall, R. E.}, 6, 20(2), 21, {\bf 56}; 
{\it Crstici, B.}, 15.

{\it {\bf D}alezman, M.}, 10, {\bf 56};
{\it Deaconescu, M.}, 16, {\bf 56}; 
{\it Dedekind, R.} (Germany, 1831--1916),  17, 59;
{\it Detlefsen M.}, 25; 
{\it  Diamond, H. G.}, 21(2), 27;
{\it Dickson, L. E.} (USA, 1874--1954), 8, 9(6), 14, 
15(2), 19(3), 22, 23, 24(2),
29(2), 30(4), 31(5), 33, 34(3), 51(7), 52(9), 53(4);
{\it Dilcher, K.}, 9;
{\it Dirichlet,  J. P. G. L.} (Germany, 1805--1859), 9;
{\it Djuki\'c, D.}
{\it Dunham, W.}, 9;
{\it Dux, E.}, 19, {\bf 54}.

{\it {\bf E}dwards,  A. W. F.} (Britain, born 1935), 12, {\bf 54}; 
{\it Elsholtz, C.}, 26(2);
{\it Engel, A.}, 13, {\bf 56};
{\it Eratosthenes} (Ancient Greece, circa 276--194 \small{B.C.}), 10, 59;
{\it Erd\H{o}s, P.} (Hungary, 1913--1996), 19, 21, 24, 31, 
 {\bf 53(3), 54, 55};
{\it Esmonde, J.}, 32(2), 35, {\bf 57};
{\it Estermann, T.} (?, 1902--1991), 32, {\bf 54};
{\it Euclid}, 3, 4, 5(2), 9, 35, {\bf 51}; 
{\it Euler, L.} (Switzerland, 1707--1783), 15(2), 18(2), 22, 23, 29, 
{\bf 51(3)};
{\it Everest, G.}, 11, 18, 26, {\bf 51(2)}.

{\it {\bf F}eng, F.}, 4;
{\it de Fermat,  P.} (Basque Country/France, 1601--1655), 8, 11(5), 
12(2), 14, 17, 20, 31, 32, 34, 51(4), 54, 55(2), 56, 57;
{\it Fibonacci $($Leonardo of Pisa$)$} (Italy, circa 1170--1250), 13(3),
14, 34, 54, 55, 56;
{\it  Fletcher, C. R.}
{\it  Forman, R.}, 26(2);
{\it Fourier, J. B. J.} (France, 1768--1830), 23;
{\it Furstenberg, H.} (USA/Israel, born 1935), 24, {\bf 54};
{\it Fuss, P.-N.} (Switzerland, 1755--1826), 11, {\bf 51}.

{\it {\bf G}allot, Y.}, 6(2);
{\it Galois, \'{E}.} (France, 1811--1832), 33;
{\it  Garrison, B.}, 26;
{\it Gauss, C. F.} (Germany, 1777--1855), 4, 6, 34, 36;
{\it Gegenbauer, L.} (Austria, 1849--1903), 19, {\bf 52};
{\it Genocchi, A.} (Italy, 1817--1889), 31, {\bf 51};
{\it  Gerst, I.}, 27; 
{\it  Giesy, D. P.}, 22;
{\it Gilchrist, M.}, 14, {\bf 57}; 
{\it Goldbach, C.} (Germany, 1690--1764), 7, 11(5), 12(2), 14, 51; 
{\it Goldblatt, R.} (New Zealand, ?), 29, {\bf 56};
{\it  Goldstein, L. J.}, 21;
{\it Golomb,  S. W.} (USA, born 1932), 14, 25, {\bf 54(2)};
{\it Gong, S.} (USA, ?), 16,
{\it Graham, R. L.} (USA, born 1935), 6, 16; 
{\it Granville, A.} (Britain, born 1962), 10, 16, 30(2), 33, 34(2), 
{\bf 57(2)};
{\it Gueron, S.}, 31; 
{\it Guy, R. K.} (England/Britain, born 1916), 10, 11, 22. 

{\it {\bf H}aas, R.}, 12, {\bf 56};
{\it Hacks, J.}, 24, {\bf 52};
{\it  Hadamard, J.} (France, 1865--1963), 18, 21;
{\it Hardy, G. H.} (England, 1877--1947), 3(2), 5, 8(2), 9, 11, 18, 
19, 20(2), 21, 22(4), 28, 34(4), 35, 36, {\bf 51(3), 53, 54(2)}; 
{\it  Harris, V. C.},  12, {\bf 54};
{\it Hartmann, F.}, 31, {\bf 53};
{\it Hasse, H.}  (Germany, 1898--1979), 13, 16, 19, 31, 32, 34(2), {\bf 51(2), 
53};
{\it Hausdorff, F.} (Germany, 1868--1942), 25; 
{\it Heath, T. L.} (Britain, 1861--1940), 3, 4, 5(2), 9, 35, {\bf 51};
  {\it Heaslet, M. A.}, 10, {\bf 53}; 
{\it Hemminiger, R. L.}, 13, {\bf 54};
{\it  Hermite, C.} (France, 1822--1901), 9, 53;
{\it Hetzl, A.}
{\it  Hildebrand, A. J.}, 28;
{\it Hirschorn, M. D.}, 24, {\bf 56}; 
{\it Honsberger, R.}, 19;
{\it de l'H\^{o}pital, G.} (France, 1661--1704), 21;
{\it Hurwitz, A.} (Germany, 1859--1919), 11.  

{\it {\bf I}rvine, S. A.}
{\it  Ishida, N.}, 14, {\bf 56}; 
{\it Ishikawa, T.}, 14, {\bf 56};
{\it Ivanov, A.}
{\it  Iwaniec, H.}, (Poland, born 1947) 31, {\bf 51}. 

{\it {\bf J}aensch, R.}, 19;
{\it Jankovi\'c. V.}
{\it  Jarden, D.}, 35, {\bf 54};
{\it  Jaroma, J. H.}, 12; 
{\it  Joyal, B.}, 10, {\bf 57};
{\it Joyce, D.}, 4.

{\it  {\bf K}eng, H. L.}, 31;
{\it  Khinchin, A. Y.} (Russia/Soviet Union, 1894--1959), 22; 
{\it Kilford, L. J. P.}, 23, {\bf 57};
{\it  Kirch, A. M.}, 25; 
{\it Klazar, M.}, 25, {\bf 56};
{\it Knopfmacher, K.}, 25,  43; 
{\it Knorr, W.} (USA, 1945--1997), 5;
{\it Knuth, D. E.} (USA, born 1938), 6, 16;
{\it Kolev, E.} 
{\it Korfhage, R. K.}, 11; 
{\it Kowalski, E.}, 31, {\bf 51};
{\it Kraus, L.}, 43;
{\it Kronecker, L.} (Germany, 1823--1891), 19, 31, 43;
{\it Kummer,  E. E.} (Germany, 1810--1893), 9, 15, 32, {\bf 52}.{
{\it Kurokawa, N.} 

{\it  {\bf L}agarias, J. C.} (USA, born 1949), 13;
{\it  Lagrange, J. L.} (France, 1736--1813) 16(2), 56;
{\it Lambek, J.} (Germany/Canada, born 1922), 12; {\bf 54};
{\it  Landry, F.}, 30, {\bf 51};
{\it  Lebesgue, V. A.} (France, 1875--1941), 31(2), 34(3), {\bf 51(4)};
{\it Leitsch}
{\it Legendre, A.-M}. (France, 1752--1833), 20, 22, 29; 
{\it L\'{e}vy, A.}, 9, {\bf 53};
{\it  Lima, F. M. S.}, 22;
{\it Loomis, S. E.}
{\it Lord, N.}, 14;
{\it Low, M. E.}, 33, 34, {\bf 54};
{\it Lubotzky, A.} (Israel, born 1956), 26;
{\it Lucas, \'{E}.} (France, 1842--1891), 32, 34(2), {\bf 52(2), 53};
{\it Luzin, N. N.} (Soviet Union/Russia, 1883--1950), 22(2); 
{\it Lyusternik, L. A.} (Soviet Union, 1899--1981), 22.

{\it {\bf M} aji}
{\it M\c{a}kowski, A.}, 10, {\bf 55};
{\it Mamangakis, S. E.}, 11;
{\it Mati\'c, I.} 
{\it Mazur, B.}, (USA, born 1937), 5, 9;
{\it Meissel, D. F. E.} (Germany, 1826--1895), 15;
{\it Mercer, I. D.}, 25;
{\it Mersenne, M.} (France, 1588--1648), 11, 12, 16(2), 56; 
{\it  Mertens, F.} (Germany, 1840--1927), 18, 20(2), {\bf 51};
{\it Me\v{s}trovi\'{c}, R.}, 10, 15(2), 31, {\bf 58(6)};
{\it M\'etrod, G.} (18??-19??), 9, {\bf 53};
{\it Mixon, D. G.}, 24, {\bf 58}; 
{\it M\"{o}bius, A. F.} (Germany, 1790--1868), 17(4), 18, 22, 31, 58(2);
{\it  Mohanty, S. P.}, 14(3), 16, {\bf 55(3)};
{\it Moll, V. H.} 
{\it  Morton, H. R.}, 14;
{\it Morton, P.}, 26(2);
{\it Moser, L.} (Canada, 1921--1970), 12, 19, {\bf 54(2)};
{\it Mullin, A. A.}, 9, 11, {\bf 54};
{\it  Murty, M. R.}, 30(2), 31, 32(2), 33, 35, {\bf 57(2)}.

{\it {\bf  N}agell, T.} (Norway, 1895--1988), 33, {\bf 54};
{\it Narkiewicz, W.}, 8, 9(2), 13, 15, 19(2), 23, 30(2), 31, 
32(4), 33(2), 34, {\bf 51, 52(3), 53(2), 54};
{\it  Nathanson, M. B.} (USA, ?), 30;
{\it Naur, T.}, 11;
{\it Neville,  C. W.}, 19, 25, {\bf 56}; 
{\it  Neville, R.}, 35(2), {\bf 57(2)};
{\it  Newman, D. J.} (USA, 1930--2007), 21;
{\it Nguyen, P.}, 29, {\bf 57};
{\it Northshield, S.}
{\it Niven, I. M.} (Canada/USA, 1915--1999), 32, {\bf 55};
{\it Nowakowski, R.}, 10, 11;
{\it Nowicki, A.}, 12, 13(2), 13, {\bf 56}.

{\it  {\bf O}doni, R. W. K.}, 10, {\bf 56}.
{\it Osler, T. J.}

{\it {\bf P}apadimitriou, I.}, 23;
{\it Paris,   J. B.} (Britain, born 1944), 28, 29, {\bf 56(2)};
{\it Patashnik, O.}, 6, 16;
{\it Pell, J.}  (Britain, 1611--1685), 3, 16, 55; 
{\it Penk, M. A.}, 6;
{\it Perott,  J.} (France, 1854--1924), 9, 10, 15, 16, 23, {\bf 52(4)};
{\it Petersen,  T. K.}, 14, {\bf 58};
{\it Petrovi\'c, N.} 
{\it Pinasco, J. P.}, 15, {\bf 57};
{\it  Plouffle, S.} (Canada, born 1956), 6;
{\it  Pocklington, H. C.} (England, 1870--1952), 31, {\bf 53};
{\it de Polignac, A.} (France, 1817--1890), 20(2), 24(2), 54(2);
{\it Pollack, P.}, 8, 9(2), 11, 12, 13(2), 14, 15, 16(2), 17, 18, 19,
20, 21, 22(3), 23, 24, 26(3), 28, 29, 30, 33, {\bf 51(2), 
52, 54(2), 55, 56(2), 58};    
{\it  P\'{o}lya, G.} (Hungary, 1887--1985), 11, 13, 15(2), 25, 
{\bf 53(3)};
{\it Pomerance, C.} (USA, born 1944), 20(2), 21, {\bf 56}; 
{\it van der Poorten, A.} (Netherlands/Australia, 1942--2010), 26; 
{\it Porubsky, \v{S}.},  25(3),  43, 46, {\bf 56};
{\it P\'{o}sa, L.} (Hungary, born 1947), 10; 
{\it Powell, B.}, 32, {\bf 55};
{\it Pythagoras} (Ancient Greece, circa 585--501 {\small B.C.}), 3, 4, 17, 55. 

{\it {\bf R}eich, S.}, 11;
{\it  Reddy, K. N.}, 12;
{\it Ribenboim, P.} (Brazil/Canada, born 1928), 4, 6, 7, 8(2), 9(3), 10, 
11, 12(2), 14, 16, 17, 18, 22, 23(2), 24(2), 35, {\bf 51(2), 52(3), 53(2),
54(3), 55(2)};
{\it Richter, C.}
{\it Riemann, B.} (Germany, 1826--1866), 21;
{\it  Robbins, N.}, 34, {\bf 56(2)};
{\it Rotkiewicz, A.}, 13, 32, {\bf 54, 55};
{\it Rubinstein, M.}, 24, {\bf 56}. 

{\it {\bf S}aidak, F.}, 14, {\bf 57};
{\it  \v{S}al\'{a}t, T.} (Slovakia, 1926--2005), 19;
{\it Samuel, P.} (France, 1921--2009), 17; 
{\it S\'{a}ndor, J.}, 10, 11, 15, 16(2), {\bf 56(2)};
{\it Satoh, T.}
{\it Schnirelman, L. G.} (Soviet Union, 1905--1938), 27, {\bf 54};
{\it Schorn, P.}, 5;
{\it  Schur, I.} (Germany, 1875--1941), 25, 33(2), {\bf 53(3)}; 
{\it Scimone, A.}, 10, {\bf 57}; 
{\it Seki, S.-I.}
{\it Sedrakian, N.}, 32, {\bf 56};
  {\it Selberg, A.} (USA/Norway, 1917--2007), 30(2);
{\it Serret, J. A.} (France, 1819--1885), 34(2), {\bf 51(2)};
{\it Shanks, D.} (USA, 1917--1996), 31(2), {\bf 55};  
{\it Shapiro, H. N.}, (USA, born 1928), 15, 28, 30(2);
{\it Shparlinski, I.} (Australia, born 1956), 26;
{\it Shklarsky, D. O.}, 34, {\bf 54};
{\it  Siebeck, H.}, 13;
{\it  Siebert, H.} 27(2); 
{\it Sierpi\'{n}ski, W.} (Poland, 1882--1969),  10(2), 12, 13, 14(2), 
24(2), 26, {\bf 54, 55(4)};
{\it Silverman, J. H.}
{\it Sloane, N. J. A.} (USA, born 1939), 6, 14, 50;
{\it Smith, R. A.}, 32, {\bf 55};
{\it Somos, M.}, 12, {\bf 56}; 
{\it Sondow, J.}, 11, 22, {\bf 52};
{\it Spohr, H.} 
{\it Srinivasan, S.}, 17, 32, {\bf 55(2)};
{\it Steinig, J.}, 32, {\bf 56};
{\it  Stephens, P. J.}, 13;
{\it  von Sterneck, R. D.}, 31, {\bf 52};
{\it Stieltjes, T. J.} (Netherlands/France, 1856--1894), 9, 13, {\bf 52};
{\it  Subbarao, M. V.} (India, 1921--2006), 12, {\bf 54};
 {\it Sun, Z.-W.} (People's Republic of China, born 1965), 6(2);
{\it Sur\'{a}nyi, J.}, 31, {\bf 55};
{\it Sylvester, J. J.} (USA/Britain, 1814--1897), 18, 19(2), 27, 31(2), 
{\bf 51, 52(4)}; 
{\it Szeg\H{o}, G.} (Hungary/USA, 1895--1985), 11, 15, 25, {\bf 53(2)}.

{\it  {\bf T}attersall, J. J.}, 29(2);
{\it Templer, M.}, 6;
{\it Tessler, R.}, 31;
{\it Tikekar, V. G.}
{\it Thompson, J. G.} (USA, born 1932), 9, {\bf 54};
{\it Thain, N.}, 30(2), 33;
{\it Thue, A.} (Norway, 1863--1922), 23, {\bf 52}; 
{\it Treiber, D.}, 19, {\bf 56};
{\it Trevi$\stackrel{\sim}{\rm n}$o, E.}
{\it Trigg, C. W.}, 9, {\bf 55}.
 
{\it {\bf U}spensky, J. V.} (Russia, 1883--1947), 10, {\bf 53}.

{\it  {\bf V}ahlen, K. Th.} (Germany, 1869--1945), 32, {\bf 52};
{\it  de la  Vall\'{e}e-Poussin, C. J.} (Belgium, 1866--1962), 18, 21;
{\it {\bf V}an der Poorten, A. J.}  
{\it   Vanden Eynden, C.}, 19, {\bf 55};
{\it Vandiver, H. S.} (USA, 1882--1973), 32, {\bf 53};
{\it  Varadarajan, V. S.}  (India/USA, born 1937), 18;
{\it  Vardi, I.}, 6, 12;
{\it Vorob'ev, N. N.}, 13.

{\it van der {\bf W}aerden, B. L.} (Netherlands, 1903--1996)
{\it  Wagstaff, Jr.,  S. S.} (USA, 1944), 22;
{\it Ward, M.} (USA, 1901--1963), 11, 13(2), 18, 26, {\bf 51(2)};
{\it Ward, T.}, 26;
{\it Washington, L. C.}, 17, 32, {\bf 55};
{\it Wegener, D. P.}, 17, {\bf 55};
{\it Weil, A.} (France, 1906--1998), 3, 16, {\bf 55}; 
{\it Wendt, E.}, 32, {\bf 52};
{\it Whang, J. P.}, 24, {\bf 57};
{\it Wildenberg, G.}, 25, {\bf 56}; 
{\it Wilkie, A. J.} (England, born 1948), 28(2), 29, {\bf 56(2)};
{\it Woods, A. R.},  28(2), {\bf 55, 56};
{\it Wright, E. M.} (England, 1906--2005), 5, 8(2), 9, 11, 18, 
19, 20(2), 21, 22(4), 28, 34(4), 35, 36, {\bf 51(3), 53, 54(2)};
{\it  Wunderlich, M.}, 13(2), {\bf 54}.

{\it {\bf Y}aglom, A. M.} (Soviet Union, 1921--2007), 16, 20(3), {\bf 51}; 
{\it  Yaglom, I. M.} (Soviet Union, 1921--1988), 16, 20(3), 34, 
{\bf 51, 54};
{\it  Yamada, T.}, 8, 9(2), 12, 17(2), 20, {\bf 52(5), 53, 55(2)};
{\it Yoo, J.}, 32;
{\it Yukimoto, Y.}, 14, {\bf 56}.

{\bf Z}{\it agier, D.} (USA, born 1951), 21;
{\it Zhang, S}., 3, 35.
{\it Ziegler, G. M.} (Germany, born 1963), 8, 16, 19(2), 23, 24, {\bf 53, 
54, 56(2)}.}


\vspace{3mm}

\centerline{{\bf E)  Subject Index}}

\vspace{3mm}

{\small
 
{\bf A}belian group, 8;
abstract, 25;
Abstract Ideal Theory, 26;
additive, 20;
additive structure, 20;
additive structure of the integers, 20;
Algebra, 17;
algebraic, 4; 
algebraic argument, 17; 
algebraic integer, 17;
algebraic number, 17;
algebraic number theory, 15;
algebraic modification, 8;
 algebraic number, 17;
 Algebraic Number Theory, 8;
algebraic number theory argument, 8;
algorithm, 10;
algorithmic, 28;
algorithmic entropy of a positive integer $(H(N))$, 28;
``Algorithmic Information Theory", 28;
almost-injective integer sequence, 26; 
alternate proof, 9;
 Analysis, 19;
analytic, 8;
 Analytic Number Theory, 19;
analytic proof, 8;
Ancient Greek mathematicians, 3;
approximant, 12;
approximation, 22;
argument, 23; 
arithmetic, 3;
arithmetic function, 8;
arithmetic progression $(AP)$, 4;
arithmetic property, 26;  
arithmetical, 27;
asymptotic, 18;
asymptotic behavior, 20; 
asymptotic density of prime numbers, 18;
asymptotic formula, 59; 
asymptotically, 22;
asymptotically equivalent, 22;
axiom, 6;
axiom scheme, 29.     

{\bf B}ehavior, 20;
Bernoulli number $(B_n)$, 23;
better approximation, 22;
binomial coefficient $({n\choose k})$, 29;
Bonse's inequality, 11;
Book VII (of ``Elements"), 3;
Book VIII (of ``Elements"), 3;
Book IX (of ``Elements"), 3;  
bound, 21; 
bounded, 21;
Bounded Arithmetic, 29;
bounded formulas, 29;
boundedness, 55.

{\bf C}ardinality of a set $S$ $(|S|)$, 21;
central binomial coefficient,  57;
character, 32;
Chebyshev's argument, 20;
Chebyshev inequalities, 27;
Chinese remainder theorem, 10;
class, 14;
closed, 25;
closed set, 25;
coefficient, 16;
collection, 15;
combination, 32;
combinatorial, 8;
combinatorial argument, 24;
combinatorial modification, 8;
combinatorial method, 54;
combinatorial proof, 23;
combinatorial version, 26;
common divisor, 13;
common factor, 34;
common multiple, 4;
common prime factor, 12;
commutative, 16;
commutative algebra, 17;
commutative group, 16;
commutative ring, 26;
complement, 25;
Complex Analysis, 22;
complex characters mod $k$, 32;
complex function, 21;
complex variable, 22;
composite, 5;
composite number, 5;
concept, 3;
congruent, 26;
conjecture, 5;
connected, 25;
connected set, 25;
connected topological space, 25;  
connected topology, 25;
consecutive, 6;
consecutive primes, 6;
consequence, 18;
constant, 19;
construction, 12;
continued fraction, 12;
contradiction, 4;
convergence, 19; 
convergence of the series, 19; 
convergent, 24;
convergent series, 23;
coprime, 12;
coprime residue class, 37;
counting argument, 23;
counting method, 8; 
cyclic, 18;
cyclic group $(\Bbb Z_m)$, 18;
cyclotomic field $(\Bbb Q(\zeta_k))$, 34;
cyclotomic polynomial, 26.

{\bf D}edekind domains, 17;
delta function $\delta(x)$, 18;
definite, 21;
definite integral, 21;
degree of a polynomial, 27;
delta function $(\delta (x))$, 18;
denominator, 23;
density of primes, 21;
density of the set, 25;
Dickson's lemma, 15;

Dirichlet $L$-function, 32;
Dirichlet $L$-series, 32;
Dirichlet series, 22;
Dirichlet's theorem, 32;
Disquisitiones Arithmeticae, 6;
distribution of the primes, 21;
divergence of the sum, 18; 
divergent, 19;
divergent infinite series, 19; 
divergent series, 19; 
divergent sum, 19;
divisibility, 3;
divisibility property, 8; 
divisible, 4;  
divisor, 5;
domain, 5;
dynamic, 4;
dynamical, 18;
dynamical system, 18;
dynamical systems proof, 18.

{\it {\bf e}} (the constant), 20; 
elementary, 8;
elementary argument, 34;
elementary number theory, 12; 
elementary proof, 8; 
``Elements'' (of Euclid), 3; 
entire function, 61;
entropy, 28;
enumerating arguments, 24;
equality, 20;
equation, 17;
equivalent, 22; 
estimate, 10;
Euclidean Criterion, 31;
Euclidean domain $\Bbb F_p[x]$, 5;
Euclid's argument, 6;
Euclid's first theorem, 38; 
Euclid's idea, 8; 
Euclid's Lemma, 38; 
Euclid's method, 33; 
Euclid's number $(E_n)$, 6; 
Euclid's proof (of $IP$), 4;  
Euclid's second theorem, 9; 
Euclid-Mullin graph, 31;
Euclid-Mullin sequence $((E_n))$, 5;
Euclid sequence, 31; 
Euclid's theorem, 3;
Euclid's theory of numbers, 6;
Euler-Legendre's theorem, 31;
Euler-Mascheroni constant  $(\gamma$), 20;
Euler theorem, 16; 
 Euler's factorization, 22; 
 Euler's formula for $\zeta(2)$, 11; 
 Euler's formula for $\zeta(4)$, 22; 
 Euler's formula for $\zeta(2n)(n=1,2,\ldots)$, 23; 
 Euler's product, 19; 
Euler's product  for the Riemann zeta function,  19;
Euler's product formula, 31; 
Euler's proof of $IP$, 19; 
Euler's second proof of $IP$, 15; 
  Euler's totient function ($\varphi(n)$), 15;
evaluation, 23;
 exponent, 20;
 exponent of prime, 20;
exponential function, 16;

{\bf F}actor, 5; 
 factorial $(n!)$, 20;
factorization, 5;
factorization theorem, 5;
 Fermat little theorem, 8;
Fermat numbers, 11; 
Fibonacci's sequence ($(f_n)$), 13;
field, 6;
finite, 8;
finite  Abelian groups, 8;
finite continued fraction, 58;
finite extension of $\Bbb Q$, 34;
finite  group, 8;
finite linear combination, 29;
finite set, 30; 
finite support, 18;
finite union, 25;
formal identity, 20;
formula, 8;
fraction, 8;
fundamental theorem of arithmetic, 6;
function, 8; 
function of a complex variable, 22;
functions defined by $\Delta_0$-formulas, 29;
fundamental theorem of unique factorization of positive integers, 24; 
Furstenberg's ideas, 26;
Furstenberg's proof, 25;
Furstenberg's topological proof, 25.

{\bf G}alois Theory, 36;
Gauss' periods of roots of unity, 35;
 generalization of  Fermat numbers, 14;
 generalization of  Sylvester's sequence, 12;
generalized Euclid sequence, 30;
generated, 61;
geometric, 4;
geometrical, 20;
 geometrical interpretation of definite integral, 21; 
geometrical proof, 20;
Goldbach's idea, 12; 
greatest common divisor, 13;
greatest integer function $([x])$, 16; 
 
{\bf H}adamard-de la Vall\'{e}e-Poussin constant $(M=0.261497\ldots)$, 19; 
harmonic series, 19;
harmonic sum, 60;
Hausdorff topology, 25;
l'H\^{o}pital's rule, 22.
  
{\bf I}deal, 17;
identity, 16;
 Inclusion-Exclusion Principle, 15;
increasing infinite  sequence, 14;
indirect proof, 5;
induction, 11; 
induction axioms, 29;
inductively, 6; 
inequality, 11;  
infinite, 3;
infinite collection, 15;
infinite continued fraction, 12;
 infinite coprime sequences, 15; 
 infinite product, 198; 
infinite sequence, 12; 
infinite series, 19 
infinite subsequence, 13; 
infinitely, 3; 
infinitely many primes, 3;
 infinitude, 3;
 infinitude of primes ($IP$), 3;
infinity, 4;
injective, 26;
integer, 3;
integer argument, 27; 
integer coefficient, 16;
integer constant, 26;
integer function, 27;
integer polynomial, 26;
integer square, 30;
integer-valued function, 34;
integer sequence, 9;
integer value,  15;
integral, 22; 
International Mathematical Olympiad (IMO 1971 and 2008), 17, 28 {\bf 60, 61}; 
irrational number, 23;
irrational part, 33;
irrationality, 4;
irrationality  measure $(\mu)$, 11; 
irrationality of $\sqrt{2}$, 4;
irrationality  of $\pi^2$, 22; 
irrationality of $e$, 22;
irreducible, 6;
irreducible factor, 6;
irreducible polynomial, 6;
isomorphism, 58;
isomorphism of the rings, 61.
 
${\mathbf k}$th cyclotomic field $(\Bbb Q(\zeta_k))$, 34; 
 Kummer's number, 6;
$k$-tuple of, 38.

{\bf L}agrange theorem, 17;
law of quadratic reciprocity, 4;
leading coefficient, 27;
least common multiple, 34;
least common multiple of polynomials, 34;
Legendre's formula, 16;
Legendre symbol ($\left(\frac{\cdot}{p}\right))$, 37;
leg, 17;
length, 3;
$L$-function, 32;
line segment, 3;
linear, 14;
linear combination, 32;
linear recurrence, 14; 
linear second order recurrence, 37;
locally connected topological space, 25;
logarithm to the base $e$ $(\log x)$, 27;
logarithm of a binomial coefficient, 29;
logarithmic integral $(\mathrm{Li}(x))$, 22; 
logarithmic complex function, 21;
lower bound, 28;
Lucas sequence $(u_n)$, 35.

{\bf M}athematical induction, 58;
measured, 3;
measuring, 3;
 Meissel-Mertens constant ($M=0.261497\ldots$), 19;
 Mersenne number, 12;
 Mertens' first theorem, 21;
 Mertens' second theorem, 19;
 Mertens' third theorem, 20;
metrizable topology, 25;
 M\"{o}bius function $(\mu(n))$, 18;
 M\"{o}bius inversion formula, 18; 
 M\"{o}bius pair, 18; 
 M\"{o}bius transform, 18; 
monomial, 15; 
monotonic sequence, 11;
Moscow school of mathematics, 22; 
multiple, 4;
multiplicative, 26;
multiplicative group $(\mathbf Z/m\mathbf Z^*)$, 16; 
multiplicative group modulo a prime, 16;
multiplicative structure, 26; 
multiplicativity, 15;
mutually prime, 11;
mutually prime integers, 11.

{\bf N}atural number, 21; 
negative Pell's equation, 17;
nonconstant, 26;
nonconstant polynomial, 27;
nonnegative integer, 27; 
non-principal character, 32;
nonstandard Analysis, 29;
nonunits factor, 31;
non-zero polynomial, 35;
normal topology, 25;
not regular topology, 25; 
Number Theory, 9;
numerator, 12;

{\bf O}dd, 11;
odd prime, 16;
open set, 25;
order of $a(\bmod{p})$, 16; 
order of subgroup, 17; 

{\bf P}airwise relatively prime, 10;
partition, 31;
partition of the positive integer, 31;
Pell's equation, 17;
period, 35;
periodic continued fraction, 17;
periodic functions on integers, 25;
$\pi$ (the constant), 11; 
 pigeonhole principle, 24;
pigeonhole principle for functions defined
by $\Delta_0$-formulas ($PHP\Delta_0$), 29; 
pole of a rational function, 61; 
 de Polignac's formula, 20;
Polish Mathematical Olympiad (2001/02), 14, {\bf 56};
polynomial, 6;
polynomial growth, 26;
polynomial in two variables, 55;
polynomial method,  14;
polynomial over the field $F$, 6;
 positive integer, 5;
 positive integer base, 10;
positive constant, 22;
primality, 6;
prime, 3;
 prime-counting function ($\pi (x)$), 16;
prime detecting function, 18;
prime divisor, 9;
prime factor, 6;
 prime number, 3;
Prime Number Theorem, 21;
prime power, 34;
prime value, 27; 
primitive $k$th root of unity $(\zeta_k)$, 34;
primitive divisor, 27;
primitive prime divisor, 27;
primitive Pythagorean triples, 17;
primorial number, 6;
primorial prime, 7;
principal character, 32;
principal ideal domains, 18;
 probability, 16;
 probability theory, 31;
product of cyclic groups, 18;
proof by contradiction (reduction ad absurdum), 4;
 proper   subgroup, 16; 
Proposition 20 (of ``Elements"), 3;
Proposition 30 (of ``Elements"), 5;
Proposition 31 (of ``Elements"), 5;
Pythagorean school (at Croton), 3;
Pythagorean theorem, 4; 
Pythagorean triples, 17.

{\bf Q}uadratic residue, 36.

{\bf R}andom integer, 31;
randomly, 16;
range, 6;
rational function, 18; 
rational multiple, 23; 
rational part, 33;
rational prime, 34;
reciprocal, 18;
recurring sequence, 58;
 recursive sequence, 34;
 recursively defined sequence, 10; 
reduction ad absurdum, 4;
regular infinite continued fraction, 12;
relatively prime, 10; 
relatively prime positive integers, 10; 
 representation as a product of primes, 5;
representation of a rational number in a positive integer base, 10;
residue, 36;
residue class, 37;
Riemann zeta function ($\zeta(s)$), 8;
ring, 17;
ring of algebraic integers, 17;
ring of polynomials with integer coefficients $(\Bbb Z[T])$, 35;
Romanian IMO Team Selection Test (1997), 13, {\bf 60};
Roth's theorem, 31;
roots of unity, 35.

{\bf S}econd order recurrence, 37;
seed, 30; 
sequence $(a_n)$, 5;
series, 19;
set of    integers ($\Bbb Z$), 15;
set of natural numbers (positive integers) ($\Bbb N$), 25;
set of nonnegative integers $(\Bbb N_0)$, 27;
${}^*$--set of positive integers, 15;
set of primes, 19;
set of rational numbers $(\Bbb Q)$, 10;
sieve of  Eratosthenes, 10;
Sloane's On-Line Encyclopedia of Integer, 9; 
Sloane's sequence, 5;
sophisticated proof, 19;
square of a prime, 18;
squarefree, 28;
squarefree integer, 36;
subadivity, 28;
subadivity of algorithmic entropy, 28;
subexponential growth, 26;
subgroup, 16;
subsequence, 13;
successive prime numbers, 19;
Sylvester's sequence, 12;
Sylvester's version of the Chebyshev inequalities, 27;

{\bf T}etrahedral number, 14;
Theory $I\Delta_{0}$, 29;
Theory of commutative groups, 16;
Theory of commutative rings,  60;
Theory of algebraic numbers, 17;
Theory of   Dedekind domains, 17;
Theory of finite   Abelian groups, 8;
Theory of  negative Pell's equations, 17; 
Theory of periodic continued fractions, 17;
Theory of principal ideal domains, 17; 
Theory of profinite groups, 26;
Theory of  unique factorization domains, 17;
topological, 8;
 topological ideas, 25;
 topological method, 54;
 topological proof of $IP$, 25;
 topological space, 25;
 topology, 25; 
totient function $(\varphi(n))$, 15;
training of the German IMO team, 13, {\bf 60};
triangular numbers, 14;

{\bf U}ncertainty principle for the M\"{o}bius transform, 18;
unique factorization, 6; 
unique factorization domains, 17;
unique factorization theorem, 6;
uniqueness theorem for Dirichlet series, 22;
upper bound, 55.

{\bf V}alue of a polynomial, 26;
van der Waerden's theorem, 29;
vocabulary $0,1,+,\cdot ,<$, 29.

{\bf W}eak pigeonhole principle for $\Delta_0$-definable functions, 29; 
 weak system of arithmetic $(I\Delta_{0})$, 29; 
weak theories of Bounded Arithmetic, 29;
Wikipedia, 12.
  \vfill\eject

\noindent{\scriptsize ROMEO ME\v{S}TROVI\'{C}

\noindent  MARITIME FACULTY KOTOR,

\noindent UNIVERSITY OF MONTENEGRO

\noindent DOBROTA, 85330 KOTOR, MONTENEGRO}

\noindent{\it E-mail address}: {\tt romeo@ucg.ac.me}

\end{document}